\documentclass[brochure,10pt]{article}
\usepackage{amssymb,amsfonts,amsmath,footnote}
\usepackage{latexsym}
\usepackage{amsmath}
\usepackage{amsthm}
\usepackage{amsfonts}
\usepackage{amssymb}
\usepackage{amsmath,amscd}
\usepackage{graphicx}
\usepackage{graphics}
\usepackage{etex}
\usepackage[matrix,arrow]{xy}
\usepackage{color}
\usepackage{tikz-cd}
\newtheorem{Th}{Theorem}
\newtheorem{Prop}{Proposition}
\newtheorem{Co}{Corollary}
\newtheorem{Lm}{Lemma}
\newtheorem{Lma}{Lemma}[section]

\newtheorem{Rm}{Remark}
\newtheorem{Con}{Open Problem}
\newcommand{\be}{\begin{equation}}
\newcommand{\ee}{\end{equation}}
\newcommand{\bes}{\begin{equation*}}
\newcommand{\ees}{\end{equation*}}

\newcommand{\R}{\mathbb{R}}
\newcommand{\N}{\mathbb{N}}
\newcommand{\C}{\mathbb{C}}
\newcommand{\Z}{\mathbb{Z}}

\newcommand\res{\mathop{\hbox{\vrule height 7pt width .5pt depth 0pt
\vrule height .5pt width 6pt depth 0pt}}\nolimits}

\def\theequation{\thesection.\arabic{equation}}
\def\theTh{\Roman{section}.\arabic{Th}}
\def\theProp{\Roman{section}.\arabic{Prop}}
\def\theCo{\Roman{section}.\arabic{Co}}

\def\theRm{\Roman{section}.\arabic{Rm}}

\newcommand{\reset}{\setcounter{equation}{0}\setcounter{Th}{0}\setcounter{Prop}{0}\setcounter{Co}{0}
\setcounter{Lm}{0}\setcounter{Rm}{0}}
\def\La{\Lambda}
\def\La{\Lambda}

\def\lf{\left}
\def\rg{\right}

\def\al{\alpha}
\def\la{\lambda}

\def\ep{\varepsilon}

\def\ds{\displaystyle}
\def\ov{\overline}
\def\Om{\Omega}
\def\om{\omega}
\def\p{\partial}

\def\va{\vec{a}}

\def\res{\mathop{\hbox{\vrule height 7pt width .5pt 
depth 0pt\vrule height .5pt width 6pt depth 0pt}}\nolimits}

\begin{document}

\title{Minmax Hierarchies, Minimal Fibrations and a PDE based Proof of the Willmore Conjecture}

\author{ Tristan Rivi\`ere\footnote{Department of Mathematics, ETH Zentrum,
CH-8093 Z\"urich, Switzerland.}}

%\date{ }
\maketitle

{\bf Abstract :}{\it We introduce a general scheme that permits to generate successive min-max problems for producing critical points of higher and higher indices to Palais-Smale Functionals in normal Banach manifolds equipped with complete Finsler structures. We call the resulting tree of minmax problems a {\bf minmax hierarchy}. We give several examples and in particular we explain how to implement this scheme in the framework of the {\bf viscosity method} introduced by the author in \cite{Riv-minmax} in order to give a new proof of the Willmore conjecture after the famous result by Marques and Neves.}

\medskip

\noindent{\bf Math. Class. 49Q05, 53A10,  58E12, 49Q10}

\section{Introduction}

While Calculus of Variations takes it's roots in the works of Bernouilli, Euler and Lagrange from the XVIIIth century it is only almost two centuries later that  variational methods have been devised to ``detect'' critical points of functionals which are not absolute minimizers. The work of Georges Birkhoff \cite{Birk} is a pioneered contribution to minmax variational strategies in proposing a full resolution of a {\it Mountain Pass Problem}. In this work the author is introducing a ``curve shortening process'' that eventually brings an initial sweep-out of an arbitrary 2-sphere $(S^2,g)$ by closed curves to a more optimal one with reduced maximal length of these curves. The iteration of this process eventually generates a non trivial closed geodesic of non zero Morse Index (for generic metrics) on any  simply connected closed surface (see also \cite{CM} for a modern presentation of Birkhoff minmax process).

In a seminal work from 1970 Richard Palais is giving the foundation of a modern general minmax theory also known as {\it Palais-Smale Deformation theory in infinite dimensional spaces} (see \cite{Palais}). Non zero Morse index critical points of suitable Lagrangians in regular\footnote{The adjective ``regular'' is refering here to separation axioms.}  Banach Manifolds equipped with {\it complete Finsler structures} are detected through the deformation of the level sets by the mean of a {\it pseudo-gradient flow}. The theory has found a wide range of applications in particular in the field of non-linear elliptic partial differential equations. {\it Palais-Smale's pseudo gradient flow} is in a sense an abstract generalisation of the somehow more constructive {\it Birkhoff curve shortening process} which is produced by the mean of combinatorics arguments while the {\it pseudo-gradient flow} comes out of an abstract construction. Palais's {\it pseudo-gradient flow}  could be seen also as the ``ancestor'' of the numerous geometric and non geometric {\it gradient flows} which have proven their efficiencies for solving long standing conjectures...

One of the limitation of the range of application of {\it Palais-Smale Theory} is the requirement that the Lagrangian, whose critical points are the object of studies, satisfies what the author's at the time called the {\it (C) condition} (known nowadays under the name of {\it Palais-Smale Condition}). This condition is a stability conditions which says roughly that  a sequence of points in the Banach manifold, below some energy level, and for which the derivative of the Lagrangian is tending to zero converges\footnote{With respect to the distance induced by the Finsler structure} (modulo extraction of a subsequence) for the distance induced by the Finsler structure towards a critical point of the Lagrangian. The problem is that many variational problems which are relevant to geometric applications does not satisfy the Palais-Smale condition due to their invariance under non compact groups of conformal transformations (Yang-Mills in 4 dimensions, Yamabe, Harmonic maps in 2 dimensions, Willmore surfaces...etc). In order to overcome the lack of Palais-Smale condition one strategy consists in regularizing the studied Lagrangian by adding some 
coercive enough relaxation which permits the use of Palais-Smale theory. This relaxation is preceded by a small ``viscosity'' parameter that is sent to zero once the minmax critical points have been obtained for the regularized Lagrangian by the Palais-Smale theory. This approach, also called {\it viscosity method}, which is very much \underbar{based on the analysis of Partial differential Equations},  has been successfully implemented for the area of immersions of surfaces in \cite{Riv-minmax, MiRi,Riv-reg,Riv-Ind,Mich,PiRi1,PiRi2,Pi1,Riv-Lag-Visc} and the authors obtain the realization of any non trivial minmax problem by a possibly branched smooth minimal immersion satisfying various properties (Morse index bound, genus bound, free-boundary property, Lagrangian property...etc). 

With an abstract existence  theorem at hand the question is then to produce minmax problems to which it can be applied. The traditional approach is of ``Morse theoretic nature'' and consists in taking advantage of the topology of the {\it configuration space} to devise minmax strategies based on the realization of non trivial homotopy, homology or cohomology classes of this space (sometimes below some energy level). For instance, for minimal surfaces, Fernando Cod\'a Marques and Andr\'e Neves with collaborators studied the infinite sequence of minmax problems generated by the so called {\it Gromov-Guth widths} which are the successive cup products of the generator of the ${\Z}_2-$cohomology of ${\Z}_2-$hypercycles in a given closed manifold \cite{MN1,MN2,MN3,MN4}. This series of works have been crowned with the proof by Antoine Song  of the Yau Conjecture about the existence of infinitely many distinct minimal hypersurfaces in low dimensions \cite{Son}.

In unpublished and non-submitted working notes \cite{Riv-hie}, the author introduced some years ago a strategy for generating successive minmax problems which is not exactly of Morse theoretic nature. This strategy called ``minmax hierarchy'' consists in taking advantage of the topology of the space of solutions to a given minmax operation to generate a or several new ones (depending on the topology) of strictly higher widths and higher dimension (see theorem~\ref{th-I.1} below). The iteration of this procedure is generating a tree of minmax problems.
The goal of the present work is to give the heuristic of this strategy. We first start by giving in section II two elementary examples of almost explicit hierarchies : a construction of the successive eigen-spaces of the Laplace Beltrami operator on a closed manifold not using the Rayleigh quotient method as well as the construction of closed geodesics on ellipsoids. In section III we give a further more intricate example of hierarchy for the Dirichlet energy of maps between $S^3$ and $S^2$. We discuss moreover the question of solving the associated minmax problems by the mean of the viscosity method based on a {\it Ginzburg-Landau relaxation}. In section IV we combine the viscosity method and the notion of hierarchy to provide with a new proof of the Willmore conjecture not using Almgren-Pitts {\it almost minimizing varifold} theory but instead  based on the use of maps and PDEs. Finally we conclude this work by a last chapter on the abstract notion of hierarchy based on ${\Z}_2-$cohomology.

\section{Two Elementary Examples of Minmax Hierarchies.}
\reset
\subsection{Minmax Hierarchies in the Linear case of Laplace Eigenvalues.}
\label{rayleigh}
A classical variational approach to the eigenvalue problem for the Laplacian on a closed oriented riemannian manifold $(N^n,h)$  is given by the {\it Rayleigh quotient method}.
It can be sketched as follows.

Introduce the {\it Hilbert Sphere}
\[
{\mathfrak S}:=\lf\{u\in W^{1,2}(N^n)\ ;\ \int_{N^n}|u|^2 \ dvol_h=1\rg\}
\]
We consider first
\[
\la_1:=\inf_{u\in {\mathfrak S}}E(u):=\int_{N^n}|du|^2_h \ dvol_h\quad\mbox{ and }\quad{\mathcal C}_1:=\lf\{u\in {\mathfrak S}\ ;\ \Delta_h u=\la_1 u\rg\}
\]
where $\Delta_h$ is the positive {\it Beltrami Laplace Operator} on $N^n$. Iteratively we introduce $E_{k-1}:=\mbox{Span}\,{{\mathcal C}_{k-1}}\oplus E_{k-2}$ (with the convention $E_0=\{0\}$)
\[
\la_k:=\inf_{u\in {\mathfrak S}\cap E_{k-1}^\perp}E(u):=\int_{N^n}|du|^2_h \ dvol_h\quad\mbox{ and }\quad{\mathcal C}_k:=\lf\{u\in {\mathfrak S}\ ;\ \Delta_h u=\la_k u\rg\}\simeq S^{n_k-1}
\]
where $E_{k-1}^\perp$ is the space orthogonal to $E_{k-1}$ for the $L^2$ scalar product and $n_{k}:=\mbox{dim} \,\mbox{Span}({\mathcal C}_{k})$.

Assume now that we don't want to make use neither of the linear nature of the problem nor on the existence of the scalar product. An alternative way to obtain the {\it Laplace Eigenspaces}
and {\it Laplace Eigenvalues} is given by what we call a {\it Minmax Hierarchy}. A {\it Minmax Hierarchy} in this  framework is the following iterative construction.
Starting from $\la_1$, ${\mathcal C}_1=\{-u^1,u^1\}$ and $1=n_1:=\mbox{dim} \,\mbox{Span}({\mathcal C}_1)$ which were obtained by a strict minimization of the {\it Dirichlet energy} in ${\mathfrak S}$, in order to produce $\la_2$ without using the scalar product one could proceed as follows. Introduce
\[
\mbox{Sweep}_1:=\lf\{ u\in \mbox{Lip}([-1,+1],{\mathfrak S}) \quad\mbox{ s. t. }\quad u_{-1}=-u^{1}\quad\mbox{ and }\quad u_{+1}:=u^1\rg\}
\]
One has
\[
\la_2=\inf_{u\in \mbox{Sweep}_1}\ \max_{y\in[-1,+1]} \, E(u_y)
\]
Using the classical {\it Palais deformation theory} one produces critical point of $E$ in ${\mathfrak S}$ realizing $\la_2$. Introduce ${\mathcal C}_2:=\lf\{u\in {\mathcal S}\ ;\ \Delta_h u=\la_2 u\rg\}$ and call $n_2:=\mbox{dim} \,\mbox{Span}({\mathcal C}_{2})$ and $N_2=n_1+n_2$. It is clear that by taking
a minmax based on the space of 1 dimensional paths connecting elements from ${\mathcal C}_2$ would lead to nothing since points in ${\mathcal C}_2$ can be connected within ${\mathcal C}_2\simeq S^{n_2-1}$ (i.e. $n_2>1$). We introduce instead $Y_2:=B^{n_2}\times B^{n_1}$
\[
\mbox{Sweep}_{N_2}:=\lf\{ 
\begin{array}{l}
\ds u\in \mbox{Lip}(Y_2,{\mathcal S}) \ ;\ \forall \, z\in \p B^{n_2}\quad u_{(z,\cdot)}\in \mbox{Sweep}_1 \\[3mm]
(u_{\cdot})^{-1}({\mathcal C}^2)\cap \p Y_2=\p B^{n_2}\times \{0\}\quad;\quad \mbox{deg}(u_{|\p B^{n_2}\times \{0\}})=+1
\end{array}
\rg\}
\]
One then proves (see subsection~\ref{rayleigh})
\[
\la_3=\inf_{u\in \mbox{Sweep}_{N_2}}\max_{y\in Y_2} E(u_y)
\]
We construct inductively 
\begin{itemize}
\item[i)] A sequence of integers $n_k\in {\N}^\ast$ 
\[
n_k:=\mbox{dim} \,({\mathcal C}_k)+1
\]
\item[ii)] A sequence $\mbox{Sweep}_{N_k}\subset \mbox{Lip}(Y_k,{\mathcal S})$ where $Y_k:= B^{n_k}\times B^{n_{k-1}}\cdots B^{n_1} $ characterized as follows
\[
\forall\, u\in\mbox{Sweep}_{N_k}\quad \forall z\in  \p B^{n_k}\quad\quad u_{(z,\cdot)}\in \mbox{Sweep}_{N_{k-1}}
\]
\item[iii)] we have $u^{-1}({\mathcal C}_{k-1})\cap \lf(\p B^{n_k}\times Y_{k-1}\rg)= \p B^{n_k}\times\{0\}$, 
\[
\mbox{deg}(u_{|_{\p B^{n_k}\times\{0\}}})=+1
\]
\item[iv)] for $\mbox{dim}\,Y_{k-1}>0$ (i.e. $k>2$)
\[
\max_{y\in B^{n_k}\times\ \p Y_{k-1}} \int_{N^n}|du_y|^2_h \ dvol_h< \la_{k-1}
\]
\end{itemize}
With this definition for $\mbox{Sweep}_{N_k}$ we have the following lemma
\begin{Lm}
\label{lm-rayleigh}
\[
\la_k=\inf_{u\in \mbox{Sweep}_{N_k}}\max_{y\in Y_k} \int_{N^n}|du_y|^2_h \ dvol_h\quad.
\]
\hfill $\Box$
\end{Lm}
\noindent{\bf Proof of Lemma~\ref{lm-rayleigh}.}
Because of theorem~\ref{th-I.1} we have a strictly increasing sequence of eigenvalues $\mu_k$ such that
\[
\mu_k:=\inf_{u\in \mbox{Sweep}_{N_k}}\max_{y\in Y_k} E(u_y)
\]
The problem is to show that we indeed ``capture'' all the successive eigenvalues of the Laplacian. It suffices to show then that
\be
\label{eigenv}
\mu_k\le\la_k\quad.
\ee
We prove (\ref{eigenv}) by induction. Let $F_l$ be the vector eigenspace associated to $\la_l$. By induction assumption we have $n_l=\mbox{dim} F_l=\mbox{dim}\,{\mathcal C}_l+1$
We take for $\om_l$ the generator of $H^{n_l-1}({\mathcal C}_l,{\Z}_2)\simeq{\Z}_2$ (we obviously have ${\mathcal C}_l\simeq S^{n_l-1}$. Let $u_1\cdots u_k$ be an arbitrary
choice of eigenfunctions in ${\mathfrak S}$ for $\la_1<\cdots<\la_k$. We denote by $SO(F_l)$ the space of positive isometries of the euclidian spaces $F_l$.

\medskip

Let $u\ :\ (-1,+1)^k\times SO(F_2)\times\cdots SO(F_{k-1})\longrightarrow {\mathfrak S}$ given by
\[
\begin{array}{l}
u(t_1\cdots t_k,R_2\cdots R_{k-1}):=\\[3mm]
\ds\lf\{
\cos(\frac{\pi t_k}{2})\, u_k+\sin(\frac{\pi t_k}{2})\lf(\cos(\frac{\pi t_{k-1}}{2})\, R_{k-1}\,u_{k-1}+\sin(\frac{\pi t_{k-1}}{2})\,\lf(\cdots\lf(\cos(\frac{\pi t_2}{2})\, R_{2}\,u_{2}+\sin(\frac{\pi t_2}{2})\, u_1\rg)\rg)\rg)
\rg\}
\end{array}
\]
It is straightforward to check that, assuming $\mu_l=\la_l$ for $l<k$ we have $u\in \mbox{Sweep}_{N_k}$ and moreover
\[
\max_{(t_1\cdots t_k,R_2\cdots R_{k-1})\in (-1,+1)^k\times SO(F_2)\times\cdots SO(F_{k-1})}E\lf(u(t_1\cdots t_k,R_2\cdots R_k)\rg)\le \la_k
\]
This proves (\ref{eigenv}).

\subsection{Minmax Hierarchy for the Length on an Ellipsoid}
We consider the ellipsoid ${\mathcal E}^2(a,b,c)$ where $0<a<b<c$ given by
\[
{\mathcal E}^2(a,b,c):=\lf\{(x_1,x_2,x_3)\in {\R}^3\quad;\quad \frac{x_1^2}{a^2}+\frac{x_2^2}{b^2}+\frac{ x_3^2}{c^2}=1\rg\}
\]
Simple closed geodesics in constant speed parametrization  on ${\mathcal E}^2$ are the critical points of the Dirichlet energy $E$ among maps in the {\it Hilbert Manifold}
\[
{\mathfrak M}:=\lf\{u\in W^{1,2}(S^1, {\R}^3)\ ; \ u(e^{i\theta})\in {\mathcal E}^2\quad \forall \,\theta\in [0,2\pi]\rg\}\quad.
\]
The first  step in the hierarchy consists in introducing the space 
\[
\mbox{Sweep}_1({\mathcal E}^2):=\lf\{u\in C^0([0,1],{\mathfrak M})\ ;\ u_\ast([0,1]\times S^1)\quad\mbox{generates } H_2({\mathcal E}^2;{\Z})\rg\}
\]
The corresponding {\it width} is given by
\[
{\mathcal W}_1({\mathcal E}^2):=\inf_{u\in \mbox{Sweep}_1}\quad\max_{t\in [0,1]}\lf[2\pi\, \int_{S^1}|\p_\theta u_t|^2\ d\theta\rg]^{1/2}
\]
and is known to be achieved by the horizontal ellipse $\{x_3=0\}$ and the $\mbox{\it inf}_{u\in {\mathfrak M}}$ is in fact a $\mbox{\it min}_{u\in {\mathfrak M}}$.
\begin{figure}
\begin{center}
%\psfrag{A}{$[0,1]\times S^1$}
\includegraphics[width=12cm,height=7cm]{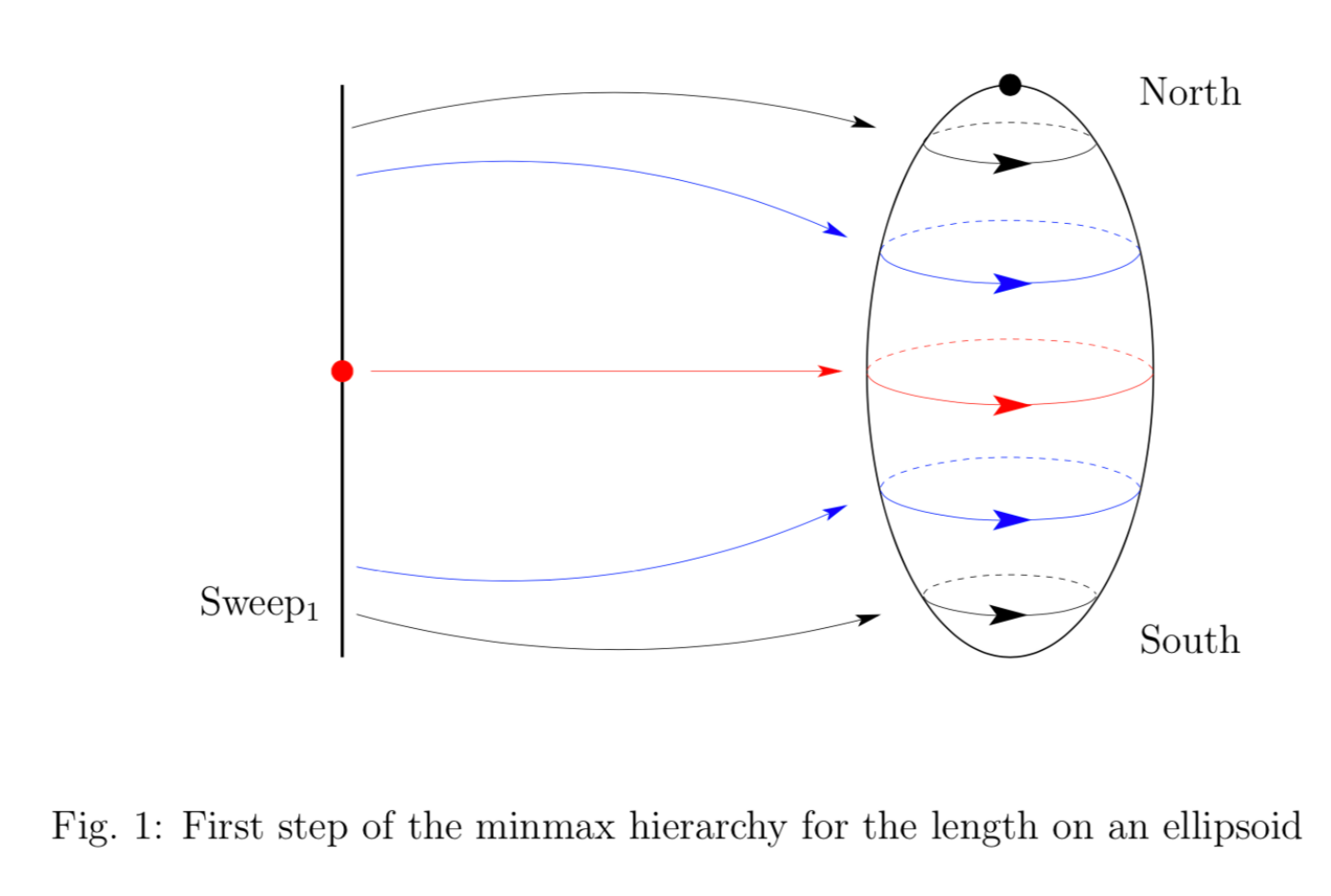}
\end{center}
\end{figure}

 The horizontal ellipse
can be given the 2 orientations. In other words we have that the space ${\mathcal C}_1$ of solutions to the previous minmax operation is disconnected and is made of \underbar{two disjoint components} each of them diffeomorphic to $S^1$ : the constant speed parametrizations of the horizontal ellipse in one direction and the constant speed parametrizations of the horizontal ellipse in the other direction. We call $u^\pm_t$ the choice of the two optimal $\mbox{Sweep}_1$ satisfying
\[
\max_{t\in [0,1]}\lf[2\pi\, \int_{S^1}|\p_\theta u_t^\pm|^2\ d\theta\rg]^{1/2}=W_1({\mathcal E}^2)
\]
and given by the constant speed parametrizations of the families of slices between the horizontal planes and the ellipsoid with the two orientations. Introduce $\mbox{Sweep}_2$ to be the following space
\[
\mbox{Sweep}_2({\mathcal E}^2):=\lf\{ U\in C^0([-1,+1],\mbox{Sweep}_1({\mathcal E}^2))\quad\mbox{ s. t. }\quad U(\pm1)=u^\pm\right\}
\]
We are taking advantage of the non-triviality of $\pi_0({\mathcal C}_1$ in order to build the next step in the hierarchy. The figure 2 illustrates the fact that $\mbox{Sweep}_2({\mathcal E}^2)\ne \emptyset$. 

The second with in the hierarchy is then given by
\[
W_2({\mathcal E}^2):=\inf_{U\in \mbox{Sweep}_2}\quad\max_{t\in [0,1], s\in [-1,+1]}\lf[2\pi\, \int_{S^1}|\p_\theta U_{s,t}|^2\ d\theta\rg]^{1/2}
\]
It is  achieved by the second smallest meridian given by $\{x_2=0\}$. Two ``optimal'' $\mbox{Sweep}_2$ that we denote by $U^\pm$ are obtained as follows : we consider the two paths  in $S^2$ joining $\vec{k}$ and $-\vec{k}$ given respectively by
$\vec{k}_s^{\,\pm}:=\cos\pm(\pi/2+s\pi/2)\ \vec{k}+\sin\pm(\pi/2+s\pi/2)\ \vec{i}$ for $s\in [-1,+1]$ and $U^\pm$ are obtained respectively by slicing the ellipsoid by the planes orthogonal to $\vec{k}_s^{\,\pm}$ as $s$ evolves in $(-1,+1)$.  
If this is the case one would construct $\mbox{Sweep}_3$ as follows
\[
\mbox{Sweep}_3({\mathcal E}^2):=\lf\{ {\mathcal U}\in C^0([-1,+1],\mbox{Sweep}_2({\mathcal E}^2))\quad\mbox{ s. t. }\quad{\mathcal U}(\pm1)=U^\pm\right\}
\]
\begin{figure}
\begin{center}
%\psfrag{A}{$[0,1]\times S^1$}
\includegraphics[width=15cm,height=9cm]{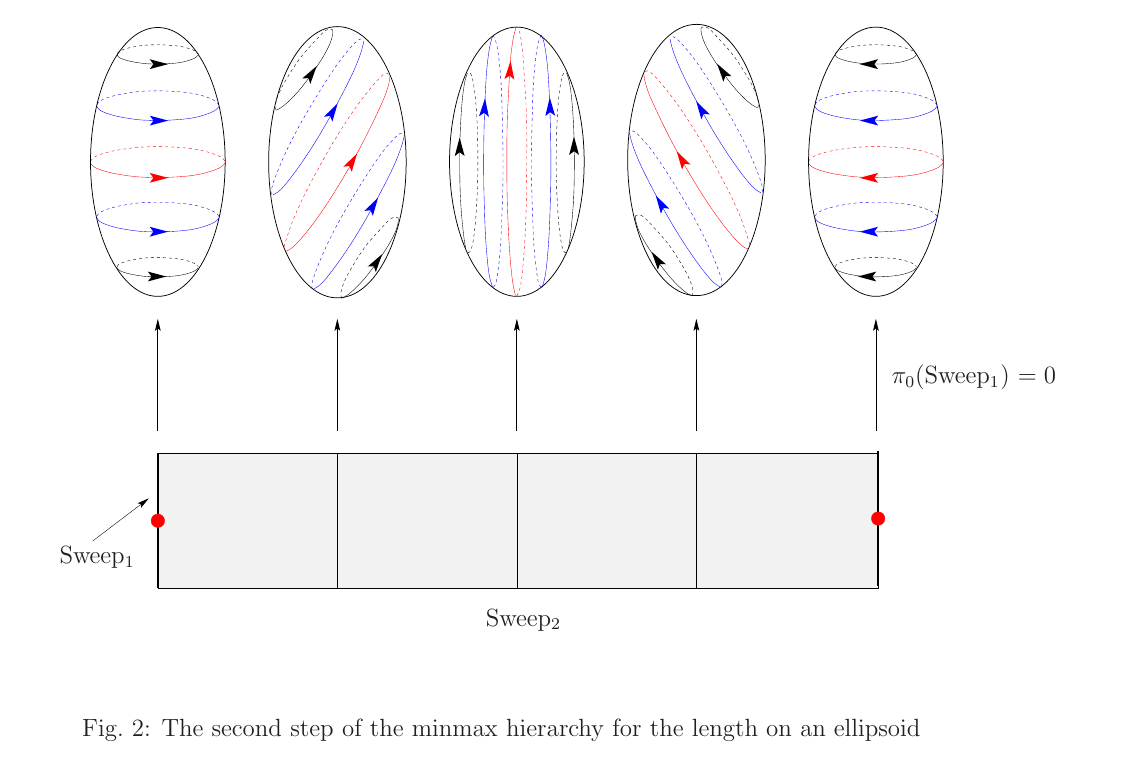}
\end{center}
\end{figure}

\begin{figure}
\begin{center}
\centerline{\includegraphics[width=17cm,height=18cm]{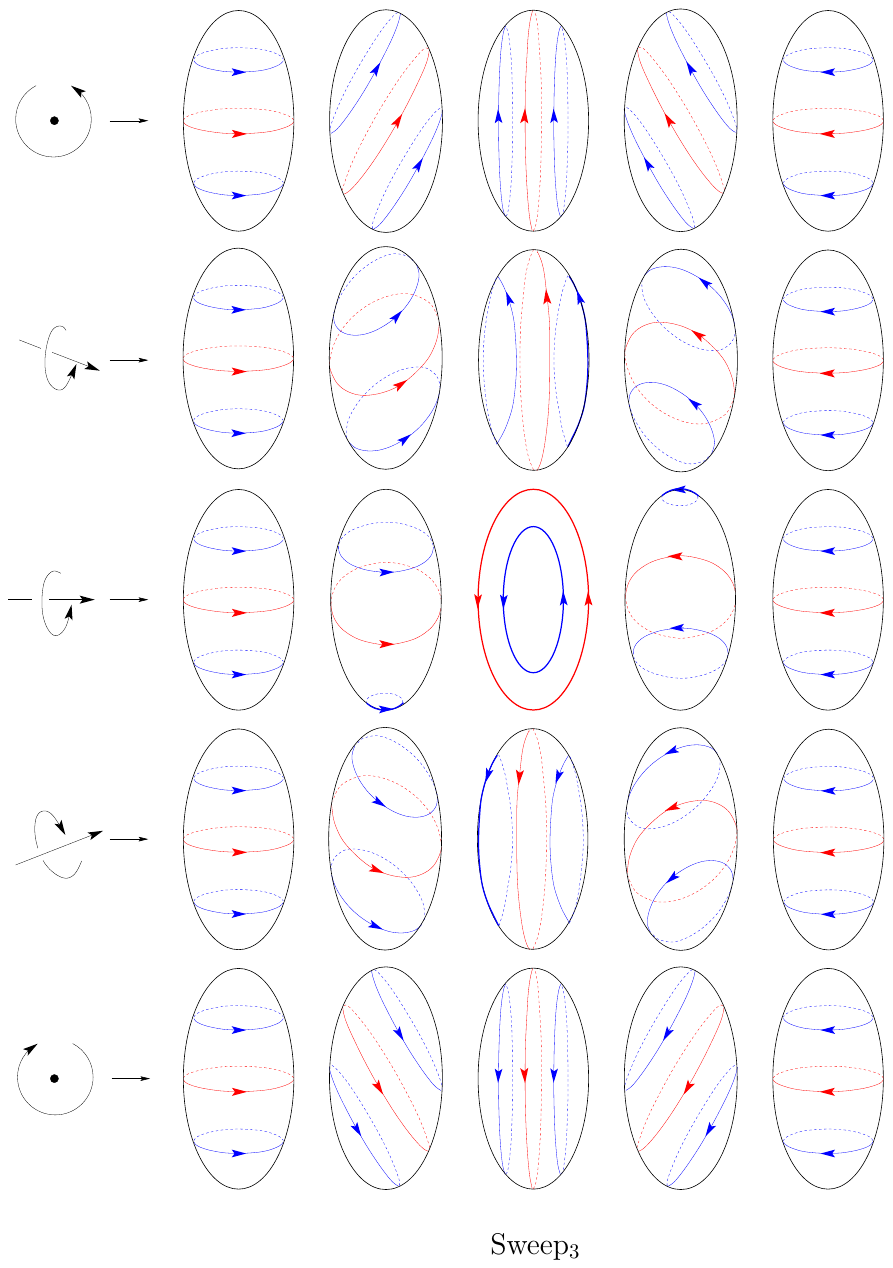}}
\end{center}
\end{figure}
Here also we take advantage of the non-triviality of $\pi_0({\mathcal C}_2)$, where ${\mathcal C}_2$ is the space of solutions to $W_2$,  in order to build the next step in the hierarchy. We have $\mbox{Sweep}_3({\mathcal E}^2)\ne \emptyset$. Indeed for $(s,\sigma)\in [-1,1]^2$ one considers the slicing of the ellipsoid by the planes orthogonal to 
\[
\vec{k}_{\sigma,s}:= \cos(\pi/2+s\,\pi/2)\ \vec{k}+\sin(\pi/2+s\,\pi/2)\ (\cos(\pi/2+\sigma\pi/2)\, \vec{i}+\sin(\pi/2+\sigma\pi/2)\, \vec{j})
\]
Finally we introduce
\[
{\mathcal W}_3({\mathcal E}^2):=\inf_{u\in \mbox{Sweep}_3}\quad\max_{t\in [0,1], s\in [-1,+1], \sigma\in[-1,+1]}\lf[2\pi\, \int_{S^1}|\p_\theta u_{t,s,\sigma}|^2\ d\theta\rg]^{1/2}
\]
It is  achieved by the largest meridian given by $\{x_1=0\}$. Hence the 3 first steps in the hierarchy is giving us 3 embedded geodesics (which happens to be the only ones for some values of $(a,b,c)$ according to a classical result by M.Morse - see also \cite{Vies}). 

The reason why we cannot get a simple $4th$ one is perhaps related to the fact that $\mbox{Sweep}_4$ defined as above is empty unless one adds a second component of $S^1$ and look for non connected geodesics.

\section{Minmax Hierarchy for  Maps in $W^{1,2}(S^3,S^2)$.}
\reset
The zero level of the hierarchy that we are considering is given by the minimization of the Dirichlet energy of maps between $S^3$ and $S^2$. There are as many minimizers as
constant maps hence we have
\[
{\mathcal C}_0\simeq S^2\quad
\]
and
\[
{\mathcal W}_0(S^3,S^2):=\inf_{u\in W^{1,2}(S^3,S^2)}E(u)=0\quad.
\]
We are now exploiting the $\pi_3({\mathcal C}_0)={\Z}\ne 0$ in order to construct the next step in the hierarchy. We introduce
\[
\mbox{Sweep}^1_1(S^3,S^2):=\lf\{ 
\begin{array}{c}
u\in C^0(\ov{B^4}, W^{1,2}(S^3,S^2))\cap C^0(B^4\times S^3, S^2)\quad; \quad\forall\ a\in \p B^4\quad u_a\in {\mathcal C}_0\\[5mm]
 a\in \p B^4\longrightarrow u_a\in S^2\quad\mbox{ is non zero homotopic}
\end{array} 
\rg\}
\]
It is first important to observe that
\begin{Lm}
\label{lm-minmax-map-1}
\[
\mbox{Sweep}^1_1(S^3,S^2)\ne\emptyset\quad,
\]
and the corresponding width satisfies
\be
\label{minmax-map-1}
{\mathcal W}^1_1(S^3,S^2):=\inf_{u\in \mbox{Sweep}^1_1} \quad\max_{a\in B^4} E(u_a)>0\quad.
\ee
\hfill $\Box$
\end{Lm}
A priori if we would have required $u$ to be continuous on $\ov{B^4}\times S^3$ we  would have got that $$\mbox{Sweep}^1_1(S^3,S^2)=\emptyset\quad.$$ indeed we have an exact sequence 
\[
\cdots \longrightarrow \pi_n(\Om_3(S^2)) \longrightarrow 
\begin{tikzcd}\pi_n(C^0(S^3,S^2)) \arrow[r, "\mathrm{ev}_{\ast}"'] &   \ar[l,bend left=80,looseness=1.8, "\iota_{\ast}"] \pi_{n}(S^2)
\end{tikzcd}\longrightarrow \pi_{n-1}(\Om_3(S^2))\longrightarrow\cdots
\]
where $\Om_3(S^2)$ denotes the space of base point preserving continuous maps from $S^3$ into $S^2$ and $\mathrm{ev}$ is the evaluation map at a based point and $\iota$ is the map which to a point assigns the constant map having this value. Since obviously $\iota$ realizes a section of the fibration given by $\mathrm{ev}$, we have that $\mathrm{ev}_\ast\circ\iota_\ast$ is the identity and the homotopy sequence splits.  This gives
\[
\pi_n(C^0(S^3,S^2))=\pi_n(S^2)\oplus\pi_{n}(\Om_3(S^2))=\pi_n(S^2)\oplus\pi_{n+3}(S^2)\quad.
\]
In the case $n=3$ the boundary condition we are imposing for the membership in $\mbox{Sweep}^1_1(S^3,S^2)$ is nothing but saying that  $a\in \p B^4\rightarrow u(a,\cdot)=\mathrm{ev}(u(a,\cdot))$ is a non zero element
in $\pi_3(C^0(S^3,S^2))$. 

\medskip

\noindent{\bf Proof of lemma~\ref{lm-minmax-map-1}.}
To any  non homotopically trivial  smooth map $u\in C^\infty(S^3,S^2)$ we associate it's {\it canonical family} given by
\be
\label{admissible-B^4}
a\in B^4\ \ \longrightarrow\ u_a:=u\circ\phi_a\in C^\infty(S^3,S^2)\quad \mbox{where}\quad \phi_a(z):=(1-|a|^2)\frac{z-a}{|z-a|^2}-a
\ee
It is not difficult to see that for any $b\in \p  B^4$  we have
\[
W^{1,2}-\lim_{a\in B^4\rightarrow b}u_a=u(-b)\in {\mathcal C}_0
\]
Since  $u$ has a non zero {\it Hopf degree} we have that $b\in \p B^4\rightarrow u(-b)$ is non zero homotopic. Hence $u_a\in \mbox{Sweep}_1(S^3,S^2)$.

\medskip

We now prove that ${\mathcal W}^1_1>0$. This is mostly coming from the Poincar\'e inequality : There exists $C_{S^3}>0$ such that, for any $u\in  W^{1,2}(S^3,S^2)$ we have
\[
\int_{S^3}\lf|u(x)-\frac{1}{|S^3|}\int_{S^3}u(y)\ dvol_{S^3}(y)\rg|^2\ dvol_{S^3}(x)\le\, C_{S^3}\ \int_{S^3}|du|^2_{S^3}\ dvol_{S^3}\quad.
\]
Let $\delta>0$ to be fixed later, assuming ${\mathcal W}^1_1<\delta$ we obtain the existence of $u\in \mbox{Sweep}^1_1$ such that
\[
\max_{a\in \ov{B^4}}E(u_a)<\delta\quad,
\]
Poincar\'e inequality gives
\[
\max_{a\in \ov{B^4}}\int_{S^3}\lf|u_a(x)-\frac{1}{|S^3|}\int_{S^3}u_a(y)\ dvol_{S^3}(y)\rg|^2\ dvol_{S^3}(x)\le\ C_{S^3}\, \delta\quad.
\]
For $\delta$ chosen small enough this implies
\[
\min_{a\in \ov{B^4}}\lf|\frac{1}{|S^3|}\int_{S^3}u_a(y)\ dvol_{S^3}(y)\rg|>1/2
\]
and we deduce that the map
\[
a\in B^4\longrightarrow\ \int_{S^3}u_a(y)\ dvol_{S^3}(y)/\lf|\int_{S^3}u_a(y)\ dvol_{S^3}(y)\rg|
\]
is  an extension of $a\in \p B^4\rightarrow \ u_a\in S^2$ in $C^0(B^4,S^2)$. This contradicts the fact that the restriction of $u_a$ on $\p B^4$ is non zero homotopic.

It is expected minmax~\ref{minmax-map-1} to be achieved by a non-constant smooth harmonic map of Morse index less or equal than $4$. The author recently proved in \cite{Riv-Morse} that they
are all of the form
\be
\label{index-4}
\varphi\circ{\mathfrak h}\circ S\quad,
\ee
where $S\in O(4)$, $\mathfrak h$ is the Hopf fibration and $\varphi$ is an arbitrary holomorphic diffeomorphism of ${\C}{\mathbb P}^1$. We conjecture that
 ${\mathcal C}^1_1$ is exactly given by this space of maps this would imply a positive answer to the following question\footnote{At this stage we only know that
 \be
 \label{lower-bound}
 {\mathcal W}_1^1\le 16\pi^2\quad.
 \ee
Indeed in \cite{ElS} it is proved that for any smooth harmonic map $u$ from $S^3$ into $S^2$ such as $u:={\mathfrak h}$ we have
\be
\label{minmax-hopf}
\max_{a\in B^4}E(u\circ\phi_a)=E(u)\quad,
\ee
and we have $|d{\mathfrak h}|_{S^3}\equiv \sqrt{2}\, 2$ and hence, since $\int_{S^3} dvol_{S^3}=2\,\pi^2$, we deduce $E({\mathfrak h})=16\,\pi^2$ and hence, thanks to (\ref{minmax-hopf}) we have (\ref{lower-bound}). }
\[
{\mathcal W}_1^1=16\,\pi^2\quad ?
\]
A positive answer to this question would imply a proof of the conjecture made by the author in \cite{Riv-Hopf} according to which the Hopf fibration minimizes the 3-energy among non zero homotopic maps from $S^3$ into $S^2$
(see \cite{Riv-Morse}).

\medskip

Observe that, starting from ${\mathcal C}_0\simeq S^2$, there is another branch for a minmax hierarchy construction based on the $\pi_2(S^2)$ this time. Indeed, we introduce\footnote{A priori if we would have required $u$ to be continuous on $\ov{B^3}\times S^3$ we  would have got also here that $$\mbox{Sweep}^2_1(S^3,S^2)=\emptyset\quad.$$ }
\[
\mbox{Sweep}^2_1(S^3,S^2):=\lf\{ 
\begin{array}{c}
u\in C^0(\ov{B^3}, W^{1,2}(S^3,S^2))\quad; \quad\forall\ a\in \p B^3\quad u_a\in {\mathcal C}_0\\[5mm]
 a\in \p B^3\longrightarrow u_a\in S^2\quad\mbox{ is non zero homotopic}
\end{array} 
\rg\}
\]
We  have the following lemma whose proof is maybe less direct than the proof of lemma~\ref{lm-minmax-map-1}.
\begin{Lm}
\label{lm-minmax-map-2}
\[
\mbox{Sweep}^2_1(S^3,S^2)\ne\emptyset\quad,
\]
and
\be
\label{minmax-map-2}
{\mathcal W}^2_1(S^3,S^2):=\inf_{u\in \mbox{Sweep}^2_1} \quad\max_{a\in B^3} E(u_a)>0\quad.
\ee
\hfill $\Box$
\end{Lm}
\noindent{\bf Proof of lemma~\ref{lm-minmax-map-2}.} We consider the following map
\[
\pi\ :\ (x_1,x_2,x_3,x_4)\in S^3\ \longrightarrow\ \frac{(x_1,x_2,x_3)}{\sqrt{1-x_4^2}}\ \in S^2
\]
This map has two point singularities on $S^3$ and $|d\pi|_{S^3}(x)\le \mbox{dist}^{-1}(x, \{North\}\cup\{South\})$. Hence $d\pi\in L^{3,\infty}(S^3)$ and in particular $d\pi\in L^2(S^3)$. We introduce the map in $C^0({B^3}, W^{1,2}(S^3,S^2))$ given by
\[
u\ :\ a=(a_1,a_2,a_3)\in B^3 \ \longrightarrow \ u_a(x)=\pi\circ \phi_{(a_1,a_2,a_3,0)}
\]
where $\phi_a$ is given by (\ref{admissible-B^4}). Observe that  for any $b=(b_1,b_2,b_3,0)\in \p  B^4\cap\{x_4=0\}$  we have
\[
W^{1,2}-\lim_{a\in B^4\rightarrow b}u_a=-\,b
\]
The map from $S^2$ into itself which to $b$ assigns $-b$ is non zero homotopic. Hence $\mbox{Sweep}^2_1(S^3,S^2)\ne\emptyset$. The proof of (\ref{minmax-map-2}) is identical to the one of (\ref{minmax-map-1}). This concludes the proof of lemma~\ref{lm-minmax-map-2}.\hfill $\Box$

\medskip
We conjecture this time that ${\mathcal C}_1^2$ is made of the maps of the form\footnote{It is also natural to expect that the Morse index of the singular harmonic morphism $\pi$ is $3$.}
\be
\label{index-3}
\varphi\circ\pi\circ S\quad,
\ee
where $S\in O(4)$ and $\varphi$ is an arbitrary holomorphic diffeomorphism of ${\C}{\mathbb P}^1$. Assuming this would be true this would answer positively to the question\footnote{At this stage we only know thanks to lemma~\ref{lma-w21} that 
\[
{\mathcal W}_1^2\le 8\,\pi^2\quad.
\]}

\[
{\mathcal W}_1^2= 8\,\pi^2\quad ?
\]
Assuming the space ${\mathcal C}^1_1$ of  solutions to the minmax problem ${\mathcal W}^1_1$ is exactly given by the space of maps of the form (\ref{index-4}),  we construct the next step in the hierarchy as follows. The space of holomorphic diffeomorphisms of ${\C}{\mathbb P}^1$ known  as $PSL(2,{\C})$ is homotopically equivalent to $SO(3)$ while the space of Hopf fibrations
is homeomorphic to ${\C}{\mathbb P}^1\times SO(3)$. We introduce
\[
\mbox{Sweep}_2^1(S^3,S^2):=\lf\{
\begin{array}{c}
\ds u\in C^0(\ov{B^3\times B^4}; W^{1,2}(S^3,S^2))\quad;\quad \forall b\in \p B^3,\quad u_{(b,a)}={\mathfrak h}_b\circ \phi_a\in Sweep_1^1(S^3,S^2)\\[5mm]
\ds \mbox{ where }b\in \p B^3\longrightarrow {\mathfrak h}_b\in {\mathcal C}^1_1\simeq {\C}{\mathbb P}^1\times SO(3)\quad\mbox{ is non zero homotopic}\\[5mm]
 \ds \sup_{(b,a)\in B^3\times\p B^4}E(u_{(b,a)})\le 8\pi^2
\end{array}
\rg\}
\]
If one can prove that $\mbox{Sweep}_2^1(S^3,S^2)\ne \emptyset$ it follows from the last section that
\[
{\mathcal W}_2^1:=\inf_{u_{b,a}\in\mbox{Sweep}_2^1} \ \max_{(b,a)\in B^3\times B^4}E(u_{(b,a)})>{\mathcal W}_1^1
\]
It is expected that ${\mathcal W}_2^1$ is achieved by an harmonic map of Morse index less or equal than $7$.
\begin{Lm}
\label{lm-sweep21}
Under the above notations we have
\be
\label{sweep21}
\mbox{Sweep}_2^1(S^3,S^2)\ne \emptyset\quad.
\ee
\hfill $\Box$
\end{Lm}
\noindent{\bf Proof of lemma~\ref{lm-sweep21}.} A generator of the $\pi_2$ of the space of Hopf fibration $\pi_2({\C}{\mathbb P}^1\times SO(3))=\pi_2({\C}{\mathbb P}^1)\oplus \pi_2(SO(3))={\Z}$ is given by the following mapping
\[
b\in S^2\quad\longrightarrow \quad v_{b}({\mathfrak q}):=-\,{\mathfrak q}\,b \, {\mathfrak q}^\ast\in C^\infty( S^2,{\mathcal C}_1^1)
\]
where we are using quaternion multiplication and ${\mathfrak q}$ denotes a unit quaternion and ${\mathfrak q}^\ast$ denotes it's conjugate. Let
\[
u_{(b,a)}({\mathfrak q}):=\phi_a({\mathfrak q})\, \pi(\phi_{(b_1,b_2,b_3,0)}({\mathfrak q}))\, (\phi_a({\mathfrak q}))^\ast\quad.
\]
We have
\[
\forall\, a \in B^4\quad \forall \,b\in \p B^3\quad W^{1,2}-\lim_{c\rightarrow b}u_{(c,a)}=v_b\circ\phi_a\quad.
\] 
Moreover, thanks to lemma~\ref{lma-w21}
\[
\limsup_{(b,a)\rightarrow B^3\times\p B^4}E(u_{(b,a)})=\max_{b\in B^3}E(\pi\circ\phi_b)=E(\pi)=8\pi^2\quad.
\]
Hence, $u_{(b,a)}\in \mbox{Sweep}_2^1(S^3,S^2)$  and lemma~\ref{lm-sweep21} is proved.\hfill $\Box$

\subsection{Analysis Tools for studying ${\mathcal W}_1^1$.}
A possibility to study the minmax problem ${\mathcal W}_1^1$ consists in introducing the {\bf Ginzburg-Landau relaxation} of the Dirichlet energy
\[
E_\ep(u):=\frac{1}{2}\int_{S^3} |du|^2+\frac{1}{2\,\ep^2}(1-|u|^2)^2\ dvol_{g_{S^3}}
\]
defined on the Hilbert vector space $W^{1,2}(S^3,{\R}^3)$. The Palais-Smale framework is available for $E_\ep$ and, by the mean of {\it Struwe Monotonicity Trick} one obtains the existence of
a sequence $\ep_k\rightarrow 0$ and a sequence of maps $u_{k}$ such that
\[
\Delta_{S^3}u_{k}+u_k\, (1-|u_k|^2)=0\quad,
\]
satisfying also
\be
\label{entropy}
\frac{1}{2\,\ep^2_k}\int_{S^3}(1-|u_k|^2)^2\ dvol_{g_{S^3}}=o\lf(\frac{1}{\log\ep_k^{-1}}\rg)\quad,
\ee
and
\[
E_{\ep_k}(u_k)=\inf_{u\in \mbox{Sweep}_1^1(S^3,{\R}^3)}\quad\max_{a\in B^4}E_{\ep_k}(u)
\]
where
\[
\mbox{Sweep}^1_1(S^3,{\R}^3):=\lf\{ 
\begin{array}{c}
u\in C^0(\ov{B^4}, W^{1,2}(S^3,{\R}^3))\cap C^0(B^4\times S^3, {\R}^3)\quad; \quad\forall\ a\in \p B^4\quad u_a\in {\mathcal C}_0\\[5mm]
 a\in \p B^4\longrightarrow u_a\in S^2\quad\mbox{ is non zero homotopic}
\end{array} 
\rg\}
\]
Moreover
\be
\label{indice}
E_{\ep_k}-Index(u_k)\le 4\quad.
\ee
We expect 
\[
\lim_{k\rightarrow +\infty} E_{\ep_k}(u_k)=16\pi^2
\]
Observe that the strong convergence of $u_k$ in $W^{1,2}$ should not be expected since ${\mathcal C}_1$ is not compact due to the action of the M\"obius group\footnote{Unless the Ginzburg-Landau relaxation combined with the entropy estimate (\ref{entropy}) is making a very special selection and is breaking the asymptotic M\"obius group action. Such gauge breaking effect by relaxation has been already observed in \cite{LMM1,LMM2}}. Hence the asymptotic analysis of $u_k$ is delicate.
The main results of \cite{LW1,LW2} should be combined with the two additional estimates  (\ref{entropy}) and (\ref{indice}).  

\medskip

Finally we would like to conclude this section  by stressing the fact that the relaxation of the minmax problem $ {\mathcal W}_1^2$ is certainly more delicate since it is not clear whether there exists elements in  $\mbox{Sweep}^2_1(S^3,S^2)$ which are strongly approximable by smooth maps in $C^0(\ov{B^3},W^{1,2}(S^3,S^2))$.

\section{Minmax Hierarchy for Lagrangian Surfaces.}
\reset
\subsection{Lagrangian Immersions in the Grassmann Manifold $Gr^+_2({\R}^4)$.}
We consider the Grassman manifold $G^+_2({\R}^4)$ of oriented 2-planes in ${\R}^4$ with the canonical metric. It is well known that this space is isometric to the product $S^2\times S^2$ that we also denote $S^2_+\times S^2_-$ for convenience and further uses. The identification goes as follows : An oriented two plane in ${\R}^4$ is given by a {\it simple unit 2-vector} of the form $\vec{a}\wedge\vec{b}$ and we consider
\[
\vec{a}\wedge\vec{b}\in G^+_2({\R}^4)\quad\longrightarrow\quad \frac{1}{\sqrt{2}}\, \lf(\vec{a}\wedge\vec{b}+\star(\vec{a}\wedge\vec{b})\,,\,\vec{a}\wedge\vec{b}-\star(\vec{a}\wedge\vec{b})\rg)\in \ S^2_+\times S^2_-
\]
where $\star$ is the canonical {\it Hodge operator} on the alternating algebra $\wedge^2{\R}^4$ into itself and $S^2_+$ and $S^2_-$ are the unit spheres of unit {\it self-dual} resp. {\it anti-self dual} 2-vectors in $\wedge^2{\R}^4$.

We equip this product with the following Symplectic form $\om_{S^2_+\times S^2_-}:=\pi_+^\ast\om_{S^2}+\pi_-^\ast\om_{S^2}$ where $\pi_\pm$ : $S^2_+\times S^2_-\rightarrow S^2_{\pm}$ are the projections onto respectively the first and the second components of the product $S^2_+\times S^2_-$ and $\om_{S^2}$ is the canonical positive volume form on $S^2$. In fact $\om_{S^2_+\times S^2_-}$ defines a K\"ahler structure whose associated complex structure
is given by $J_{(x,y)}(v,w)=(x\times v, w\times y)$. The tangent space  at a point $\vec{\mathfrak G}:=\vec{a}\wedge\vec{b}\in G^+_2({\R}^4)$ where $\vec{a}$ and  $\vec{b}$ are two unit vectors of ${\R}^4$ orthogonal to each other is given by
\[
T_{\vec{\mathfrak G}}S^2_+\times S^2_-=\lf\{\vec{\Gamma}:=\vec{v}\wedge\vec{b}+\vec{a}\wedge\vec{w}\quad;\quad \vec{v}\,  ,\, \vec{w}\in \lf(\mbox{Span}\{\vec{a},\vec{b}\}\rg)^\perp\rg\}
\]
and with this notation we have
\be
\label{sign-complex}
J(\vec{\mathfrak G})\vec{\Gamma}=\star\lf(\vec{v}^\perp\wedge\vec{b}+\vec{a}\wedge\vec{w}^\perp\rg)\quad\mbox{ where }\quad \vec{v}^\perp:=\star(\vec{v}\wedge\vec{a}\wedge\vec{b})\quad(\mbox{i. e. } \vec{a}\wedge \vec{b}\wedge \vec{v}\wedge \vec{v}^\perp=|v|^2\ \star 1 ).
\ee
An immersion  $\vec{\mathfrak G}=(\vec{G}^{\,+},\vec{G}^{\,-})$ from a surface $\Sigma$ into $S^2_+\times S^2_-$ being given, we say that $ \vec{\mathfrak G}$ defines a {\it lagrangian immersion} if
\[
\vec{\mathfrak G}^{\,\ast}\om_{S^2_+\times S^2_-}=(\vec{G}^{\,+})^\ast\om_{S^2}+(\vec{G}^{\,-})^\ast\om_{S^2}\equiv 0\quad\mbox{ on }\Sigma\quad.
\]
We have the following lemma
\begin{Lm}
\label{lm-lag-cond}
Let $\vec{\mathfrak G}$ be a smooth map from the disc $D^2$ into $G^+_2({\R}^4)$ and let $(\vec{a},\vec{b})\in S^3\times S^3$ be a smooth lift of $\vec{\mathfrak G}$ that is $\vec{a}\cdot\vec{b}\equiv 0$ and 
$\vec{\mathfrak G}=\lf((\vec{a}\wedge\vec{b})^+\,,\,(\vec{a}\wedge\vec{b})^-\rg)$. We have
\be
\label{lag-cond-0}
\vec{\mathfrak G}^{\,\ast}\om_{S^2_+\times S^2_-}=2\,d\vec{a}\dot{\wedge}d\vec{b}=2\, (\p_{x_1}\vec{a}\cdot\p_{x_2}\vec{b}-\p_{x_2}\vec{a}\cdot\p_{x_1}\vec{b})\ dx_1\wedge dx_2\quad.
\ee
Moreover, if $\vec{\mathfrak G}$ is lagrangian if and only if
\be
\label{lag-cond}
d\vec{a}\dot{\wedge}d\vec{b}=(\p_{x_1}\vec{a}\cdot\p_{x_2}\vec{b}-\p_{x_2}\vec{a}\cdot\p_{x_1}\vec{b})\ dx_1\wedge dx_2\ =0
\ee
and we have
\be
\label{degree}
\lf((\vec{a}\wedge\vec{b})^+\rg)^\ast\om_{S^2}=\star\ \frac{1}{2}\vec{a}\wedge\vec{b}\wedge\lf(d\vec{a}\wedge d\vec{a}+d\vec{b}\wedge d\vec{b}\rg)\quad.
\ee
\hfill $\Box$
\end{Lm}
\noindent{\bf Proof of lemma~\ref{lm-lag-cond}.} 
Let $(\vec{c},\vec{d})$ be a smooth lift of $\vec{\mathfrak G}^{\,\perp}$ the orthogonal to $\vec{\mathfrak G}$ and denote by $\pi_\perp$ the projection onto $\vec{\mathfrak G}^{\,\perp}$ the span of $\{\vec{c},\vec{d}\}$. From (\ref{sign-complex}) we have
\[
\forall\, \vec{X}\in \,\mbox{Span}(\vec{c},\vec{d})\quad\quad J_{S^2}({\vec{G}^{\,\pm}})\lf((\vec{a}\wedge\vec{X})^\pm\rg)=\pm(\vec{a}\wedge\vec{X}^{\perp})^\pm
\]
A direct computation gives\footnote{The - sign comes from the fact that the relation between the K\"ahler metric and it's associated K\"ahler form is given by $g(u,v)=\om(u,Jv)$ and hence $-g(u,Jv)=\om(u,v)$.}
\be
\label{pull+}
\begin{array}{l}
\ds-\lf((\vec{a}\wedge\vec{b})^+\rg)^\ast\om_{S^2}\,(X,Y)=\lf(d_X(\vec{a}\wedge\vec{b})^+,J_{\vec{G}^+}(d_Y(\vec{a}\wedge\vec{b})^+)\rg)\\[5mm]
\ds =\lf(\pi_\perp(d_X\vec{a})\wedge \vec{b}+\vec{a}\wedge \pi_{\perp}(d_X\vec{b}))^+, ((\pi_\perp(d_Y\vec{a}))^\perp\wedge\vec{b}+\vec{a}\wedge (\pi_\perp(d_Y\vec{b}))^\perp)^+\rg)\\[5mm]
\ds=2^{-1}\lf( \pi_\perp(d_X\vec{a})\wedge\vec{b}-\vec{a}\wedge (\pi_\perp(d_X\vec{a}))^\perp +\vec{a}\wedge \pi_\perp(d_X\vec{b})+(\pi_\perp(d_X\vec{b}))^\perp\wedge\vec{b}\rg.\\[5mm]
\lf., (\pi_\perp(d_Y\vec{a}))^\perp\wedge\vec{b}+\vec{a}\wedge \pi_\perp(d_Y\vec{a} )+\vec{a}\wedge (\pi_\perp(d_Y\vec{b}))^\perp-\pi_\perp(d_Y\vec{b})\wedge\vec{b}       \rg)\\[5mm]
\ds =d_X\vec{a}\cdot(\pi_{\perp}d_Y\vec{a})^\perp+d_X\vec{b}\cdot(\pi_\perp(d_Y\vec{b}))^\perp-d_X\vec{a}\cdot d_Y\vec{b}+d_Y\vec{a}\cdot d_X\vec{b}
  \end{array}
\ee
Similarly we have
\be
\label{pull-}
\begin{array}{l}
\ds-\lf((\vec{a}\wedge\vec{b})^-\rg)^\ast\om_{S^2}\,(X,Y)=\lf(d_X(\vec{a}\wedge\vec{b})^-,J_{\vec{G}^-}(d_Y(\vec{a}\wedge\vec{b})^-)\rg)\\[5mm]
\ds =-d_X\vec{a}\cdot(\pi_\perp(d_Y\vec{a}))^\perp-d_X\vec{b}\cdot(\pi_\perp(d_Y\vec{b}))^\perp-d_X\vec{a}\cdot d_Y\vec{b}+d_Y\vec{a}\cdot d_X\vec{b}
  \end{array}
\ee
Summing (\ref{pull+}) and (\ref{pull-}) gives (\ref{lag-cond-0}).
Assume now $\vec{\mathfrak G}$ is lagrangian, then the condition  $d\vec{a}\dot{\wedge} d\vec{b}=0$ implemented in (\ref{pull+}) gives (\ref{degree}) and this concludes the proof of lemma~\ref{lm-lag-cond}. \hfill $\Box$

\medskip

Let $\vec{\mathfrak G}$ be a lagrangian immersion. We define the associated {\it Lagrangian Jacobian} to be the function $C$ on $\Sigma$ given by
\[
(\vec{G}^{\,+})^\ast\om_{S^2}= C(x)\ dvol_{\vec{\mathfrak G}}\quad.
\]
It is proved
%\footnote{Observe that in \cite{CaUr} a different convention is taken. In their papers the authors consider $S^2_+\times S^2_-$ equipped with the K\"ahler form $\pi_+^\ast\om_{S^2}-\pi_-^\ast\om_{S^2}$. The main reason for our present choice and for the difference between the two conventions is due to the fact that in \cite{CaUr} the Gauss map of an immersion in $S^3$ is represented by $\vec{e}_1\wedge\vec{e}_2$ where $(\vec{e}_1,\vec{e}_2)$ is a local tangent frame.
%We decided instead to represent the Gauss map by the \underbar{orthogonal} to $\vec{e}_1\wedge\vec{e}_2$  which is equal in our notations to $\vec{\Phi}\wedge\vec{n}_{\vec{\Phi}}$ - this is what is denoted  $\hat{\mathbf {\Phi}}$ in \cite{CaUr} if ${\mathbf \Phi}$ is the Gauss map - The main reason for this switch of notations is due to the fact that we are going to define the {\it Hamiltonian isotopy classes} via the existence of a global lift to $V_2({\R}^4)$ (see lemma~\ref{lm-lie-surf} below). Such a lift for the Gauss maps of an immersion in $S^3$ always exists by definition while considering $\vec{\Phi}\wedge \vec{n}_{\vec{\Phi}}$ while it never exists while considering instead $\vec{e}_1\wedge\vec{e}_2$ unless $\Sigma$ is a torus.} 
in \cite{CaUr} that 
\[
\quad \quad\forall\, x\in\Sigma\quad 4\,|C(x)|^2\le 1\quad \quad,
\]
and equality holds at $x$ if and only if $\vec{G}^{\,+}$ (or equivalently $\vec{G}^{\,-}$) is conformal at $x$ with respect to the metric in $\Sigma$ induced by the immersion $\vec{\mathfrak G}$. It is also proved in \cite{CaUr} that, if $\Sigma$ is connected then
\[
\quad \quad\forall\, x\in\Sigma\quad 4\, |C(x)|^2\equiv 1\quad 
\]
if and only if there exists an element ${\mathcal P}\in SO(3)$  such that $\vec{\mathfrak G}(\Sigma)$ is included into the ${\mathcal P}-$geodesic  lagrangian sphere in $S^2_+\times S^2_-$ given explicitly by
\[
S^2_{\mathcal P}:=\lf\{(x,-{\mathcal P}(x))\in S^2_+\times S^2_-\quad;\quad x\in S^2\rg\}\quad.
\]
From now on in this section $\Sigma$ is closed without boundary. We define the {\it degree} {deg}$(\vec{\mathfrak G})$ of a lagrangian immersion  $\vec{\mathfrak G}$ to be the integer given by
\[
\mbox{deg}(\vec{\mathfrak G}):=\frac{1}{4\pi}\int_{\Sigma}(\vec{G}^{\,+})^\ast\om_{S^2}= \frac{1}{4\pi}\int_{\Sigma}C(x)\ dvol_{\vec{\mathfrak G}}=\frac{1}{4\pi}\int_{\Sigma}(\vec{G}^{\,-})^\ast\om_{S^2}
\]
in other words this is the degree of the map  $\vec{G}^{\,+}$ which is also the degree of the map $\vec{G}^{\,-}$. 

\medskip

 We consider the map $\vec{\mathcal{G}}_{\pm}:=\sqrt{2^{-1}}\, (\vec{G}^{\,+}\pm\vec{G}^{\,-})\in \wedge^2{\R}^4$. We have $\vec{\mathcal{G}}_{\pm}\wedge\vec{\mathcal{G}}_\pm=0$. We are interested in Lagrangian immersions\footnote{The space of $C^l$ lagrangian maps from a surface $\Sigma$ into $G_2^+({\R}^4)=S^2_+\times S^2_-$ will be denoted $C^l_{lag}(\Sigma,G_2^+({\R}^4))$. By an abuse of notations we shall sometimes
 mix the corresponding maps $\vec{\mathfrak G}\in S^2_+\times S^2_-$ and  $\vec{\mathcal G}_+\in \wedge^2{\R}^4$ since one can be recovered from the other in a tautological way due to the fact that the space of self dual  2-vectors $(\wedge^2{\R}^4)_+$ is orthogonal to the space of anti-self-dual 2-vectors $(\wedge^2{\R}^4)_-$.}
 $\vec{\mathfrak G}$ such that  $\vec{\mathcal G}_+$ admits a lifting in the Stiefel $S^1-$bundle $V_2({\R}^4)\simeq S^3\times S^2$ over $Gr^+_2({\R}^4)$ of orthonormal 2-frames in ${\R}^4$. The condition for admitting such a lift can easily be described. The K\"ahler form $\om_{S^2_+\times S^2_-}$ is the curvature for the canonical connection $\nabla_0$ on the tautological bundle $V_2({\R}^4)$ over $Gr^+_2({\R}^4)$. In order to lift a map $\vec{\mathfrak G}$ one could  first take the pull-back bundle $\vec{\mathfrak G}^{\,-1}V_2({\R}^4)$ over $\Sigma$ as well as the pullback-connection $\vec{\mathfrak G}^{\,-1}\nabla_0$. Since the immersion is assumed to be Lagrangian, the pull-back connection is flat. Then one could start from one point and consider the lift given by the parallel transport with respect to this flat connection along paths in $\Sigma$. As long as two such paths are defining a zero homology loop,
 since the connection is flat, the parallel transport gives the same lift. The problem comes when two such paths do not define a zero homology loop. In that case, in order for this operation to be uni-valued and in order to garantee
 the existence of a lift we will assume that the parallel transport along any generator of the $\pi_1(\Sigma)$ realizes a closed loop. In other words, we assume that
 \be
 \label{lift}
 \forall\, \Gamma\in \mbox{Loop}(\Sigma) \quad \forall D\subset S^2_+\times S^2_-\quad\mbox{ s. t. }\quad\p D=\vec{\mathfrak G}_\ast[\Gamma] \quad\int_D\om_{S^2_+\times S^2_-}\in 4\pi{\Z}
 \ee
This is for instance not the case for the {\it Clifford torus} $S^1\times S^1=\{ (x,y)\in S^2_+\times S^2_-\quad x_3=y_3=0\}$ while it's double cover satisfies (\ref{lift}) and admits a lift in $V_2({\R}^4)$ which is going to play an important role in the last part of the present section.

Examples of Lagrangian immersions admitting a lift and satisfying (\ref{lift}) are given by the {\it Gauss Maps } of immersions into $S^3$ (see \cite{Law,Pal,CaUr}) : Let $\vec{\Phi}$ be an immersion of an orientable  surface $\Sigma$ into $S^3$ and denote by $\vec{n}$ it's unit normal vector within the tangent space to $S^3$. We take    $\vec{\mathcal G}_{+}:=\vec{\Phi}\wedge \vec{n}$. The associated map $\vec{\mathfrak G}_{\vec{\Phi}}$ defines a {\it Lagrangian immersion} for the K\"ahler form given by $\om_{S^2\times S^2}$. Indeed, the symmetry of the {\it second fundamental form} gives $d\vec{\Phi}\dot{\wedge}d\vec{n}=0$ which is the lagrangian condition as given by lemma~\ref{lm-lag-cond}. In the particular case of a lagrangian immersion issued from a Gauss map of an immersion into $S^3$ the degree has the following topological interpretation.
\begin{Lm}
\label{lm-degree-gauss} {\bf [Degree of the Gauss Map.]}
Let $\Sigma$ be a closed oriented surface and let $\vec{\Phi}$ be an immersion of $\Sigma$ into $S^3$ then
\be
\label{deg-gauss}
\mbox{deg}\,(\vec{\mathfrak G}_{\vec{\Phi}})=\frac{1}{4\pi}\int_\Sigma\lf((\vec{\Phi}\wedge\vec{n}_{\vec{\Phi}})^+\rg)^\ast\om_{S^2}=\frac{1}{4\pi}\int_\Sigma C_{\vec{\Phi}}(x) \ dvol_{\vec{\mathfrak G}_{\vec{\Phi}}}= 1-g(\Sigma)\quad,
\ee
where $g(\Sigma)$ is the genus of $\Sigma$.
\hfill $\Box$
\end{Lm}
\noindent{\bf Proof of lemma~\ref{lm-degree-gauss}.} From (\ref{degree}) we have (omitting the subscript $\vec{\Phi}$ when there is no ambiguity)
\[
\frac{1}{4\pi}\int_\Sigma\lf((\vec{\Phi}\wedge\vec{n}_{\vec{\Phi}})^+\rg)^\ast\om_{S^2}=\frac{1}{8\pi}\int_\Sigma\star\ \vec{\Phi}\wedge\vec{n}\wedge\lf(d\vec{\Phi}\wedge d\vec{\Phi}+d\vec{n}\wedge d\vec{n}\rg)
\]
We have in one hand
\be
\label{deg-1}
\frac{1}{8\pi}\int_\Sigma\star\ \vec{\Phi}\wedge\vec{n}\wedge\lf(d\vec{\Phi}\wedge d\vec{\Phi}\rg)=\frac{1}{4\pi}\int_\Sigma\star\ \vec{\Phi}\wedge\vec{n}\wedge\lf(\p_{x_1}\vec{\Phi}\wedge \p_{x_2}\vec{\Phi}\rg)\ dx_1\wedge dx_2=\frac{1}{4\pi}\int_\Sigma\ dvol_{g_{\vec{\Phi}}}\quad.
\ee
In the other hand we have, using a local positive orthonormal frame $(\vec{e}_1,\vec{e}_2)$ of the tangent space to the immersion
\be
\label{deg-2}
\begin{array}{l}
\ds\frac{1}{8\pi}\int_\Sigma\star\ \vec{\Phi}\wedge\vec{n}\wedge\lf(d\vec{n}\wedge d\vec{n}\rg)=\frac{1}{4\pi}\int_\Sigma\p_{x_1}\vec{n}\cdot\vec{e}_1\ \p_{x_2}\vec{n}\cdot\vec{e}_2-(\p_{x_1}\vec{n}\cdot\vec{e}_2)^2\ dx_1\wedge dx_2\\[5mm]
\ds=\frac{1}{4\pi}\int_{\Sigma} K_{ext}\ dvol_{g_{\vec{\Phi}}}=\frac{1}{4\pi}\int_{\Sigma} K_{int}\ dvol_{g_{\vec{\Phi}}}-\frac{1}{4\pi}\int_{\Sigma} \ dvol_{g_{\vec{\Phi}}}
\end{array}
\ee
where $K_{ext}$ and $K_{int}$ are respectively the {\it extrinsic} and {\it intrinsic Gauss  curvatures}. We obtain (\ref{deg-gauss}) by summing (\ref{deg-1}) and (\ref{deg-2}) and applying Gauss Bonnet theorem. This concludes the proof of lemma~\ref{lm-degree-gauss}.\hfill $\Box$

\medskip

The following lemma is establishing some kind of reciproque of the above statement according to which the Gauss map of an immersion into $S^3$ defines a Lagrangian immersion into $G_2^+({\R}^4)$. Precisely we prove the existence of an underlying {\it Lie surface} to any Lagrangian immersion satisfying (\ref{lift}) (This is the global counterpart to proposition 3.2 \cite{Pal})
\begin{Lm}
\label{lm-lie-surf}
Let $\vec{\mathfrak G}$ be a smooth Lagrangian immersion into  $G_2^+({\R}^4)$. Assume $\vec{\mathfrak G}$ satisfies (\ref{lift}), then there exists a smooth pair of maps  $(\vec{a},\vec{b})$ from $\Sigma$ into $S^3$ such that
\[
\vec{\mathfrak G}=\sqrt{2^{-1}}\,\lf(\vec{a}\wedge \vec{b}+\star (\vec{a}\wedge \vec{b}),\vec{a}\wedge \vec{b}-\star (\vec{a}\wedge \vec{b})\rg) \quad\mbox{ such that }\quad \vec{a}\cdot\vec{b}=0
\]
and satisfying 
\be
\label{lift-1}
d\vec{a}\cdot\vec{b}=0\quad.
\ee
The choice of $(\vec{a},\vec{b})$ is unique modulo a global rotation constant on each connected component of $\Sigma$ : any other choice is given by $(\vec{a}^{\,\theta},\vec{b}^{\,\theta})$ such that $\vec{a}^{\,\theta}+\,i\,\vec{b}^{\,\theta}=e^{i\theta} (\vec{a}+i\,\vec{b})$ where $\theta$ is constant on any connected component of $\Sigma$.
\hfill $\Box$
\end{Lm}
\noindent{\bf Proof of Lemma~\ref{lm-lie-surf}.} The assumption (\ref{lift}) implies that  the $S^1-$bundle $\vec{\mathfrak G}^{\,-1}V_2({\R}^4)$ (pull back by $\vec{\mathfrak G}$ of the tautological bundle $V_2({\R}^4)\rightarrow G^+_2({\R}^4)$) is trivial. We take a global trivialization represented by a pair of unit vectors in ${\R}^4$ orthogonal to each other $(\vec{\al},\vec{\beta})$ and such that $\vec{\al}\wedge\vec{\beta}=\vec{\mathcal G}_+$.
The pull-back connection $\vec{\mathfrak G}^\ast\nabla^0$ which happens to be flat (since it's curvature is the pull-back by the Lagrangian immersion of the K\"ahler form) is represented in this trivialization by the following $i\,{\R}$ valued  closed form
\[
A:=i\,\vec{\al}\cdot d\vec{\beta}
\]
By (\ref{lift}) the circulation of this form on any closed loop is in $2\pi\,{\Z}$. Then there exists an $S^1$ valued map $g\ :\ \Sigma\rightarrow S^1$ whose degree on any loop $\Gamma$ is given by $(2\pi\,i)^{-1}\int_\Gamma\ g^{-1}\, dg$ and satisfying
\[
g^{-1}\, dg=i\,\vec{\al}\cdot d\vec{\beta}
\]
This gives that $d( g(\vec{\al}+i\vec{\beta}))=0$ and the pair of unit orthogonal vectors $(\vec{a},\vec{b})$ such that $(\vec{a}+i\vec{b})=g(\vec{\al}+i\vec{\beta})$ is a solution of the lemma.\hfill $\Box$

\medskip

Finally we have the following lemma

\begin{Lm}
\label{lm-gauss-imm}
The space of Gauss maps of immersions of a given closed oriented surface $\Sigma$ in $S^3$ realizes an open subspace of the space of Lagrangian immersions satisfying (\ref{lift}) which is itself an open subset of the space of Lagrangian maps from $\Sigma$ into $S^2_+\times S^2_-$ for the Frechet $C^\infty$ topology.\hfill $\Box$
\end{Lm} 
\noindent{\bf Proof of lemma~\ref{lm-gauss-imm}.} Let $\vec{\Phi}$ be an immersion from $\Sigma$ into $S^3$ and denote by $\vec{n}_{\vec{\Phi}}$ it's associated unit normal. Let $(\vec{\mathcal G}_{\vec{\Phi}})_+:=\vec{\Phi}\wedge\vec{n}_{\vec{\Phi}}$ be the associated Gauss map (we shall now omit the subscript $+$). We have in local coordinates (that we take to be conformal for $\vec{\Phi}$)
\[
\p_{x_i}(\vec{\Phi}\wedge\vec{n}_{\vec{\Phi}})=\p_{x_i}\vec{\Phi}\wedge\vec{n}_{\vec{\Phi}}+\vec{\Phi}\wedge\p_{x_i}\vec{n}_{\vec{\Phi}}
\]
Let $(\lambda_1,\lambda_2)$ such that
\[
\sum_{i=1}^2\la_i\, \p_{x_i}(\vec{\Phi}\wedge\vec{n}_{\vec{\Phi}})=0
\]
The scalar product of this identity respectively with $\p_{x_1}\vec{\Phi}\wedge\vec{n}_{\vec{\Phi}}$ and $\p_{x_2}\vec{\Phi}\wedge\vec{n}_{\vec{\Phi}}$ gives respectively
\[
\la_1\, |\p_{x_1}\vec{\Phi}|^2=0\quad\mbox{ and }\quad\la_2\, |\p_{x_2}\vec{\Phi}|^2=0\
\]
This implies that $(\la_1,\la_2)=(0,0)$ and we deduce that rank$\{\p_{x_1}(\vec{\Phi}\wedge\vec{n}_{\vec{\Phi}}),\p_{x_2}(\vec{\Phi}\wedge\vec{n}_{\vec{\Phi}})\}=2$ which implies that $\vec{\mathcal G}_{\vec{\Phi}}$ is an immersion.

We claim that for any $\ep>0$  there exists $\delta>0$ such that for any $\vec{\mathcal G}\in C^1_{lag}(\Sigma,G^+_2({\R}^4)$ satisfying (\ref{lift}) and 
\be
\label{lag-rep}
\|\vec{\mathcal G}-\vec{\mathcal G}_{\vec{\Phi}}\|_{C^1(\Sigma)}<\delta\quad\Longrightarrow\quad\exists\ \vec{w}\in \Gamma_\ep(\vec{\Phi}^{-1}TS^3) \quad\mbox{ s. t. }\quad\vec{\mathcal G}=\vec{\Phi}(\vec{w})\wedge \vec{n}_{\vec{\Phi}(\vec{w})}\quad,
\ee
where $\Gamma_\ep(\vec{\Phi}^{-1}TS^3)$ denotes the sections $\vec{w}$ of the pull-back bundle by $\vec{\Phi}$ of the tangent bundle $TS^3$ and satisfying\footnote{We take the $C^{0,\al}$ norm w.r.t. the metric induced by $\vec{\Phi}$ for some $0<\al<1$. We could reinforce the norm by taking the $C^1_{loc}-$norm away from the umbilic points of $\vec{\Phi}$ but since we are stating our result in the Frechet $C^\infty$ topology having $C^1$-closeness implies $C^1$-closeness is not required.} $\|\vec{w}\|_{C^{0,\al}(\Sigma)}<\ep$ and $\vec{\Phi}(\vec{w}):=\vec{\Phi}+\vec{w}/|\vec{\Phi}+\vec{w}|$. Assume by contradiction that (\ref{lag-rep}) does not hold. Then we would find $\vec{\mathcal G}_k$ converging in $C^1$ norm towards $\vec{\mathcal G}_{\vec{\Phi}}$, according to lemma~\ref{lm-lie-surf} there exists a global lift $\vec{\mathcal G}_k=\vec{a}_k\wedge\vec{b}_k$ such that $d\vec{a}_k\cdot\vec{b}_k=0$ but no representative
$\cos\,t\ \vec{a}_k+\sin\, t\ \vec{b}_k$ is close to $\vec{\Phi}$ in $C^{0,\al}$ norm. We have
\[
d\vec{\mathcal G}_k= d\vec{a}_k\wedge\vec{b}_k+\vec{a}_k\wedge d\vec{b}_k\longrightarrow d\vec{\mathcal G}_{\vec{\Phi}}\quad\quad\mbox{ in }C^0(\Sigma)\quad.
\]
Since $|d\vec{\mathcal G}_k|^2=|d\vec{a}_k|^2+|d\vec{b}_k|^2$ (where we crucially use the fact that $d\vec{a}_k\cdot\vec{b}_k=0$ which is also equivalent to $d\vec{b}_k\cdot\vec{a}_k=0$), we have that modulo a subsequence that
we also denote $(\vec{a}_k,\vec{b}_k)$ we have that 
\[
d\vec{a}_k\rightharpoonup d\vec{a}_\infty\quad\mbox{ and }\quad d\vec{b}_k\rightharpoonup d\vec{b}_\infty \quad\mbox{ weakly in }(L^\infty(\Sigma))^\ast
\]
Hence $\vec{a}_k$ (resp. $\vec{b}_k$) converge strongly in $C^{0,\al}(\Sigma)$ for any $\al<1$ towards $\vec{a}_\infty$ (resp. $\vec{b}_\infty$) and we have
\[
\vec{a}_\infty,\vec{b}_\infty\in W^{1,\infty}(\Sigma,S^3)\quad,\quad\vec{a}_\infty\cdot\vec{b}_\infty=0\quad,\quad d\vec{a}_\infty\cdot\vec{b}_\infty=0\quad\mbox{ and }\quad \vec{a}_\infty\wedge\vec{b}_\infty=\vec{\mathcal G}_{\vec{\Phi}}\quad.
\]
Then, according to lemma~\ref{lm-lie-surf}  there exists $t\in [0,2\pi)$ such that   $\cos\,t\ \vec{a}_\infty+\sin\, t\ \vec{b}_\infty=\vec{\Phi}$ and $-\sin\,t\ \vec{a}_\infty+\cos\, t\ \vec{b}_\infty=\vec{n}_{\vec{\Phi}}$. We can always modify the sequence $(\vec{a}_k,\vec{b}_k)$ by a rotation in such a way that
\[
\vec{a}_k\rightharpoonup\vec{\Phi}\quad\mbox{ and }\quad\vec{b}_k\rightharpoonup\vec{n}_{\vec{\Phi}}\quad\mbox{ weakly in }(W^{1,\infty}(\Sigma))^\ast
\]
This gives a contradiction. Hence we have deduced that (\ref{lag-rep}) holds true. Observe that if $\vec{\mathcal G}\in C^l$ for any $l\ge 1$, the procedures described above converting  the condition (\ref{lift}) into the existence of a Lagrangian lift imply  that $\vec{\Phi}(\vec{w})$ as well as $\vec{w}$ are both in $C^l$. This permits to bootstrap and deduce lemma~\ref{lm-gauss-imm} from (\ref{lag-rep}).\hfill $\Box$

\subsection{The two first steps of the hierarchy for Lagrangian immersions in $G^+_2({\R}^4)$.}

 Let $\Sigma$ be a closed oriented surface. We shall denote by $\mbox{Imm}^{4,2}_{\,0}(\Sigma, S^3)$ the space of $W^{4,2}$ immersions of a closed oriented 2-manifold $\Sigma$ into $S^3$ which are regular homotopic to a minimal immersion. A result of Pinkall, \cite{Pink}, gives that for any genus $g(\Sigma)>0$ there are exactly 2 classes of regular homotopic immersions (i.e. there are exactly 2 connected components  in the class of $C^1$-immersions) of the closed oriented surface of genus $g$ into $S^3$ while the space of immersions of $S^2$ into ${\R}^3$  is path connected (and hence into $S^3$ as well) due to a famous result by Smale \cite{Sm} . From the result of Lawson \cite{Law} there is at least one of the two classes which contains a minimal embedding. The main reason why we are restricting to the isotopy classes possessing a minimal immersion is first because we are interested in the area of minimal embedded surfaces in $S^3$ (with area strictly less than $8\pi$ in fact) and the second main reason is due to the following lemma. 
 \begin{Lm}
 \label{reg-iso}{\bf [Hamiltonian Stationary Gauss Maps.]}
 Let $\vec{\Phi}$ in the class $\mbox{Imm}^{4,2}_{\,0}(\Sigma, S^3)$ such that the associated Gauss map is a critical point of the area under variations of $\vec{\Phi}$. Then on each component of $\Sigma$ there exists $t\in [0,2\pi]$ such that
 $\vec{\Phi}_t=\cos t\, \vec{\Phi}+\sin t \ \vec{n}_{\vec{\Phi}}$ defines a possibly branched minimal immersion into $S^3$ as well as $\vec{\Phi}_{t+\pi/2}$. Moreover on each component there holds
 \be
 \label{id-aire}
 \mbox{Area}(\vec{G}_{\vec{\Phi}})=2\, \mbox{Area}({\vec{\Phi}}_t)+2\, \mbox{Area}({\vec{\Phi}}_{t+\pi/2})\quad.
 \ee
 \hfill $\Box$
 \end{Lm}
\noindent{\bf Proof of lemma~\ref{reg-iso}.}  We denote by $\vec{H}_{\vec{\mathfrak G}_{\vec{\Phi}}}$ the mean curvature vector of the associated Gauss map $\vec{\mathfrak G}_{\vec{\Phi}}$ in $G_2^+({\R}^4)$. Consider it's contraction\footnote{This contraction is the so called {\it Maslov form}.} with the K\"ahler form $\om_{S^2\times S^2}$ along the lagrangian immersion given by $\vec{\mathfrak G}_{\vec{\Phi}}$ :   
 \[
\forall \ X\in T_x\Sigma\quad\quad \al_{\vec{\Phi}}(X):=\frac{1}{\pi}\,\om_{S^2\times S^2}\lf(d\vec{\mathfrak G}_{\vec{\Phi}}(X),\vec{H}_{\vec{\mathfrak G}_{\vec{\Phi}}}(x)\rg)
 \]
 A computation discovered by Dazor \cite{Daz} gives
 \be
 \label{daz}
 d\al_{\vec{\Phi}}=  (\vec{\mathfrak G}_{\vec{\Phi}})^\ast\mbox{Ric}
 \ee
where Ric denotes the Ricci 2-form on the K\"ahler manifold $G_2^+({\R}^4)$. Since the manifold is {\it K\"ahler-Einstein}, the Ricci form is proportional to the K\"ahler form and since $\vec{\mathfrak G}_{\vec{\Phi}}$ is lagrangian (\ref{daz}) implies that $\al_{\vec{\Phi}}$ is closed and defined a {\it real cohomology class} on $\Sigma$. This class is invariant under {\it hamiltonian isotopies}  (see \cite{Oh3}). Regular homotopies at the level of $\vec{\Phi}$ generate obviously {\it hamiltonian isotopies} at the level of $\vec{\mathfrak G}_{\vec{\Phi}}$.  For a minimal immersion  $\vec{\xi}$, $\vec{\mathfrak G}_{\vec{\xi}}$ is also minimal (\cite{Law}) that is $\vec{H}_{\vec{\mathfrak G}_{\vec{\xi}}}\equiv 0$ hence $\al_{\vec{\xi}}=0$. Then for any immersion $\vec{\Phi}$ in the same component as $\vec{\xi}$ the class $\al_{\vec{\Phi}}$ is trivial.

Assuming now that $\vec{\Phi}\in \mbox{Imm}^{4,2}_{\,0}(\Sigma, S^3)$ and that $\vec{\mathfrak G}_{\vec{\Phi}}$ is a critical point of the area under variations of $\vec{\Phi}$. Because of lemma~\ref{lm-gauss-imm} the Gauss map $\vec{\mathfrak G}_{\vec{\Phi}}$ is a critical point under variations which preserve the condition (\ref{lift}). In particular it is a critical point of the area under {\it Hamiltonian perturbations} i.e. $\vec{\mathfrak G}_{\vec{\Phi}}$ is {\it Hamiltonian stationary}. This implies that $\al_{\vec{\Phi}}$ is co-closed for the metric induced by $\vec{\mathfrak G}_{\vec{\Phi}}$. But $\al_{\vec{\Phi}}$ is also exact since the associated cohomology class is zero. Hence $\al_{\vec{\Phi}}$ is zero and $\vec{\mathfrak G}_{\vec{\Phi}}$ is also critical for Lagrangian perturbations (i.e. $\vec{\mathfrak G}_{\vec{\Phi}}$ is {\it lagrangian stationary}). It is proved in \cite{SW} that  branched lagrangian immersions which are
{\it lagrangian stationary} in a K\"ahler Einstein manifolds are minimal. Hence we have $\vec{H}_{\vec{\mathfrak G}_{\vec{\Phi}}}=0$ on $\Sigma$. 

We assume $\Sigma$ is connected. Using proposition 5 of \cite{CaUr} we deduce that either $4\, |C|\equiv 1$ on $\Sigma$  or the zeros of $4\,|C|(x)-1$ are isolated. In the first case $\vec{\mathfrak G}_{\vec{\Phi}}(\Sigma)\subset S^2_{\mathcal P}$ for some ${\mathcal P}$ in $SO(3)$. Since $\vec{\Phi}$ is assumed to be an immersion, this is also the case for $\vec{\mathfrak G}_{\vec{\Phi}}$, Such an immersion would realize a covering of $S^2$ by $\Sigma$. Such a covering induces an injection on the first homotopy groups and then $\Sigma\simeq S^2$ moreover the number of pre-image by such a covering, which is constant, has to be one since the path connecting two distinct preimages would be pushed  to a non trivial element of $\pi_1(S^2)=0$. Hence,   modulo the action of a diffeomorphism of $S^2$,  $\vec{\mathfrak G}_{\vec{\Phi}}(x)=(x,{\mathcal P}(x))$. From the previous subsection we know that $x\rightarrow (x,{\mathcal P}(x))$ is the Gauss map of  the immersion of a geodesic two sphere of $S^3$. Using now lemma~\ref{lm-lie-surf} we deduce the lemma in this case.

 In the second case, when the zeros of $x\in \Sigma\rightarrow\, 4\, |C(x)|-1$ are isolated, away from the zeros, from theorem 3 of \cite{CaUr}, there exists ${\mathfrak I}\in O(S^2_+\times S^2_-)$ and a minimal immersion $\vec{\xi}$ into $S^3$ such that locally 
 \[
 \vec{\mathfrak G}_{\vec{\Phi}}={\mathfrak I}\circ \vec{\mathfrak G}_{\vec{\xi}}
 \]
 Using lemma~\ref{lm-iso-s2s2}, there exists then locally, away from the zeros of $4\, |C(x)|-1$, an isometry ${\mathcal I}$ in $O(4)$ and $i,j\in{1,2}$ such that
 \[
 \lf((\vec{\Phi}\wedge\vec{n}_{\vec{\Phi}})^+,(\vec{\Phi}\wedge\vec{n}_{\vec{\Phi}})^-\rg)=(-1)^i\lf(({\mathcal I}\circ\vec{\xi}\wedge\vec{n}_{{\mathcal I}\circ\vec{\xi}})^+, (-1)^j ({\mathcal I}\circ\vec{\xi}\wedge\vec{n}_{{\mathcal I}\circ\vec{\xi}})^-\rg)
 \]
 Observe that $j=2$ since we have respectively $(\vec{G}^+_{\vec{\Phi}})^\ast\om_{S^2}=(\vec{G}^-_{\vec{\Phi}})^\ast\om_{S^2}$ and $(\vec{G}^+_{{\mathcal I}\circ\vec{\xi}})^\ast\om_{S^2}=(\vec{G}^-_{{\mathcal I}\circ\vec{\xi}})^\ast\om_{S^2}$. Obviously ${{\mathcal I}\circ\vec{\xi}}$ is minimal as well as $\vec{n}_{{\mathcal I}\circ\vec{\xi}}$ thanks to a classical result in \cite{Law}. Thus we have the local existence of  a minimal immersion $\vec{\sigma}$ such that locally
 \[
 \vec{\mathfrak G}_{\vec{\Phi}}=(-1)^i\, \vec{\mathfrak G}_{\vec{\sigma}}\quad.
 \]
Using lemma~\ref{lm-lie-surf}  we deduce the existence there exists locally of $t\in [0,2\pi)$ such that $\vec{\Phi}_t:=\cos t\, \vec{\Phi}+\sin t \ \vec{n}_{\vec{\Phi}}$ is a minimal immersion. The set of $t\in [0,2\pi)$ such that $\vec{\Phi}_t$ is minimal known to be is discrete. Hence a simple continuity argument implies that $t$ is constant on each component of $\Sigma$ and lemma~\ref{reg-iso} is proved.\hfill $\Box$

\medskip 
 
 We denote by $\mbox{Imm}^{3,2}_{Lag}(\Sigma, G_2^+({\R}^4))$ the space of Sobolev $W^{3,2}$ lagrangian immersions from $\Sigma$ into $G_2^+({\R}^4)$. Following \cite{Riv-Ind} we introduce $\mbox{Diff}_+^\ast(\Sigma)$ the subspace of positive $W^{4,2}-$ diffeomorphisms isotopic to the identity and fixing $3-$points if $g(\Sigma)=0$, 1 point if $g(\Sigma)=1$ and no point if $g(\Sigma)>1$. The configuration space for the hierarchy that we are going to consider is given by
\[
{\mathfrak M}(\Sigma):=\mbox{Imm}^{4,2}_{\, 0}(\Sigma, S^3)/\mbox{Diff}_+^\ast(\Sigma)
\]
and
\[
{\mathcal M}(\Sigma):=\lf\{\vec{\mathfrak G}\in \mbox{Imm}^{3,2}_{Lag}(\Sigma, G_2^+({\R}^4))/\mbox{Diff}_+^\ast(\Sigma)\quad;\quad (\ref{lift})\mbox{ holds }\rg\}\quad.
\]
It is proved in \cite{Riv-Lag-Visc} that ${\mathfrak M}(\Sigma)$ defines a Hilbert manifold complete for the {\it Palais distance}. And we denote by ${\mathfrak M}_b$ (resp. ${\mathcal M}_b$) to be the disjoint union of the ${\mathfrak M}(\Sigma)$ (resp. ${\mathcal M}_b$) where $b(\Sigma)$ is the total {\it Betti number} of $\Sigma$, $b(\Sigma)=b_0(\Sigma)+b_1(\Sigma)+b_2(\Sigma)\le b$ and ${\mathfrak M}$ (resp. ${\mathcal M}$) denotes the unions of the ${\mathfrak M}_b$ (resp. ${\mathcal M}_b$) for any $b\in {\N}$. On the space ${\mathcal M}$ of lagrangian immersions we shall be considering the following energy
\[
A(\vec{\mathfrak G}):=\int_\Sigma dvol_{\vec{\mathfrak G}}+8\,\pi\ \mbox{deg}(\vec{\mathfrak G})
\]
We have the following lemma
\begin{Lm}
\label{lm-lagran} We have
\be
\label{lagran}
\min_{\vec{\mathfrak G}\in {\mathcal M}}\ A(\vec{\mathfrak G})=0
\ee
and $A(\vec{\mathfrak G})=0$  exactly when $\vec{\mathfrak G}:=-\,\vec{\mathfrak G}_{\vec{\Phi}}$ where $\vec{\mathfrak G}_{\vec{\Phi}}$ is the Gauss map of  a multiple cover of a geodesic 2-sphere in $S^3$.
\hfill$\Box$
\end{Lm}
\noindent{\bf Proof of lemma~\ref{lm-lagran}.} Observe that, denoting $|\cdot|$ the metric induced by $\vec{\mathfrak G}$, we have
\[
\frac{1}{4\pi}\lf|(\vec{G}^+)^\ast\om_{S^2}\rg|=\frac{1}{8\pi}\lf(|(\vec{G}^+)^\ast\om_{S^2}|+|(\vec{G}^-)^\ast\om_{S^2}|\rg)\le \frac{1}{8\pi}\lf|\vec{\mathfrak G}^\ast(\pi_+^\ast\om_{S^2_+}+\pi_-^\ast\om_{S^2_-})\rg|
\] 
We claim that the {\it co-mass} of $\pi_+^\ast\om_{S^2_+}+\pi_-^\ast\om_{S^2_-}$ is 1. Indeed we have
\[
\|\pi_+^\ast\om_{S^2_+}+\pi_-^\ast\om_{S^2_-}\|_\ast:=\max_{y=(y_+,y_-)\in G_2^+{\R}^4= S^2_+\times S^2_-}\ \max_{\vec{E},\vec{F}\in U_yG_2^+{\R}^4}<\pi_+^\ast\om_{S^2_+}+\pi_-^\ast\om_{S^2_-},\vec{E}\wedge\vec{F}>
\]
where $UG_2^+{\R}^4$  is the unit sphere bundle in $TG_2^+{\R}^4$. Denoting $\vec{E}_\pm=\pi^\pm_\ast \vec{E}$ and $\vec{F}_\pm=\pi^\pm_\ast \vec{F}$, we have
\[
\begin{array}{l}
\ds<\pi_+^\ast\om_{S^2_+}+\pi_-^\ast\om_{S^2_-},\vec{E}\wedge\vec{F}>=\om_{S^2}(\vec{E}^+,\vec{F}^+)+\om_{S^2}(\vec{E}^-,\vec{F}^-)\\[5mm]
\ds\le |\vec{E}_+\wedge\vec{F}_+|+|\vec{E}_-\wedge\vec{F}_-|\le |\vec{E}_+|\,|\vec{F}_+|+|\vec{E}_-|\,|\vec{F}_-|\le 2^{-1}\,(|\vec{E}_+|^2+|\vec{F}_+|^2+|\vec{E}_-|^2+|\vec{F}_-|^2)\le 1
\end{array}
\]
This proves the claim and we finally obtain
\[
- \frac{1}{4\pi}(\vec{G}^+)^\ast\om_{S^2}\le dvol_{\vec{\mathfrak G}}
\]
which gives  $A\ge 0$ on the space of Lagrangian immersions into $G_2^+({\R}^4)$. Observe that $A(\vec{\mathfrak G})=0$ if and only if the above inequalities are equalities and
\[
- \frac{1}{4\pi}(\vec{G}^+)^\ast\om_{S^2}=-\frac{1}{8\pi}\vec{\mathfrak G}^\ast(\pi_+^\ast\om_{S^2_+}+\pi_-^\ast\om_{S^2_-})= \frac{1}{8\pi} dvol_{\vec{\mathfrak G}}
\]
which means that $\vec{\mathfrak G}$ is calibrated by $-\pi_+^\ast\om_{S^2_+}+\pi_-^\ast\om_{S^2_-}$. This implies in particular that $4\,|C|^2\equiv 1$ on $\Sigma$ and this gives that $\vec{\mathfrak G}(\Sigma)=S^2_{\mathcal P}$ for some ${\mathcal P}\in SO(3)$ (see proposition 2 of \cite{CaUr}). We shall see in the next subsection that the geodesic immersion $\vec{\Phi}$ : $S^2\rightarrow S^3$ whose image is the oriented geodesic 2-sphere orthogonal and oriented by the unit vector $\vec{g}\in S^3$ generates a Gauss map $\vec{\mathfrak G}_{\vec{\Phi}}$ which is given (modulo reparametrisation) by $y\rightarrow (y, {\mathcal P}(\vec{g})(y))$ where the rotation ${\mathcal P}(\vec{g})$
is given  by ${\mathbf y}\rightarrow {\mathbf g}^\ast{\mathbf y}{\mathbf g}$ where $\mathbf y$ and $\mathbf g$ are the quaternionic representatives of respectively $\vec{y}\in S^3$ and $\vec{g}\in S^3$. 
This concludes the proof of lemma~\ref{lm-lagran}.\hfill $\Box$ 

\medskip

We shall denote by ${\mathbf F}$ the distance
%\footnote{ Obviously ${\mathfrak M}$ is not complete for ${\mathbf F}$ while it is complete for the Palais distance (see \cite{Riv-minmax}). Moreover the Palais distance defines clearly a much finer topology than the ${\mathbf F}-$ distance. Hence the continuity with respect to ${\mathbf F}$ required in the definition of the admissible families will be obviously preserved by the pseudo-gradient flow.  } 
obtained by summing the flat distance for currents with the varifold distance. Precisely, for any immersion $\vec{\mathfrak G}\in {\mathfrak M}(\Sigma)$ we can define the corresponding {\it oriented varifold} in $G_2({\R}^4)=S^2_+\times S^2_-$ as follows
\[
\forall \varphi\in C^0({G}_2(S^2_+\times S^2_-))\quad V_{\vec{\mathfrak G}}(\varphi):=\int_{\Sigma}\varphi\lf(\vec{\mathfrak G}(x),\vec{\mathfrak G}_\ast T_x\Sigma\rg)\ dvol_{g_{\vec{\mathfrak G}}}
\]
where ${G}_2(S^2_+\times S^2_-)$ is the Grassmann bundle of oriented 2 planes in $T(S^2_+\times S^2_-)$ over $S^2_+\times S^2_-$. When $\varphi$ is just a function in $S^2_+\times S^2_-$ we keep denoting
\[
V_{\vec{\mathfrak G}}(\varphi):=\int_{\Sigma}\varphi\lf(\vec{\mathfrak G}(x)\rg)\ dvol_{g_{\vec{\mathfrak G}}}
\]
We call the {\it ${\mathbf F}-$distance} between 2 immersions $\vec{\mathfrak G}$ and $\vec{\mathfrak H}$ of respectively 2 oriented closed surfaces $\Sigma$ and $\Sigma'$
\[
{\mathbf F}(\vec{\mathfrak G},\vec{\mathfrak H}):=\sup_{\|\varphi\|_{Lip}\le 1}V_{\vec{\mathfrak G}}(\varphi)-V_{\vec{\mathfrak H}}(\varphi)+{\mathcal F}\lf(\vec{\mathfrak G}_\ast[\Sigma]-\vec{\mathfrak H}_\ast[\Sigma']\rg)
\]
where ${\mathcal F}$ is the usual Flat norms between 2-cycles and $\vec{\mathfrak G}_\ast[\Sigma]$ and $\vec{\mathfrak H}_\ast[\Sigma']$ denote respectively
the push forwards  by $\vec{\mathfrak G}$ and by $\vec{\mathfrak H}$ of the currents of integration along respectively $\Sigma$ and  $\Sigma'$. Observe that $V$ and ${\mathbf F}$ are independent of the oriented parametrization and hence these two functions ``descend'' to ${\mathfrak M}$.

\medskip

If now one restricts to the Lagrangian immersions issued from Gauss maps the lower bound is increased. We have the following theorem proved in \cite{OrRi}.
\begin{Th}
\label{lm-gauss-low} {\bf [The ground states]} We have 
\be
\label{gauss-low}
\min_{\vec{\Phi}\in {\mathfrak M}}A(\vec{\mathfrak G}_{\vec{\Phi}})=16\pi
\ee
The minimum is exactly achieved by the immersions of $S^2$ covering a geodesic 2-sphere in  $S^3$ exactly once. Moreover for any $b\in {\N}^\ast$  there exists $\vec{\Phi}_k\in {\mathfrak M}_b$ such that
\[
\lim_{k\rightarrow +\infty}A(\vec{\mathfrak G}_{\vec{\Phi}_k})=16\pi
\]
such that $\Pi\circ\vec{\mathfrak G}_{\vec{\Phi}_k}$ ${\mathbf F}$-converges towards a finite number of $\Pi( S^2_{\mathcal P})$ spheres 
\footnote{Observe that only the positive spheres are Lagrangian maps issued from a Gauss map of an immersion into $S^3$. The space of Lagrangian maps issued from Gauss maps is clearly non closed under
${\mathbf F}$ convergence. We shall see concrete examples in the next subsection} and where $\Pi$ is denoting the canonical projection from $S^2\times S^2$ into ${\R}P^2\times{\R}P^2$ .\hfill $\Box$ 
\end{Th}

\medskip

\begin{Con}
\label{la-conjecture}
The only way to approach the ground state value $16\pi$ is the one given by theorem~\ref{gauss-low}.
\hfill $\Box$
\end{Con}

Hence using the notations of the previous sections the first stage of the hierarchy for $A\circ\vec{\mathfrak G}$ in ${\mathfrak M}$ is given by 
\[
{\mathcal W}_0=\inf_{\vec{\Phi}\in {\mathfrak M}}A(\vec{\mathfrak G}_{\vec{\Phi}})=16\pi\quad,
\]
and
\[
{\mathcal C}_0=\{S^2_{\mathcal P}\quad;\quad {\mathcal P}\in SO(3)\}\simeq SO(3)\quad.
\]
We are now exploiting the $H_3({\mathcal C}_0,{\Z})=\Z$ in order to move up to the second stage of  the hierarchy. 

\medskip

For any $b\in {\N}^\ast$ we define
%\footnote{Observe that  $\vec{\Lambda}\in \mbox{Sweep}_1({\mathfrak M})$ does not imply that
%\[
%\max_{z\in \p Z}A(\vec{\Lambda}(z))=0\quad.
%\]
%Indeed it is not excluded  that $\vec{\Lambda}(z)$ is the immersion of $S^2\sqcup S^2$. The first immersion is $\vec{\mathfrak G}_1\ :\ x\in S^2\rightarrow (x,x)\in {\mathcal C}_0$ while the second is $\vec{\mathfrak G}_1\ :\ x\in S^2\rightarrow (-x,-x)\in -\,{\mathcal C}_0$. Hence $\mbox{deg}(\vec{\Lambda}(z))=\mbox{deg}(\vec{\mathfrak G}_1)+\mbox{deg}(\vec{\mathfrak G}_2)=0$ but 
%\[
%\int_{S^2\sqcup S^2} dvol_{\vec{\Lambda}(z)}=\int_{S^2} dvol_{\vec{\mathfrak G}_1}+\int_{S^2} dvol_{\vec{\mathfrak G}_2}=16\pi
%\] }
\[
\mbox{Sweep}_1({\mathfrak M}_b):=\lf\{
\begin{array}{c}
\ds\vec{\Phi}\in C^0(Z,{\mathfrak M}_b)\quad\mbox{ where $Z$ is a finite simplicial 4-complex }\\[5mm]
\ds \exists \, \vec{\mathfrak G}\in C^0_{\mathbf F}(\ov{Z},{\mathcal M})\quad\quad\mbox{ s. t. }\quad\vec{\mathfrak G}=\vec{\mathfrak G}_{\vec{\Phi}}\quad\mbox{ in }Z\\[5mm]
%\ds \forall \, z\in \dot{Z}\quad g(\Sigma)\ge 1\quad\mbox{ where } \Lambda(z)\in \mbox{Imm}^{3,2}_{Lag}(\Sigma, G^+_2({\R}^4))\quad\\[5mm]
\ds A(\vec{\mathfrak G})\equiv 16\pi\quad\mbox{ on }\p Z\quad;\quad \vec{\mathfrak G}(\p Z)\subset {\Z}\cdot{\mathcal C}_0 \\[5mm]
\ds\vec{\mathfrak G}_\ast [\p Z]\ne 0\quad\mbox{ in } H_3(SO(3),{\Z})
\end{array}
\rg\}
\]
We shall see in the next subsection that 
\be
\label{nontriv}
\mbox{Sweep}_1({\mathfrak M}_4)\ne \emptyset\quad.
\ee
 We then define the corresponding {\it width} of this second stage of the hierarchy
\[
{\mathcal W}_1(b):=\inf_{\vec{\Phi}\,\in\, \mbox{Sweep}_1({\mathfrak M}_b)}\ \max_{z\in Z}A(\vec{\mathfrak G}_{\vec{\Phi}}(z))\quad.
\]
Finally we denote
\[
{\mathcal W}_1=\inf_{b\in {\N}}{\mathcal W}_1(b)
\]
The non triviality of the second stage of the hierarchy is given by the following lemma
\begin{Lm}
\label{lm-width}{\bf [Non triviality of the second step.]} Under the previous notations, for any $b\in {\N}$ such that $\mbox{Sweep}_1({\mathfrak M}_b)\ne \emptyset$, assume the conjecture~\ref{la-conjecture} is true we have
\be
\label{with-16-pi}
{\mathcal W}_1(b)>16\pi\quad.
\ee
\hfill $\Box$
\end{Lm}
\noindent{\bf Proof of lemma~\ref{lm-width}.} Assume ${\mathcal W}_1(b)=16\pi$. We consider a minimizing sequence $\vec{\Phi}_k\in C^0(Z,{\mathfrak M}_b)$ such that
\[
\lim_{k\rightarrow +\infty} \max_{z\in Z}A(\vec{\mathfrak G}_{\vec{\Phi}_k}(z))=16\pi\quad.
\]
Assuming conjecture~\ref{la-conjecture}  we have then
\[
\lim_{k\rightarrow +\infty} \max_{z\in Z}\lf|A(\vec{\mathfrak G}_{\vec{\Phi}_k}(z))-16\pi\rg|=0\quad.
\]
Using lemma~\ref{lm-degree-gauss} we have
\[
 \max_{z\in Z}\lf|\mbox{deg}(\vec{\mathfrak G}_{\vec{\Phi}_k}(z))\rg|\le b\quad.
\]
Hence there exists an $N\in {\N}$ such that
\be
\label{conv-with-16pi}
\lim_{k\rightarrow +\infty}\max_{z\in Z}\ {\mathbf F}(\vec{\mathfrak G}_{\vec{\Phi}_k}(z),{\mathcal C}_0^N)=0\quad\mbox{where }\quad{\mathcal C}_0^N:=\bigsqcup_{|i|\le N}i\,{\mathcal C}_0\quad.
\ee
Since ${\mathcal C}_0^N$ is a smooth finite dimensional closed sub-manifold of the Banach space of flat chains, we can construct in the ${\mathbf F}-$neighbourhood of ${\mathcal C}_0^N$ a continuous projection map $\Pi_N$ onto
${\mathcal C}_0^N$. For $k$ large enough the map
\[
z\in Z\ \longrightarrow\ \Pi_N(\vec{\mathfrak G}_{\vec{\Phi}_k}(z))\in {\mathcal C}_0^N
\]
is well defined, continuous and defines a cycle in ${\mathcal Z}(SO(3),{\Z})$, the space of integer singular chains in $SO(3)$, whose boundary equals $(\vec{\mathfrak G}_{\vec{\Phi}_k})_\ast [\p Z]\ne 0$ in $H_3(SO(3),{\Z})$ which is a contradiction.\hfill $\Box$

\subsection{Optimal Lagrangian Canonical Families for Minimal Embeddings.}
The goal of this section is to establish the following lemma.
\begin{Lm}
\label{lm-canonical family}{\bf [Optimal Lagrangian Canonical Families]}
Let $\vec{\Phi}$ be a minimal embedding from a closed oriented surface $\Sigma$ into $S^3$ which is not a geodesic 2-sphere or a multiple covering of a geodesic sphere. Then there exists a 4-dimensional simplex $Z$  and an element $\vec{\mathfrak G}\in C^0_{\mathbf F}(\ov{Z},{\mathcal M})$ such that $\vec{\mathfrak G}\in\mbox{Sweep}_1({\mathfrak M}_{b(\Sigma)})$, $\vec{\mathfrak G}_\ast [\p Z]$ realizes a $2\,g(\Sigma)$ multiple of a generator of $H_3(SO(3),{\Z})$. Moreover $\vec{\mathfrak G}_{\vec{\Phi}}\in \vec{\mathfrak G}(Z)$ and
\be
\label{linegalite}
 \max_{z\in Z}\mbox{Area}(\vec{\mathfrak G}(z))+8\pi\ \mbox{deg}(\vec{\mathfrak G}(z))=\mbox{Area}(\vec{\mathfrak G}_{\vec{\Phi}})+8\pi\ \mbox{deg}(\vec{\mathfrak G}_{\vec{\Phi}})=4\ \mbox{Area}(\vec{\Phi})\quad.
\ee
\hfill $\Box$
\end{Lm}

As a consequence of lemma~\ref{lm-canonical family} we get in particular the following corollary :
\begin{Co}
\label{co-lemma}
Under the previous notations we have 
\be
\label{non-trivial}
\mbox{Sweep}_1({\mathfrak M}_{4})\ne\emptyset\quad.
\ee
Moreover, for any minimal immersion $\vec{\Phi}$ of a surface $\Sigma$ whose total Betti number $b(\Sigma)$ is bounded by $b$ and which is not a geodesic $2-$sphere in $S^3$ one has
\be
\label{fund-ine-a}
{\mathcal W}_1(b)\le 4\ \mbox{Area}(\vec{\Phi})\quad.
\ee
In particular
\be
\label{fund-ine}
{\mathcal W}_1(b)\le 8\,\pi^2\quad.
\ee
\hfill $\Box$
\end{Co}
The last assertion of the corollary is obtained by taking $\vec{\Phi}$ to be the {\it Clifford Torus}. In the next sub-section we establish that (\ref{fund-ine}) is in fact an equality that will imply the Willmore conjecture.

\medskip

\noindent{\bf Proof of lemma~\ref{lm-canonical family}.}
Let $\Sigma$ be a closed oriented surface and $\vec{\Phi}\in \mbox{Imm}(\Sigma,S^3)$. Recall that $\vec{\mathfrak G}_{\vec{\Phi}}=\vec{G}^+_{\vec{\Phi}}+\vec{G}^-_{\vec{\Phi}} = \sqrt{2}\, \vec{\Phi}\wedge\vec{n}_{\vec{\Phi}}\in S^5_{\sqrt{2}}$ where we use the isometric immersion of $S^2_+\times S^2_-$ into the sphere of radius $\sqrt{2}$ in ${\R}^3\times {\R}^3$. We have in conformal coordinates for $\vec{\Phi}$
\[
\begin{array}{l}
\ds|\nabla\vec{\mathfrak G}_{\vec{\Phi}}|^2=|\p_{x_1}\vec{\mathfrak G}_{\vec{\Phi}}|^2+|\p_{x_2}\vec{\mathfrak G}_{\vec{\Phi}}|^2\\[5mm]
=2\, |\p_{x_1}\vec{\Phi}\wedge\vec{n}_{\vec{\Phi}}+\vec{\Phi}\wedge\p_{x_1}\vec{n}_{\vec{\Phi}}|^2+2\, |\p_{x_2}\vec{\Phi}\wedge\vec{n}+\vec{\Phi}\wedge\p_{x_2}\vec{n}|^2=2\,|\nabla\vec{\Phi}|^2+2\, |\nabla\vec{n}|^2=4\, e^{2\la}+2\, |\nabla\vec{n}|^2
\end{array}
\]
We have
\be
\label{can-1}
\begin{array}{l}
\ds\int_{\Sigma}\ dvol_{\vec{\mathfrak G}_{\vec{\Phi}}}= \int_{\Sigma}|\p_{x_1}\vec{\mathfrak G}_{\vec{\Phi}}\wedge \p_{x_2}\vec{\mathfrak G}_{\vec{\Phi}}|\ dx_{1}\wedge dx_2\\[5mm]
\ds\quad\quad\le\frac{1}{2}\int_{\Sigma}|\nabla\vec{\mathfrak G}_{\vec{\Phi}}|^2\ dx_1\wedge dx_2= \int_{\Sigma} (2+|d\vec{n}|^2_{g_{\vec{\Phi}}})\ dvol_{g_{\vec{\Phi}}}=\int_{\Sigma}2+2 \kappa_1\,\kappa_2+|\kappa_1-\kappa_2|^2\ dvol_{g_{\vec{\Phi}}}\\[5mm]
\ds\quad\quad=2\,\int_{\Sigma} K_{int} \ dvol_{g_{\vec{\Phi}}}+2\, \int_\Sigma |A^0_{\vec{\Phi}}|^2\, dvol_{g_{\vec{\Phi}}}=8\pi\ \mbox{deg}(\vec{\mathfrak G}_{\vec{\Phi}})+2\, \int_\Sigma |A^0_{\vec{\Phi}}|^2\, dvol_{g_{\vec{\Phi}}}
\end{array}
\ee
where $A^0_{\vec{\Phi}}$ is the trace free part of the second fundamental form of the immersion $\vec{\Phi}$ : ${\mathbb I}- H\, g_{\vec{\Phi}}$. $A^0_{\vec{\Phi}}$ has eigenvalues respectively  $(\kappa_1-\kappa_2)/2$ and $(\kappa_2-\kappa_1)/2$ and $|A^0_{\vec{\Phi}}|^2=2^{-1}\,|\kappa_1-\kappa_2|^2$. Observe that if $\vec{\Phi}$ is a minimal immersion into $S^3$, $\vec{n}_{\vec{\Phi}}$ is also a branched minimal immersion (see \cite{Law}), moreover, not only $\vec{\Phi}$ and $\vec{n}_{\vec{\Phi}}$ define the same conformal structure but in any local coordinates where $\vec{\Phi}$ is conformal, $\vec{n}_{\vec{\Phi}}$ is conformal as well.
Hence we have in particular
\[
|\p_{x_1}\vec{\mathfrak G}_{\vec{\Phi}}|^2=2\, |\p_{x_1}\vec{\Phi}\wedge\vec{n}_{\vec{\Phi}}+\vec{\Phi}\wedge\p_{x_1}\vec{n}_{\vec{\Phi}}|^2=2\, |\p_{x_1}\vec{\Phi}|^2+2\, |\p_{x_1}\vec{n}_{\vec{\Phi}}|^2=2\, |\p_{x_2}\vec{\Phi}|^2+2\, |\p_{x_2}\vec{n}_{\vec{\Phi}}|^2=|\p_{x_2}\vec{\mathfrak G}_{\vec{\Phi}}|^2
\]
and
\[
\p_{x_1}\vec{\mathfrak G}_{\vec{\Phi}}\cdot\p_{x_2}\vec{\mathfrak G}_{\vec{\Phi}}=2\, \lf[\p_{x_1}\vec{\Phi}\wedge\vec{n}_{\vec{\Phi}}+\vec{\Phi}\wedge\p_{x_1}\vec{n}_{\vec{\Phi}}\rg]\cdot \lf[\p_{x_2}\vec{\Phi}\wedge\vec{n}_{\vec{\Phi}}+\vec{\Phi}\wedge\p_{x_2}\vec{n}_{\vec{\Phi}}\rg]=0
\]
In other words, ${\vec{\mathfrak G}_{\vec{\Phi}}}$ is also conformal and the inequality in (\ref{can-1}) is an equality. Since $\vec{\Phi}$ is minimal we have
\be
\label{min-gauss-prep}
\begin{array}{l}
\ds0=\int_{\Sigma}|\kappa_1+\kappa_2|^2\ dvol_{g_{\vec{\Phi}}}=4\,\int_{\Sigma} K_{ext} \ dvol_{g_{\vec{\Phi}}}+2\, \int_\Sigma |A^0_{\vec{\Phi}}|^2\, dvol_{g_{\vec{\Phi}}}\\[5mm]
\ds\quad=4\,\int_{\Sigma} K_{int} \ dvol_{g_{\vec{\Phi}}}+2\, \int_\Sigma |A^0_{\vec{\Phi}}|^2\, dvol_{g_{\vec{\Phi}}}-4\, \int_\Sigma\ dvol_{g_{\vec{\Phi}}}\quad
\end{array}
\ee
Thus combining (\ref{can-1}) and (\ref{min-gauss-prep}), we have for $\vec{\Phi}$ minimal
\be
\label{minimal-gauss}
A(\vec{\Phi})=\int_{\Sigma}\ dvol_{\vec{\mathfrak G}_{\vec{\Phi}}}+8\pi\ \mbox{deg}(\vec{\mathfrak G}_{\vec{\Phi}})=4\, \int_\Sigma\ dvol_{g_{\vec{\Phi}}}
\ee

\medskip

Introduce the following M\"obius conformal transformations of $S^3$ for any $\vec{a}\in B^4$
\be
\label{mobius}
\Psi_{\vec{a}}({z}):=(1-|\vec{a}|^2)\ \frac{z-\vec{a}}{|z-\vec{a}|^2}-\vec{a}\quad.
\ee
It is well known that $\int_\Sigma |A^0_{\vec{\Phi}}|^2\, dvol_{g_{\vec{\Phi}}}$ is conformally invariant. Hence, for $\vec{\Phi}$ being a minimal immersion and for any $\vec{a}\in B^3$ we have
\be
\label{max-can}
\begin{array}{l}
\ds A(\Psi_{\vec{a}}\circ\vec{\Phi})\le\int_\Sigma |A^0_{\Psi_{\vec{a}}\circ\vec{\Phi}}|^2\, dvol_{g_{\Psi_{\vec{a}}\circ\vec{\Phi}}}+16\,\pi\ \mbox{deg}(\vec{\mathfrak G}_{\Psi_{\vec{a}}\circ\vec{\Phi}})=\int_\Sigma |A^0_{\vec{\Phi}}|^2\, dvol_{g_{\vec{\Phi}}}+16\,\pi\ \mbox{deg}(\vec{\mathfrak G}_{\vec{\Phi}})\\[5mm]
\ds\quad\quad=A(\vec{\Phi})=4\,\mbox{Area}(\vec{\Phi})
\end{array}
\ee
 We now work locally away from umbilic points.
\[
e^{2\,\la}:=|\p_{x_1}\vec{\Phi}|^2=|\p_{x_2}\vec{\Phi}|^2\quad\mbox{ and }\quad\vec{n}_{\vec{\Phi}}:=-\, e^{-2\la}\, \star\vec{\Phi}\wedge\p_{x_1}\vec{\Phi}\wedge\p_{x_2}\vec{\Phi}\quad.
\]
We have
\[
0=\vec{H}_{\vec{\Phi}}=H_{\vec{\Phi}}\ \vec{n}_{\vec{\Phi}}=2^{-1}\ (\Delta_{g_{\vec{\Phi}}}\vec{\Phi}+2\,\vec{\Phi})=\frac{e^{-2\la}}{2}\, (\Delta\vec{\Phi}+\vec{\Phi}\, |\nabla\vec{\Phi}|^2)
\]
Denote $\p_{z}:= 2^{-1}\, (\p_{x_1}-\,i\,\p_{x_2})$. We have
\[
\p_{z}^2\p_{\ov{z}}\vec{\Phi}=\p_z(\p^2_{z\ov{z}}\vec{\Phi})=-{4}^{-1}\,\p_z(\vec{\Phi}\,|\nabla\vec{\Phi}|^2)
\]
Hence, using $\vec{n}_{\vec{\Phi}}\cdot\p_{z}\vec{\Phi}=0$ and the conformality condition $\p_{z}\vec{\Phi}\cdot\p_z\vec{\Phi}=0$ we deduce
\be
\label{II.1}
\vec{n}\cdot\Delta\p_{z}\vec{\Phi}=0\quad.
%\mbox{ and }\quad \p_{z}\vec{\Phi}\cdot\p_{z}^2\p_{\ov{z}}\vec{\Phi}=0\quad.
\ee
This implies
\be
\label{II.2}
\p_{\ov{z}}\lf(\vec{n}\cdot\p^2_{z^2}\vec{\Phi}\rg)=\p_{\ov{z}}\vec{n}\cdot\p^2_{z^2}\vec{\Phi}=e^{2\la}\, \p_{\ov{z}}\vec{n}\cdot\p_z\lf(e^{-2\la}\,\p_z\vec{\Phi}\rg)+2\, \p_z\la\ \p_{\ov{z}}\vec{n}\cdot\p_z\vec{\Phi}\quad.
\ee
We recall the definition of the {\it Weingarten quadratic form}
\[
\vec{h}^{\,0}:=e^{2\la}\, \p_z\lf(e^{-2\la}\,\p_z\vec{\Phi}\rg)\ dz\otimes dz
\]
It is well known that $\p_z\lf(e^{-2\la}\,\p_z\vec{\Phi}\rg)$ is parallel to $\vec{n}$, indeed we have respectively
\[
\vec{\Phi}\cdot \p_z\lf(e^{-2\la}\,\p_z\vec{\Phi}\rg)=\p_z\lf(e^{-2\la}\,\vec{\Phi}\cdot\p_z\vec{\Phi}\rg)-e^{-2\la}\,\p_z\vec{\Phi}\cdot\p_z\vec{\Phi}=0\quad,
\]
moreover
\[
2\, e^{-2\la}\, \p_{z}\vec{\Phi}\cdot\p_z\lf(e^{-2\la}\,\p_z\vec{\Phi}\rg)=\p_z\lf(e^{-4\la}\, \p_z\vec{\Phi}\cdot\p_z\vec{\Phi}\rg)=0\quad,
\]
and finally
\[
2\,\p_{\ov{z}}\vec{\Phi}\cdot\p_z\lf(e^{-2\la}\,\p_z\vec{\Phi}\rg)=\p_z\lf(e^{-2\la}\, e^{2\la}\rg)-e^{-2\la}\, \p_{\ov{z}}\lf(\p_z\vec{\Phi}\cdot\p_z\vec{\Phi}\rg)=0\quad.
\]
We have also
\[
\p_{\ov{z}}\vec{n}\cdot\p_z\vec{\Phi}=-\,2^{-1}\ e^{2\la}\ H_{\vec{\Phi}}=0
\]
Hence (\ref{II.2}) implies
\be
\label{II.3}
\p_{\ov{z}}\lf(\vec{n}\cdot\p^2_{z^2}\vec{\Phi}\rg)=0\quad\Longrightarrow \vec{h}^{\,0}\ \mbox{ is an {\it holomorphic quadratic differential}.}
\ee
The complex valued function $f(z):=\vec{n}\cdot\p^2_{z^2}\vec{\Phi}$ is holomorphic. We shall now study the canonical family away from umbiic points of $\vec{\Phi}$ since they are isolated and contribute in an inessential way to
the homology at the boundary. Away from the umbilic points the function $f(z)$ is non zero and by replacing the coordinates $z$ by $w$ such that $w'(z)=\sqrt{f(z)}$ we can choose it to be equal to $1/2$ .
 Recall that in such a case the principal curvatures $k_1$ and $k_2$ satisfy
 \[
 4\ \vec{n}_{\vec{\Phi}}\cdot\p^2_{z^2}\vec{\Phi}= e^{2\la}\, (k_1-k_2)=2\, e^{2\la}\, k_1
 \]
 Hence we have
 \[
  k_1=e^{-2\la}\quad\Rightarrow\quad -\, e^{-2\la}\ \Delta\la=K_{int}=1-k_1^2=1-e^{-4\la}
 \]
 which gives that $\la$ is a solution of the following {\it sinh-Gordon equation}
 \[
 -\Delta\, \la=2\, \sinh(2\la)\quad.
 \]
 Observe that we have
 \be
 \label{II.4}
 \lf\{
 \begin{array}{l}
 \ds\p_{x_1}\vec{n}_{\vec{\Phi}}=\ -\, k_1\ \p_{x_1}\vec{\Phi}=-\ e^{-2\la}\ \p_{x_1}\vec{\Phi}\\[5mm]
 \ds\p_{x_2}\vec{n}_{\vec{\Phi}}=\ k_1\ \p_{x_2}\vec{\Phi}=e^{-2\la}\ \p_{x_2}\vec{\Phi}
 \end{array}
 \rg.
 \ee
 which implies that $\vec{n}_{\vec{\Phi}}$ is also a conformal immersion and we have
 \be
 \label{II.5}
 |\p_{x_1}\vec{n}_{\vec{\Phi}}|^2=|\p_{x_2}\vec{n}_{\vec{\Phi}}|^2=e^{-2\la}\quad\mbox{ and }\quad\p_{x_1}\vec{n}_{\vec{\Phi}}\cdot\p_{x_2}\vec{n}_{\vec{\Phi}}=0\quad.
 \ee
 We can rewrite (\ref{II.4}) in the form
 \[
 \p_{\ov{z}}\vec{n}_{\vec{\Phi}}=-\,e^{-2\la}\, \p_z\vec{\Phi}
 \]
We have in particular
 \be
 \label{II.6}
 dvol_{g_{\vec{n}}}=2^{-1}\,|\nabla \vec{n}_{\vec{\Phi}}|^2\ dx^2= e^{-2\la}\ dx^2
 \ee
We take now the M\"obius transformations $\Psi_{\vec{a}}$ given by  (\ref{mobius}). We have for any $\vec{Y}\in T_z S^3$
\[
d\Psi_{\vec{a}}({z})\cdot\vec{Y}=\, e^{\mu(z)}\ \lf(\vec{Y}+2\, \vec{a}\cdot\vec{Y}\ \frac{z-\vec{a}}{|{z}-\vec{a}|^2}\rg)\quad\mbox{ where }\quad e^{\mu(z)}:=\frac{1-|\vec{a}|^2}{|z-\vec{a}|^2}=\frac{1-|\vec{a}|^2}{1+|\vec{a}|^2-2\,\vec{a}\cdot z}
\]
and
\[
\lf|d\Psi_{\vec{a}}({z})\cdot\vec{Y}\rg|=e^{\mu(z)}\ |\vec{Y}|\quad.
\]
In order to match with the notations in \cite{MR} for instance we introduce $\vec{g}=-\, 2 \, \vec{a}/(1+|\vec{a}|^2)$ Observe that with this notation one has for instance
\be
\label{ag}
\frac{1-|\vec{a}|^2}{1+|\vec{a}|^2}=\sqrt{1-|\vec{g}|^2}\quad,\quad e^{\mu(z)}=\frac{\sqrt{1-|\vec{g}|^2}}{1+\vec{g}\cdot z}\quad\mbox{ and }\quad \vec{\Phi}_{\vec{g}}+\vec{a}=\frac{\sqrt{1-|\vec{g}|^2}}{1+\vec{g}\cdot\vec{\Phi}}\ (\vec{\Phi}-\vec{a})\quad,
\ee
Let $\vec{\Phi}_{\vec{g}}:=\Psi_{\vec{a}}\circ\vec{\Phi}$ and denote $\vec{n}_{\vec{\Phi}_{\vec{g}}}:=\vec{n}_{\vec{g}}$ the Gauss unit vector associated to $\vec{\Phi}_{\vec{g}}$. In order to simplify notations we shall also write $\vec{n}$ for $\vec{n}_0=\vec{n}_{\vec{\Phi}}$. We have
\be
\label{vecn}
\vec{n}_{\vec{g}}= \vec{n}_{\vec{\Phi}}+2\, \vec{a}\cdot\vec{n}_{\vec{\Phi}}\ \frac{\vec{\Phi}-\vec{a}}{1+|\vec{a}|^2-2\,\vec{\Phi}\cdot\vec{a}}=\vec{n}_{\vec{\Phi}}-\vec{g}\cdot\vec{n}_{\vec{\Phi}}\ \frac{\vec{\Phi}-\vec{a}}{1+\vec{g}\cdot\vec{\Phi}}=\vec{n}_{\vec{\Phi}}-\frac{\vec{g}\cdot\vec{n}_{\vec{\Phi}}}{\sqrt{1-|\vec{g}|^2}}\ (\vec{\Phi}_{\vec{g}}+\vec{a})
\ee
%We have in particular
%\be
%\label{II.8}
%\begin{array}{l}
%\ds(\vec{\Phi}_{\vec{g}}+\vec{a})\wedge\vec{n}_{\vec{g}}=\frac{1-|\vec{a}|^2}{|\vec{\Phi}-\vec{a}|^2}\lf(\vec{\Phi}\wedge\vec{n}_{\vec{\Phi}}-\vec{a}\wedge\vec{n}_{\vec{\Phi}}\rg)+\frac{\vec{g}\cdot\vec{n}_{\vec{\Phi}}}{1+\vec{g}\cdot\vec{\Phi}}\ \vec{a}\wedge\vec{\Phi}\\[5mm]
%\ds \quad=e^{\mu_{\vec{g}}(\vec{\Phi})}\ \lf(\vec{\Phi}\wedge\vec{n}_{\vec{\Phi}}-\vec{a}\wedge\vec{n}_{\vec{\Phi}}+\frac{\vec{g}\cdot\vec{n}_{\vec{\Phi}}}{\sqrt{1-|\vec{g}|^2}}\ \vec{a}\wedge\vec{\Phi}\rg)
%\end{array}
%\ee
We have
\[
\begin{array}{l}
\ds\p_{x_i}\vec{n}_{\vec{g}}=\p_{x_i}\vec{n}-\vec{g}\cdot\p_{x_i}\vec{n}\ \frac{\vec{\Phi}-\vec{a}}{1+\vec{g}\cdot\vec{\Phi}}-\vec{g}\cdot\vec{n}\ \frac{\p}{\p x_i}\lf(\frac{\vec{\Phi}-\vec{a}}{1+\vec{g}\cdot\vec{\Phi}}\rg)\\[5mm]

\ds \quad=e^{-\mu(\vec{\Phi})}\ d\Psi_{\vec{a}}\cdot\p_{x_i}\vec{n}-\vec{g}\cdot\vec{n}\ \frac{\p}{\p x_i}\lf(\frac{\vec{\Phi}-\vec{a}}{1+\vec{g}\cdot\vec{\Phi}}\rg)\\[5mm]
\ds\quad=e^{-\mu(\vec{\Phi})}\ d\Psi_{\vec{a}}\cdot\p_{x_i}\vec{n}-\frac{\vec{g}\cdot\vec{n}}{1+\vec{g}\cdot\vec{\Phi}} e^{-\mu(\vec{\Phi})}\ d\Psi_{\vec{a}}\cdot\p_{x_i}\vec{\Phi}\\[5mm]
\ds\quad=\lf((-1)^i\ e^{-2\la-\mu}-\frac{\vec{g}\cdot\vec{n}}{\sqrt{1-|\vec{g}|^2}}\rg)\ \p_{x_i}\vec{\Phi}_{\vec{g}}=\lf((-1)^i\ e^{-\la}-\frac{\vec{g}\cdot\vec{n}}{{1+\vec{\Phi}\cdot\vec{g}}}\ e^{\la}\rg)\ e^{-\mu-\la}\, \p_{x_i}\vec{\Phi}_{\vec{g}}
\end{array}
\]
We have in particular
\[
\vec{H}_{\vec{g}}:=\vec{H}_{\vec{\Phi}_{\vec{g}}}=e^{-\mu}\frac{\vec{g}\cdot\vec{n}}{{1+\vec{\Phi}\cdot\vec{g}}}\ \vec{n}_{\vec{a}}=\frac{\vec{g}\cdot\vec{n}}{\sqrt{1-|\vec{g}|^2}}\ \vec{n}_{\vec{g}}
\]
Denote 
\[
\vec{e}_{\vec{g}}^{\,i}:=e^{-\mu-\la}\, \p_{x_i}\vec{\Phi}_{\vec{g}}\quad\mbox{ and }\quad G_{\vec{g}}:=\vec{\Phi}_{\vec{g}}\wedge\vec{n}_{\vec{g}}\quad.
\]
We also introduce
\[
 \vec{G}^{\,\pm}_{\vec{g}}:=\frac{1}{\sqrt{2}}\, \lf(\vec{\Phi}_{\vec{g}}\wedge\vec{n}_{\vec{g}}\pm \star\,\vec{\Phi}_{\vec{g}}\wedge\vec{n}_{\vec{g}}\rg)=\frac{1}{\sqrt{2}}\, \lf(\vec{\Phi}_{\vec{g}}\wedge\vec{n}_{\vec{\vec{g}}}\pm \vec{e}_{\vec{g}}^{\,1}\wedge\vec{e}_{\vec{g}}^{\,2}\rg)
\]
We identify the unit self-dual (resp. anti-self-dual) 2-vectors in $\wedge^2{\R}^4$ to $S^2$. Each of the 2-spheres will be denoted respectively $S^2_\pm$. Let
\[
\vec{\mathfrak{G}}_{\vec{g}}:=(\vec{G}^{\,+}_{\vec{g}},\vec{G}^{\,-}_{\vec{g}})\in S^2_+\times S^2_-\quad.
\]
We have the following positive``Frenet Frame'' given by
\[
\lf\{
\begin{array}{l}
\ds \vec{E}_{1}^\pm=\frac{1}{\sqrt{2}} \lf(\vec{\Phi}_{\vec{g}}\wedge\vec{n}_{\vec{g}}\pm\vec{e}^{\,1}_{\vec{g}}\wedge\vec{e}^{\,2}_{\vec{g}}\rg)\\[5mm]
\ds \vec{E}_2^\pm=\frac{1}{\sqrt{2}} \lf(\vec{\Phi}_{\vec{g}}\wedge\vec{e}^{\,1}_{\vec{g}}\pm\vec{e}^{\,2}_{\vec{g}}\wedge\vec{n}_{\vec{g}}\rg)\\[5mm]
\ds \vec{E}_3^\pm=\frac{1}{\sqrt{2}} \lf(\vec{\Phi}_{\vec{g}}\wedge\vec{e}^{\,2}_{\vec{g}}\pm\vec{n}_{\vec{g}}\wedge\vec{e}^{\,1}_{\vec{g}}\rg)
\end{array}
\rg.
\]
Denote\footnote{By an abuse of notations we shall sometimes simply write $A$ (resp. $B$) for $A_{\vec{g}}(x)$ (resp. $B_{\vec{g}}(x)$). Observe that $A_{\vec{g}}(x)$ is nothing but the opposite of the mean curvature
 $H_{\vec{\Phi}_{\vec{g}}}$ of $\vec{\Phi}_{\vec{g}}$ at $x$.}
\[
A_{\vec{g}}(x):=-\frac{\vec{g}\cdot\vec{n}}{\sqrt{1-|\vec{g}|^2}}\quad\mbox{ and }\quad B_{\vec{g}}(x)=-\,e^{-2\la-\mu_{\vec{g}}(\vec{\Phi})}=-e^{-2\la}\ \frac{(1+\vec{g}\cdot\vec{\Phi})}{\sqrt{1-|\vec{g}|^2}}
\]
This notations are introduced in order to simplify the writing of $\p_{x_i}\vec{n}_{\vec{g}}$ in terms of $\p_{x_i}\vec{\Phi}_{\vec{g}}$ : we have respectively
\[
\p_{x_1}\vec{n}_{\vec{g}}= (A+B)\ \p_{x_1}\vec{\Phi}_{\vec{g}}\quad\mbox{ and }\quad \p_{x_2}\vec{n}_{\vec{g}}= (A-B)\ \p_{x_2}\vec{\Phi}_{\vec{g}}
\]

We compute first
\[
\lf\{
\begin{array}{l}
\ds\p_{x_1}\vec{G}^{\,+}_{\vec{g}}= e^{\la+\mu_{\vec{g}}}\ \lf[(A+B)\, \vec{E}_2^+-\vec{E}_3^+\rg]\quad\mbox{and }\quad \p_{x_2}\vec{G}^{\,+}_{\vec{g}}= e^{\la+\mu_{\vec{g}}}\ \lf[ \vec{E}_2^++(A-B)\, \vec{E}_3^+\rg]\\[5mm]
\ds\p_{x_1}\vec{G}^{\,-}_{\vec{g}}= e^{\la+\mu_{\vec{g}}}\ \lf[(A+B)\, \vec{E}_2^-+\vec{E}_3^-\rg]\quad\mbox{and }\quad \p_{x_2}\vec{G}^{\,-}_{\vec{g}}= e^{\la+\mu_{\vec{g}}}\ \lf[ -\vec{E}_2^-+(A-B)\, \vec{E}_3^-\rg]
\end{array}
\rg.
\]
This gives
\[
\lf\{
\begin{array}{l}
\ds(\vec{G}^{\,+}_{\vec{g}})^\ast\om_{S^2}=\vec{G}^{\,+}_{\vec{g}}\cdot \p_{x_1}\vec{G}^{\,+}_{\vec{g}}\times \p_{x_2}\vec{G}^{\,+}_{\vec{g}}\ dx_1\wedge dx_2=\ (1+A^2-B^2)\ dvol_{\vec{\Phi}_{\vec{g}}}\\[5mm]
\ds(\vec{G}^{\,-}_{\vec{g}})^\ast\om_{S^2}=\vec{G}^{\,-}_{\vec{g}}\cdot \p_{x_1}\vec{G}^{\,-}_{\vec{g}}\times \p_{x_2}\vec{G}^{\,-}_{\vec{g}}\ dx_1\wedge dx_2=\ (1+A^2-B^2)\ dvol_{\vec{\Phi}_{\vec{g}}}
\end{array}
\rg.
\]
We have also
\[
\lf\{
\begin{array}{l}
\ds G^\pm_{\vec{g},11}:=\p_{x_1}\vec{G}^{\,\pm}_{\vec{g}}\cdot\p_{x_1}\vec{G}^{\,\pm}_{\vec{g}}=\ e^{2\la+2\mu_{\vec{g}}}\ (1+(A+B)^2)\\[5mm]
\ds G^\pm_{\vec{g},22}:=\p_{x_2}\vec{G}^{\,\pm}_{\vec{g}}\cdot\p_{x_2}\vec{G}^{\,\pm}_{\vec{g}}=\ e^{2\la+2\mu_{\vec{g}}}\ (1+(A-B)^2)\\[5mm]
\ds G^\pm_{\vec{g},12}:=\p_{x_1}\vec{G}^{\,\pm}_{\vec{g}}\cdot\p_{x_2}\vec{G}^{\,\pm}_{\vec{g}}=\pm\, 2\, \ e^{2\la+2\mu_{\vec{g}}}\ B
\end{array}
\rg.
\]
This gives
\[
dvol_{\vec{G}^{\,+}_{\vec{g}}}= \sqrt{\lf|(1+(A-B)^2)\,(1+(A+B)^2)-4\, B^2\rg|}\ dvol_{\vec{\Phi}_{\vec{g}}}
\]
The metric induced by the {\it Gauss immersion} $\vec{\mathfrak{G}}_{\vec{g}}$ is given by
\[
g_{\vec{\mathfrak{G}}_{\vec{g}}}=2\, e^{2\la+2\mu_{\vec{g}}}\ \lf((1+(A+B)^2)\ dx_1^2+(1+(A-B)^2)\ dx^2_2\rg)
\]
and the corresponding volume form is then
\be
\label{forme-volume-gauss}
dvol_{\vec{\mathfrak{G}}_{\vec{g}}}=2\, \sqrt{(1+(A-B)^2)\,(1+(A+B)^2)}\ dvol_{\vec{\Phi}_{\vec{g}}}=2\, B^{-2}\, \sqrt{(1+(A-B)^2)\,(1+(A+B)^2)}\ e^{-2\la}\ dx^2
\ee
The associated {\it Jacobian to the Lagrangian mapping} $\vec{\mathfrak{G}}_{\vec{g}}$ is  by definition the function $C_{\vec{g}}$ given by
\[
(\vec{G}^{\,+}_{\vec{g}})^\ast\om_{S^2}=C_{\vec{g}}(x)\,dvol_{\vec{\mathfrak{G}}_{\vec{g}}}
\]
It ``measures'' in particular how far is the map $\vec{G}^{\,+}_{\vec{g}}$ for being conformal from $(\Sigma,g_{\vec{\mathfrak{G}}_{\vec{g}}})$ into $S^2$ : we have obviously $|C_{\vec{g}}(x)|\le 1/2$ with $|C_{\vec{g}}(x)|=1/2$ if and only if $\vec{G}^{\,+}_{\vec{g}}$ is conformal at $x$. The above computations give
\[
C_{\vec{g}}(x)=\frac{1+A^2-B^2}{2\,\sqrt{(1+(A-B)^2)\,(1+(A+B)^2)}}
\]
%Inserting the exact values of $A=-\vec{g}\cdot\vec{n}/\sqrt{1-|\vec{g}|^2}$ and $B=-e^{-2\la}\ (1+\vec{g}\cdot\vec{\Phi})/\sqrt{1-|\vec{g}|^2}$ and denoting $f:=e^{-2\la}\ (1+\vec{g}\cdot\vec{\Phi})$ give
%\be
%\label{II.9}
%C_{\vec{g}}(x)=\frac{1-|\vec{g}|^2+|\vec{g}\cdot\vec{n}|^2-f^2}{2\,\sqrt{(1-|\vec{g}|^2+(\vec{g}\cdot\vec{n}-f)^2)\ (1-|\vec{g}|^2+(\vec{g}\cdot\vec{n}+f)^2)}}
%\ee
We interrupt at this stage the proof of  lemma~\ref{lm-canonical family} in order to establish 4 intermediate results. We prove the following ``no neck energy'' lemma
\begin{Lm}
\label{lm-neck-region}{\bf [ No ``neck energy'' lemma.]} Under the previous notations we have that
\be
\label{noneck}
\lim_{\delta\rightarrow 0^+}\lim_{|\vec{g}|\rightarrow 1^-}\int_\Sigma {\mathbf 1}_{{\mathfrak B}_{\vec{g}}^\delta}\ dvol_{\vec{\mathfrak{G}}_{\vec{g}}}=0
\ee
where ${\mathbf 1}_{{\mathfrak B}_{\vec{g}}^\delta}$ is the characteristic function of the following set
\[
{\mathfrak B}_{\vec{g}}^\delta:=\lf\{x\in \Sigma\quad;\quad 4\ |C_{\vec{g}}(x)|^2<1-\delta\rg\}\quad.
\]
For $\vec{g}$ contained in a neighborhood of $\vec{\Phi}(\Sigma)$ such  that $\vec{g}$ admits a unique projection onto $\vec{\Phi}(\Sigma)$ we denote by $x_{\vec{g}}$ the pre-image by $\vec{\Phi}$ on $\Sigma$. 
and let $t_{\vec{g}}\in (-\pi/2,\pi/2)$ such that
\[
-\vec{g}=|\vec{g}|\ \lf(\cos\,t_{\vec{g}}\ \vec{\Phi}(x_{\vec{g}})+\sin\,t_{\vec{g}}\ \vec{n}(x_{\vec{g}})\rg)
\]
Denote $d_{\vec{g}}:=\sqrt{(1-|\vec{g}|)+t_{\vec{g}}^2}$. For any $\eta>0$ we define the $\eta-${\bf bubble} at $x_{\vec{g}}$ to be subset of $\Sigma$ (or indifferently it's image by $\vec{\Phi}$)
\[
{\mathcal B}_{\eta,{\vec{g}}}:=\lf\{x\in \Sigma\quad;\quad |x-x_{\vec{g}}|_{\vec{\Phi}}\le \eta^{-1}\ d_{\vec{g}}\rg\}\quad.
\]
where $d_{\vec{g}}:=\sqrt{(1-|\vec{g}|)+t_{\vec{g}}^2}$. Then, under these notations we have respectively
\be
\label{neck-lag}
\limsup_{\vec{g}\rightarrow \vec{\Phi}(\Sigma)}\|4\ |C_{\vec{g}}|^2-1\|_{L^\infty({\mathcal B}_{\eta,{\vec{g}}})}=0\quad,\quad \limsup_{\vec{g}\rightarrow \vec{\Phi}(\Sigma)}\|4\ |C_{\vec{g}}|^2-1\|_{L^\infty(\Sigma\setminus B_\eta(x_{\vec{g}}))}=0
\ee
and
\be
\label{noannular}
\lim_{\eta\rightarrow 0^+}\limsup_{\vec{g}\rightarrow \vec{\Phi}(\Sigma)}\int_{B_\eta(x_{\vec{g}}))\setminus{\mathcal B}_{\eta,{\vec{g}}}} dvol_{\vec{\mathfrak{G}}_{\vec{g}}}=0
\ee
where the conformally degenerating annular region ${\mathcal A}_{\eta,{\vec{g}}}:=B_\eta(x_{\vec{g}}))\setminus{\mathfrak B}_{\eta,{\vec{g}}}$ is the $\eta-${\bf neck}.
\hfill $\Box$
\end{Lm}
\noindent{\bf Proof of lemma~\ref{lm-neck-region}.} We have
\[
4\ |C_{\vec{g}}(x)|^2<1-\delta\quad \Longleftrightarrow\quad\delta\ \frac{1+(A^2-B^2)^2}{A^2+B^2}+2\ \frac{A^2-B^2}{A^2+B^2}<2\ (1-\delta)
\]
Hence if $4\ |C_{\vec{g}}(x)|^2<1-\delta$ we have
\be
\label{II.10}
x^2+tx-t<0\quad\Longleftrightarrow\quad -\sqrt{t^2+4t}\le 2\,x+t\le \sqrt{t^2+4t},
\ee
where $x:=(A^2-B^2)/(A^2+B^2)$ and $t:=2\delta^{-1}(A^2+B^2)^{-1}$. Observe
\[
\frac{2}{t}=\delta\ (A^2+B^2)=\delta\ \frac{|\vec{g}\cdot\vec{n}|^2+e^{-4\la}\ (1+\vec{g}\cdot\vec{\Phi})^2}{1-|\vec{g}|^2}\quad.
\]
First we consider for any $\eta>0$ the set 
\[
E_{\vec{g}}^\eta:=\lf\{x\in \Sigma\quad;\quad e^{-2\la(x)}\, (1+\vec{\Phi}(x)\cdot \vec{g})>\eta\rg\}
\]
On this set we have
\[
\frac{\eta^2}{1-|\vec{g}|^2}\le B^2\le \frac{C}{1-|\vec{g}|^2}
\]
and $t\simeq 1-|\vec{g}|$ and the condition (\ref{II.10}) implies
\[
\lf|\frac{A^2-B^2}{A^2+B^2}\rg|\le\, C\ \sqrt{1-|\vec{g}|}\quad \Longleftrightarrow\quad|A^2-B^2|=O\lf(\frac{1}{\sqrt{1-|\vec{g}|}}\rg)
\]
Hence using the explicit expression of $dvol_{\vec{\mathfrak{G}}_{\vec{g}}}$ given by (\ref{forme-volume-gauss}) we obtain
\be
\label{bad-far}
\int_{{\mathfrak B}_{{\vec{g}}}^\delta\cap E_{\vec{g}}^\eta}dvol_{\vec{\mathfrak{G}}_{\vec{g}}}\le C_{\delta,\eta}\ \sqrt{1-|\vec{g}|}
\ee
We consider now the case where $-\vec{g}$ is ``close'' to $\vec{\Phi}(\Sigma)$. Since $\vec{\Phi}$ is assumed to be an embedding there exits a unique
$x_{\vec{g}}$ such that $\vec{\Phi}(x_{\vec{g}})$ is the closest point in $\vec{\Phi}(\Sigma)$ for the $S^3$ distance to $-\vec{g}/|\vec{g}|$. Then there exists $t_{\vec{g}}\in (-\pi/2,\pi/2)$ such that
\[
-\vec{g}=|\vec{g}|\ \lf(\cos\,t_{\vec{g}}\ \vec{\Phi}(x_{\vec{g}})+\sin\,t_{\vec{g}}\ \vec{n}(x_{\vec{g}})\rg)
\]
Let $d_{\vec{g}}:=\sqrt{(1-|\vec{g}|)+t_{\vec{g}}^2}$. Observe that we have
\be
\label{II.11}
1+\vec{g}\cdot\vec{\Phi}(x)=1-|\vec{g}| \lf(\cos\,t_{\vec{g}}\ \vec{\Phi}(x_{\vec{g}})\cdot\vec{\Phi}(x)+\sin\,t_{\vec{g}}\ \vec{n}(x_{\vec{g}})\cdot\vec{\Phi}(x)\rg)
\ee
The Taylor expansion of $\vec{\Phi}$ gives
\[
\vec{\Phi}(x)=\vec{\Phi}(x_{\vec{g}})+\sum_{i=1,2} (x^i-x_{\vec{g}}^i)\ \p_{x_i}\vec{\Phi}(x_{\vec{g}})+2^{-1}\sum_{i,j=1,2} (x^i-x_{\vec{g}}^i)\ (x^j-x_{\vec{g}}^j)\ \p^2_{x_ix_j}\vec{\Phi}(x_{\vec{g}})+O_x(|x-x_{\vec{g}}|^3)
\]
Hence
\[
\vec{\Phi}(x_{\vec{g}})\cdot\vec{\Phi}(x)=1-2^{-1}\, \ |x-x_{\vec{g}}|^2_{\vec{\Phi}}+O_x(|x-x_{\vec{g}}|^3_{\vec{\Phi}})\quad\mbox{ and }\quad \vec{n}(x_{\vec{g}})\cdot\vec{\Phi}(x)=O_x(|x-x_{\vec{g}}|^2)
\]
This gives
\be
\label{II.12}
1+\vec{g}\cdot\vec{\Phi}(x)=1-|\vec{g}|\,\cos t_{\vec{g}}+ (2^{-1}+ O(d_{\vec{g}}))\ |x-x_{\vec{g}}|^2_{\vec{\Phi}}+O_x(|x-x_{\vec{g}}|^3_{\vec{\Phi}})\
\ee
which implies
\be
\label{II.13}
B=-\ e^{-2\la}\lf(\sqrt{\frac{1-|\vec{g}|}{2}}+\frac{t_{\vec{g}}^2\ (1+o_g(1))}{2\,\sqrt{2}\, \sqrt{1-|g|}  }+\frac{(1+o_g(1))}{2\,\sqrt{2}\, \sqrt{1-|g|}  }\ \lf[|x-x_{\vec{g}}|^2_{\vec{\Phi}}+O_x(|x-x_{\vec{g}}|^3_{\vec{\Phi}})\rg]\rg)
\ee
and
\be
\label{II.14}
A_{\vec{g}}(x)=\frac{(1+o_{\vec{g}}(1))}{\sqrt{1-|\vec{g}|^2}}\ (t_g+O_x(|x-x_{\vec{g}}|^2)
\ee
which gives
\be
\label{II.15}
\frac{A_{\vec{g}}(x)}{B_{\vec{g}}(x)}=-\frac{e^{2\la}}{\sqrt{2}}\ \frac{(1+o_{\vec{g}})\ (t_g+O_x(|x-x_{\vec{g}}|^2))}{{2-2|\vec{g}|}+{t_{\vec{g}}^2\ (1+o_g(1))}+{\,(1+o_g(1))}\ \lf[|x-x_{\vec{g}}|^2_{\vec{\Phi}}+O_x(|x-x_{\vec{g}}|^3_{\vec{\Phi}})\rg]}
\ee
We argue by contradiction and we assume that there exists a sequence $\vec{g}_k\rightarrow \vec{\Phi}(\Sigma)$ and
\be
\label{II.15.0a}
\lim_{\delta\rightarrow 0}\lim_{k\rightarrow +\infty}\int_{{\mathfrak B}_{\vec{g}_k}^\delta}dvol_{\vec{\mathfrak{G}}_{\vec{g}_k}}>0
\ee
We shall be considering two cases. Denote $t_k$ for $t_{\vec{g}_k}$ as well as $d_k:=d_{\vec{g}_k}$ or $A_k$ (resp. $B_k$) for $A_{\vec{g}_k}$  (resp. $B_{\vec{g}_k}$)  and w.l.o.g we can assume $\vec{g}_k\rightarrow \vec{g}_\infty$ and denote $x_\infty:=x_{\vec{g}_\infty}$.

\medskip

\noindent{\bf Case 1.} 
\be
\label{II.15a}
\lim_{k\rightarrow +\infty} \frac{t_k}{1-|\vec{g}_k|}<+\infty
\ee
Let $1>\eta>0$ independent of $k$ but depending on $\delta$ that we are going to fix later on.  Under the assumption (\ref{II.15a} we have respectively
\be
\label{II.15b}
\limsup_{k\rightarrow +\infty}\lf\|\frac{A_k(x)}{B_k(x)}\rg\|_{L^\infty(|x-x_k|\le \eta)} <D
\ee
where $D>0$ is independent of $k$ and $\eta$. Under the assumption (\ref{II.15a}) we have for $|x-x_k|<\eta$
\be
\label{II.15c}
\frac{1}{|B_k(x)|}\simeq \frac{{d_k}}{d^2_k+|x-x_k|^2}\quad\mbox{ and }\quad |A_k(x)|\leq \frac{d^2_k+|x-x_k|^2}{{d_k}}\quad.
\ee
We denote 
$${\mathcal A}_{\eta,k}:=\lf\{x\in \Sigma\quad;\quad \eta^{-2}\,{(1-|\vec{g}_k|)}\le|x-x_k|^2\le \eta^2\rg\}\quad,$$
which is the $\eta-$neck region while the ``bubble'' is given by
\[
{\mathcal B}_{\eta,k}:=\lf\{x\in \Sigma\quad ;\quad |x-x_k|^2\le \eta^{-2}\ (1-|\vec{g}_k|)\rg\}
\]
We treat ${\mathcal A}_{\eta,k}$ and ${\mathcal B}_{\eta,k}$ separately.
In ${\mathcal B}_{\eta,{k}}$ the following holds :
\be
\label{II.16}
4\ |C_{\vec{g}_k}|^2(x)=\frac{\ds\lf(\frac{1}{B^2_k}+\frac{A^2_k}{B^2_k}-1\rg)^2}{\ds \lf(\frac{1}{B^2_k}+\lf(1-\frac{A_k}{B_k}\rg)^2\rg)\ \lf(\frac{1}{B^2_k}+\lf(1+\frac{A_k}{B_k}\rg)^2\rg)}=1- O(d_k)
\ee
We now estimate the area of the  annular neck region
\be
\label{II.17}
\int_{{\mathcal A}_{\eta,k}}dvol_{\vec{\mathfrak{G}}_{\vec{g}}}\leq\int_{{\eta}^{-1} \sqrt{d_k}}^{{\eta}} \lf(\frac{d_k}{r^4}+1 \rg)\ r\,dr\leq 3\,\eta^2\quad.
\ee
Hence combining (\ref{bad-far}), (\ref{II.16}) and (\ref{II.17}) we contradict (\ref{II.15.0a}) in the {\bf case 1} that is when (\ref{II.15a})  is assumed.

\medskip

\noindent{\bf Case 2.} 
\be
\label{II.23}
\lim_{k\rightarrow +\infty} \frac{t_k}{1-|\vec{g}_k|}=+\infty\quad\mbox{ and }\quad\lim_{k\rightarrow +\infty} \frac{t^2_k}{1-|\vec{g}_k|}=0\quad.
\ee
In that case we have for $|x-x_k|<\eta$
\be
\label{II.24}
\frac{1}{|B_k(x)|}\simeq \frac{\sqrt{1-|\vec{g}_k|}}{1-|\vec{g}_k|+|x-x_k|^2}\quad\mbox{ and }\quad |A_k(x)|=\frac{t_k+O(|x-x_k|^2)}{\sqrt{1-|\vec{g}_k|}}\quad.
\ee
For the bubble we take again this time
\[
{\mathcal B}_{\eta,k}:=\lf\{x\in \Sigma\quad ;\quad |x-x_k|^2\le \eta^{-2} \, (1-|\vec{g}_k|)\rg\}\quad\mbox{ and }\quad {\mathcal A}_{\eta,k}:=B_\eta(x_k)\setminus {\mathcal B}_{\eta,k}
\]
and we have
\be
\label{II.25}
\lf|\frac{A_k(x)}{B_k(x)}\rg|\simeq\frac{t_k+O(|x-x_k|^2)}{1-|\vec{g}_k|+|x-x_k|^2}\ge\eta^2\,\frac{t_k}{1-|\vec{g}_k|} (1-o_k(1))\quad \mbox{ in }\quad{\mathcal B}_{\eta,k}
\ee
Hence in ${\mathcal B}_{\eta,{k}}$ the following holds :
\be
\label{II.26}
4\ |C_{\vec{g}_k}|^2(x)=\frac{\ds\lf(\frac{1}{B^2_k}+\frac{A^2_k}{B^2_k}-1\rg)^2}{\ds \lf(\frac{1}{B^2_k}+\lf(1-\frac{A_k}{B_k}\rg)^2\rg)\ \lf(\frac{1}{B^2_k}+\lf(1+\frac{A_k}{B_k}\rg)^2\rg)}=1- O\lf(\frac{1-|\vec{g}_k|}{\eta^2\,t_k}\rg)
\ee
Moreover we have
\be
\label{II.27}
\int_{{\mathcal A}_{\eta,k}}dvol_{\vec{\mathfrak{G}}_{\vec{g}}}\leq(1-|\vec{g}_k|)\int_{{\eta}^{-1} \sqrt{1-|\vec{g}_k|}}^{{\eta}} \frac{1}{r^4}\lf(\frac{t^2_k+r^4}{ 1-|\vec{g}_k|}  +1 \rg)\ r\,dr\leq 3\,\eta^2\quad.
\ee

\medskip

\noindent{\bf Case 3.} 
\be
\label{II.18}
\lim_{k\rightarrow +\infty} \frac{t^2_k}{1-|\vec{g}_k|}\in {\R_+^\ast}\cup\{+\infty\}
\ee
Let $\eta>0$ independent of $k$ to be fixed later. In that case we have for $|x-x_k|<\eta$
\be
\label{II.19}
\frac{1}{|B_k(x)|}\simeq \frac{\sqrt{1-|\vec{g}_k|}}{d_k^2+|x-x_k|^2}\quad\mbox{ and }\quad |A_k(x)|\leq \frac{d_k+|x-x_k|^2}{\sqrt{1-|\vec{g}_k|}}\quad.
\ee
For the bubble we take this time
\[
{\mathcal B}_{\eta,k}:=\lf\{x\in \Sigma\quad ;\quad |x-x_k|\le \eta^{-1} \max\{1, K_\infty\}\,{d_k}\rg\}\quad\mbox{ and }\quad {\mathcal A}_{\eta,k}:=B_\eta(x_k)\setminus {\mathcal B}_{\eta,k}
\]
where
\[
0\le K_\infty=\lim_{k\rightarrow +\infty}\frac{\sqrt{1-|\vec{g}_k|}}{d_k}<+\infty\quad.
\]
We have
\be
\label{II.20}
\lf|\frac{A_k(x)}{B_k(x)}\rg|\simeq \frac{d_k}{d_k^2+|x-x_k|^2}\ge \frac{\eta^2}{d_k\ \max\{1,K_\infty\}}\quad \mbox{ in }\quad{\mathcal B}_{\eta,k}
\ee
Hence in ${\mathcal B}_{\eta,{k}}$  the following holds :
\be
\label{II.21}
4\ |C_{\vec{g}_k}|^2(x)=\frac{\ds\lf(\frac{1}{B^2_k}+\frac{A^2_k}{B^2_k}-1\rg)^2}{\ds \lf(\frac{1}{B^2_k}+\lf(1-\frac{A_k}{B_k}\rg)^2\rg)\ \lf(\frac{1}{B^2_k}+\lf(1+\frac{A_k}{B_k}\rg)^2\rg)}=1- O(\eta^{-2}\,\max\{1,K_\infty\}\,d_k)
\ee
Moreover we have
\be
\label{II.22}
\int_{{\mathcal A}_{\eta,k}}dvol_{\vec{\mathfrak{G}}_{\vec{g}}}\leq(1-|\vec{g}_k|)\int_{{{\eta}^{-1} d_k}}^{{\eta}} \frac{1}{r^4}\lf(\frac{d^2_k+r^4}{ 1-|\vec{g}_k|}  +1 \rg)\ r\,dr\leq \,\lf(1+\frac{1-|\vec{g}_k|}{d_k^2}\rg)\,\eta^2\quad.
\ee
Hence combining (\ref{bad-far}), (\ref{II.21}) and (\ref{II.22}) we contradict (\ref{II.15.0a}) in the {\bf case 3} that is when (\ref{II.18})  is assumed.

\medskip

In all cases we have obtained a contradiction to the assumption (\ref{II.15.0a}). This implies that (\ref{noneck}) holds.

\medskip

Regarding the second part of the lemma, assuming that (\ref{neck-lag}) or (\ref{noannular}) do not hold then one could extract a subsequence $\vec{g}_{k}$ converging towards $\vec{\Phi}(\Sigma)$ and such that either case 1, case 2 and case 3 hold. In the 3 cases we have established contradictions above. This concludes the proof of the lemma. \hfill $\Box$

\begin{Lm}
\label{lm-polarisation}{\bf[The polarization map away from the embedded surface.]}
Let $\vec{g}\in S^3$ we consider the following map which for any $\vec{n}\in S^3\cap T_{\vec{g}}S^3$ is defined by
\[
\begin{array}{rcl}
\ds{\mathcal P}_{\p B^4}({\vec{g}}) :\  S^2_+\ & \longrightarrow &\ S^2_-\\[3mm]
\ds\frac{1}{\sqrt{2}}\lf(\vec{g}\wedge\vec{n}+\ast (\vec{g}\wedge\vec{n})\rg)&\longrightarrow &\ds \frac{1}{\sqrt{2}}\lf(\vec{g}\wedge\vec{n}-\ast (\vec{g}\wedge\vec{n})\rg)
\end{array}
\]
The map ${\mathcal P}_{\p B^4}({\vec{g}})$ is a positive isometry of ${\R}^3$ and $$\vec{g}\in S^3\rightarrow {\mathcal P}_{\p B^4}({\vec{g}})\in SO(3)$$ identifies with the canonical negative\footnote{The map which to ${\mathbf g}$ assigns
\[
 \ {\mathbf y}\in {\Im m}({\mathbb H})\rightarrow {\mathbf g}{\mathbf  y}{\mathbf g}^{-1}={\mathbf g}{\mathbf  y}{\mathbf g}^\ast\]
is the degree $+2$ map from $S^3$ onto $SO(3)\simeq{\R}P^3$ while the map ${\mathbf g}\rightarrow {\mathbf g}^\ast$ is a degree $-1$ map from $S^3$ into itself.} double covering of $SO(3)$ by $S^3$ given by
${\mathcal P}_{\p B^4}({\vec{g}})({\mathbf y})=\mathbf{g}^\ast\,{\mathbf y}\,\mathbf{g}$ where $\mathbf g$ is the quaternion corresponding to $\vec{g}$ and ${\R}^3$ is identified with  $\Im m({\mathbb H})$. In particular we have
\be
\label{deg-sphere}
\mbox{deg}\lf({\mathcal P}_{\p B^4}\rg)=2
\ee
.\hfill $\Box$
\end{Lm}
\noindent{\bf Proof of lemma~\ref{lm-polarisation}.}
Let $(\vec{\ep}_i)_{i=1\cdots 4}$ be the canonical basis of ${\R}^4$. We take for ${\R}^3\simeq (\wedge^2{\R}^4)_\pm$ the canonical basis 
\[
\lf\{
\begin{array}{l}
\ds\vec{E}_1^\pm:=\frac{1}{\sqrt{2}}\lf(\vec{\ep}_1\wedge\vec{\ep}_2\pm\vec{\ep}_3\wedge\vec{\ep}_4\rg)\\[3mm]
\ds\vec{E}_1^\pm:=\frac{1}{\sqrt{2}}\lf(\vec{\ep}_1\wedge\vec{\ep}_3\pm\vec{\ep}_4\wedge\vec{\ep}_2\rg)\\[3mm]
\ds\vec{E}_1^\pm:=\frac{1}{\sqrt{2}}\lf(\vec{\ep}_1\wedge\vec{\ep}_4\pm\vec{\ep}_2\wedge\vec{\ep}_3\rg)
\end{array}
\rg.
\]
and we identify $\vec{E}_i^+$ with $\vec{E}_i^-$. The map ${\mathcal P}({\vec{g}})$ is explicitly given by
\[
{\mathcal P}_{\p B^4}({\vec{g}})\ :\ 
\lf(\begin{array}{c}
g_1n_2-g_2n_1+g_3n_4-g_4n_3\\[3mm]
g_1n_3-g_3n_1+g_4n_2-g_2n_4\\[3mm]
g_1n_4-g_4n_1+g_2n_3-g_3n_2
\end{array}\rg)
\ \longrightarrow\ 
\lf(\begin{array}{c}
g_1n_2-g_2n_1-g_3n_4+g_4n_3\\[3mm]
g_1n_3-g_3n_1-g_4n_2+g_2n_4\\[3mm]
g_1n_4-g_4n_1-g_2n_3+g_3n_2
\end{array}\rg)
\]
Using quaternionic representatives the map ${\mathcal P}({\vec{g}})$ is equal to the following map  $$(<{\mathbf n},{\mathbf i}{\mathbf g}>, <{\mathbf n},{\mathbf j}{\mathbf g}>,<{\mathbf n},{\mathbf k}{\mathbf g}>)\longrightarrow (<{\mathbf n},{\mathbf g}{\mathbf i}>, <{\mathbf n},{\mathbf g}{\mathbf j}>,<{\mathbf n},{\mathbf g}{\mathbf k}>)$$
where $<\cdot,\cdot>$ denotes the scalar product in ${\mathbb H}$. The map ${\mathcal P}({\vec{g}})$ corresponds in passing from the coordinates of ${\mathbf n}$ in the basis $({\mathbf i}{\mathbf g},{\mathbf j}{\mathbf g},{\mathbf k}{\mathbf g})$ to the coordinates of the same vector in the basis $({\mathbf g}{\mathbf i},{\mathbf g}{\mathbf j},{\mathbf g}{\mathbf k})$. In other words this is the map
\[
(y_1\,{\mathbf i}+y_2\,{\mathbf j}+y_3\,{\mathbf k})\ \longrightarrow \ {\mathbf g}^\ast(y_1\,{\mathbf i}+y_2\,{\mathbf j}+y_3\,{\mathbf k})\,{\mathbf g}\quad.
\]
This concludes the proof of lemma~\ref{lm-polarisation}.\hfill $\Box$

\medskip

\begin{Lm}
\label{lm-bubble}
{\bf[The Bubbles Shapes]} For any $0<\eta<1$ and  $X=(X_1,X_2)\in B_{\eta^{-1}}(0)\subset {\C}$ we introduce
\[
\vec{{\mathfrak T}}_{(x_{\vec{g}},t_{\vec{g}},|\vec{g}|)}(X)=\lf({\mathcal T}^+_{(x_{\vec{g}},t_{\vec{g}},|\vec{g}|)}(X),{\mathcal T}^-_{(x_{\vec{g}},t_{\vec{g}},|\vec{g}|)}(X)\rg):=\vec{\mathfrak G}_{\vec{g}}\lf(x_{\vec{g}}+e^{\la(x_{\vec{g}})}\, \sqrt{t_{\vec{g}}^2+2-2|\vec{g}|}\,X\rg)\quad.
\]
Assume as $\vec{g}\rightarrow \vec{\Phi}(\Sigma)$ that there exists $\al\in [-\pi/2,\pi/2]$ such that
\[
\frac{\sqrt{2-2|\vec{g}|}}{\sqrt{t_{\vec{g}}^2+2-2|\vec{g}|}}\longrightarrow \cos\al\quad\mbox{ and }\quad\frac{t_{\vec{g}}}{\sqrt{t_{\vec{g}}^2+2-2|\vec{g}|}}\longrightarrow \sin\al
\]
then
\be
\label{bubbleconv+}
{\mathcal T}^\pm_{(x_{\vec{g}},t_{\vec{g}},|\vec{g}|)}(X)\longrightarrow {\mathcal T}^\pm_\al(X)
\ee
where, using the quaternion counterparts ${\mathbf n}$, $\mathbf \Phi$ and $\mathbf X$ of respectively $\vec{n}(x_{\vec{g}})$, $\vec{\Phi}(x_{\vec{g}})$ and $\vec{X}:=x_1\vec{e}_1(x_{\vec{g}})+x_2\vec{e}_2(x_{\vec{g}})$
\[
{\mathcal T}^+_{(\al,x)}(X):=\frac{1-|X|^2}{1+|X|^2}
\lf(\begin{array}{c}
<{\mathbf n},{\mathbf i}{\mathbf \Phi}>\\[3mm]
<{\mathbf n},{\mathbf j}{\mathbf \Phi}>\\[3mm]
<{\mathbf n},{\mathbf k}{\mathbf \Phi}>
\end{array}\rg)
+2\, \frac{\cos\al}{1+|X|^2}\lf(\begin{array}{c}
<{\mathbf n},{\mathbf i}{\mathbf X}>\\[3mm]
<{\mathbf n},{\mathbf j}{\mathbf X}>\\[3mm]
<{\mathbf n},{\mathbf k}{\mathbf X}>
\end{array}\rg)
+2\, \frac{\sin\al}{1+|X|^2}\lf(\begin{array}{c}
<{\mathbf \Phi},{\mathbf i}{\mathbf X}>\\[3mm]
<{\mathbf \Phi},{\mathbf j}{\mathbf X}>\\[3mm]
<{\mathbf \Phi},{\mathbf k}{\mathbf X}>
\end{array}\rg)
\]
and
\[
{\mathcal T}^-_{(\al,x)}(X):=\frac{1-|X|^2}{1+|X|^2}
\lf(\begin{array}{c}
<{\mathbf n},{\mathbf \Phi}{\mathbf i}>\\[3mm]
<{\mathbf n},{\mathbf \Phi}{\mathbf j}>\\[3mm]
<{\mathbf n},{\mathbf \Phi}{\mathbf k}>
\end{array}\rg)
+2\, \frac{\cos\al}{1+|X|^2}\lf(\begin{array}{c}
<{\mathbf n},{\mathbf X}{\mathbf i}>\\[3mm]
<{\mathbf n},{\mathbf X}{\mathbf j}>\\[3mm]
<{\mathbf n},{\mathbf X}{\mathbf k}>
\end{array}\rg)
+2\, \frac{\sin\al}{1+|X|^2}\lf(\begin{array}{c}
<{\mathbf \Phi},{\mathbf X}{\mathbf i}>\\[3mm]
<{\mathbf \Phi},{\mathbf X}{\mathbf j}>\\[3mm]
<{\mathbf \Phi},{\mathbf X}{\mathbf k}>
\end{array}\rg)
\]
\end{Lm}
\noindent{\bf Proof of lemma~\ref{lm-bubble}.} Using the notations above, (\ref{ag}) and (\ref{vecn}) imply
\be
\label{II.28}
\lf\{
\begin{array}{l}
\ds\vec{\Phi}_{\vec{g}}=-\vec{a}-\frac{1}{e^{2\la}\, B_{\vec{g}}(x)}\ (\vec{\Phi}-\vec{a})\\[5mm]
\ds\vec{n}_{\vec{g}}=\vec{n}_{\vec{\Phi}}-\frac{A_{\vec{g}}(x)}{e^{2\la}\, B_{\vec{g}}(x)}\ (\vec{\Phi}-\vec{a})\\[5mm]
\ds \vec{a}=-\frac{\vec{g}}{1+\sqrt{1-|\vec{g}|^2}}
\end{array}
\rg.
\ee
We then deduce
\be
\label{II.29}
\vec{\Phi}_{\vec{g}}\wedge \vec{n}_{\vec{g}}=-\,\vec{a}\wedge\vec{n}_{\vec{\Phi}}+\lf(\vec{n}_{\vec{\Phi}}+A_{\vec{g}}(x)\,\vec{a}\rg)\wedge\frac{\vec{\Phi}-\vec{a}}{e^{2\la} B_{\vec{g}}(x)}
\ee
We have respectively in the $\eta-$bubbles (i.e. in the domain $|x-x_{\vec{g}}|_{\vec{\Phi}}\le \eta^{-1}\ \max(d_{\vec{g}},\sqrt{1-|\vec{g}|})$) we have
\be
\label{II.30}
-\,\vec{a}\wedge\vec{n}_{\vec{\Phi}}=-\,\vec{\Phi}(x_{\vec{g}})\wedge\vec{n}_{\vec{\Phi}}(x_{\vec{g}})+o_{d_{\vec{g}}}(1)\quad.
\ee
moreover
\be
\label{II.31}
\begin{array}{l}
\ds\vec{\Phi}(x)-\vec{a}=\vec{\Phi}(x_{\vec{g}})\lf[\frac{1-|\vec{g}|\ \cos\,t_{\vec{g}}+\sqrt{1-|\vec{g}|^2}}{1+\sqrt{1-|\vec{g}|^2}}\rg]-\sin\,t_{\vec{g}}\ \vec{n}_{\vec{g}}(x_{\vec{g}})\\[5mm]
\ds\quad+\sum_{i=1}^2(x^i-x^i_{\vec{g}})\ \p_{x_i}\vec{\Phi}(x_{\vec{g}})+O_x(|x-x_{\vec{g}}|^3)
\end{array}
\ee
and
\[
\frac{1}{e^{2\la} B_{\vec{g}}(x)}=-\frac{\sqrt{1-|\vec{g}|^2}}{1-|\vec{g}|\ \cos\, t_{\vec{g}}+(2^{-1}+O(d_{\vec{g}}))\ |x-x_{\vec{g}}|^2_{\vec{\Phi}}+O_x(|x-x_{\vec{g}}|^3)}
\]
moreover
\[
A_{\vec{g}}(x)=\frac{(1+o_{\vec{g}}(1))}{\sqrt{1-|\vec{g}|^2}}\ (t_g+O_x(|x-x_{\vec{g}}|^2)
\]
We deduce
\be
\label{II.32}
\begin{array}{l}
\ds\frac{\vec{\Phi}(x)-\vec{a}}{e^{2\la}\ B_{\vec{g}}(x)}=-\ \frac{{2-2|\vec{g}|^2}+t_{\vec{g}}^2\,\sqrt{1-|\vec{g}|^2}}{2-2|\vec{g}|+ t_{\vec{g}}^2+\ |x-x_{\vec{g}}|^2_{\vec{\Phi}}}\ (1+o(1))\ \vec{\Phi}(x_{\vec{g}})\\[5mm]
\ds+\ \frac{2\,t_{\vec{g}}\ \sqrt{1-|\vec{g}|^2}}{  2-2\,|\vec{g}|+  t_{\vec{g}}^2+\,|x-x_{\vec{g}}|^2_{\vec{\Phi}}}\ (1+o(1))\ \vec{n}(x_{\vec{g}})\\[5mm]
\ds-\ \frac{ 2\, \sqrt{1-|\vec{g}|^2}}{  2-2\,|\vec{g}|+ \, t_{\vec{g}}^2+\,|x-x_{\vec{g}}|^2}\ (1+o(1))\ \sum_{i=1}^2(x^i-x^i_{\vec{g}})\ \p_{x_i}\vec{\Phi}(x_{\vec{g}})
\end{array}
\ee
This gives in one hand
\be
\label{II.33}
\begin{array}{l}
\ds\vec{n}_{\vec{\Phi}}\wedge \frac{\vec{\Phi}(x)-\vec{a}}{e^{2\la}\ B_{\vec{g}}(x)}=\frac{{2-2|\vec{g}|^2}+t_{\vec{g}}^2\,\sqrt{1-|\vec{g}|^2}}{2-2|\vec{g}|+ t_{\vec{g}}^2+\ |x-x_{\vec{g}}|^2_{\vec{\Phi}}}\ \vec{\Phi}(x_{\vec{g}})\wedge \vec{n}_{\vec{\Phi}}(x_{\vec{g}})\\[5mm]
\ds+\ \frac{2\, \sqrt{1-|\vec{g}|^2}}{  2-2\,|\vec{g}|+ \, t_{\vec{g}}^2+\,|x-x_{\vec{g}}|^2_{\vec{\Phi}}}\ \ \sum_{i=1}^2(x^i-x^i_{\vec{g}})\ \p_{x_i}\vec{\Phi}(x_{\vec{g}})\wedge \vec{n}_{\vec{\Phi}}(x_{\vec{g}}) +o(1)
\end{array}
\ee
and in the other hand
\be
\label{II.34}
\begin{array}{l}
\ds A_{\vec{g}}(x)\ \vec{a}\wedge \frac{\vec{\Phi}(x)-\vec{a}}{e^{2\la}\ B_{\vec{g}}(x)}=\frac{2\,(t^2_{\vec{g}}+t_{\vec{g}}\,O(|x-x_{\vec{g}}|^2_{\vec{\Phi}})}{  2-2\,|\vec{g}|+  t_{\vec{g}}^2+\,|x-x_{\vec{g}}|^2_{\vec{\Phi}}}\ \vec{\Phi}(x_{\vec{g}})\wedge\vec{n}(x_{\vec{g}})  \\[5mm]
\ds- \frac{ 2\, (t_{\vec{g}}+O(|x-x_{\vec{g}}|^2_{\vec{\Phi}})}{  2-2\,|\vec{g}|+ \, t_{\vec{g}}^2+\,|x-x_{\vec{g}}|^2_{\vec{\Phi}}}\ \ \sum_{i=1}^2(x^i-x^i_{\vec{g}})\ \vec{\Phi}(x_{\vec{g}})\wedge\p_{x_i}\vec{\Phi}(x_{\vec{g}})+o(1)
\end{array}
\ee
Combining (\ref{II.29}), (\ref{II.30}), (\ref{II.33}) and (\ref{II.34}) gives
\be
\label{II.35}
\begin{array}{l}
\ds\vec{\Phi}_{\vec{g}}\wedge \vec{n}_{\vec{g}}=\lf[\frac{{2-2|\vec{g}|}+t_{\vec{g}}^2- |x-x_{\vec{g}}|^2_{\vec{\Phi}}\,}{2-2|\vec{g}|+ t_{\vec{g}}^2+\ |x-x_{\vec{g}}|^2_{\vec{\Phi}}}\rg]\ \vec{\Phi}(x_{\vec{g}})\wedge\vec{n}_{\vec{\Phi}}(x_{\vec{g}})\\[5mm]
\ds+\ \frac{2\, \sqrt{2}\, \sqrt{1-|\vec{g}|}}{  2-2\,|\vec{g}|+ \, t_{\vec{g}}^2+\,|x-x_{\vec{g}}|^2_{\vec{\Phi}}}\ \ \sum_{i=1}^2(x^i-x^i_{\vec{g}})\ \p_{x_i}\vec{\Phi}(x_{\vec{g}})\wedge \vec{n}_{\vec{\Phi}}(x_{\vec{g}}) \\[5mm]
\ds+\frac{ 2\, t_{\vec{g}}}{  2-2\,|\vec{g}|+ \, t_{\vec{g}}^2+\,|x-x_{\vec{g}}|^2_{\vec{\Phi}}}\ \ \sum_{i=1}^2(x^i-x^i_{\vec{g}})\ \p_{x_i}\vec{\Phi}(x_{\vec{g}})\wedge\vec{\Phi}(x_{\vec{g}})+o(1)
\end{array}
\ee
This concludes the proof of lemma~\ref{lm-bubble}.\hfill $\Box$
\begin{Lm}
\label{lm-polar-bubble}{\bf[The polarization map on the bubbles.]}
There exists a map
\[
\begin{array}{rcl}
\ds{\mathcal P}_{S^1\times \Sigma}\quad:\ [-\pi/2,\pi/2]\times\Sigma\quad &\longrightarrow&\quad SO(3)\\[5mm]
(\al,x)\quad&\longrightarrow & {\mathcal P}_{S^1\times \Sigma}(\al,x)
\end{array}
\]
such that
\be
\label{polar}
\forall \,(\al,x)\in [-\pi/2,\pi/2]\times\Sigma\quad\forall\, X\in {\C}\quad {\mathcal P}_{S^1\times \Sigma}(\al,x)\lf ({\mathcal T}^+_{(\al,x)}(X)\rg)={\mathcal T}^-_{(\al,x)}(X)\quad.
\ee
Moreover, using quaternionic notations the map ${\mathcal P}_{S^1\times \Sigma}(\al,x)$ is given by
\be
\label{expli-polar-form}
{\mathcal P}_{S^1\times \Sigma}(\al,x)(y_1{\mathbf i}+y_2{\mathbf j}+y_3{\mathbf k})=(\cos\al \, {\mathbf n}(x)+\sin\al\ {\mathbf {\Phi}}(x))^\ast(y_1{\mathbf i}+y_2{\mathbf j}+y_3{\mathbf k})(\cos\al \, {\mathbf n}(x)+\sin\al\ {\mathbf {\Phi}}(x))
\ee
In particular ${\mathcal P}_{S^1\times \Sigma}(-\pi/2,x)={\mathcal P}_{S^1\times \Sigma}(\pi/2,x)$ and ${\mathcal P}_{S^1\times \Sigma}$ extends into a map from $S^1\times \Sigma$ into $SO(3)$ and it's degree is given by
\be
\label{deg-sigma}
\mbox{deg}\,({\mathcal P}_{S^1\times \Sigma})=2\, (1-g(\Sigma))
\ee
\hfill $\Box$
\end{Lm}
\noindent{\bf Proof of lemma~\ref{lm-polar-bubble}} We consider the map from ${\R}^3$ into itself such that 
\[
{\mathcal P}(\al,x)\lf(\begin{array}{c}
<\cos\al \,{\mathbf n}+\sin\al\, {\mathbf \Phi},{\mathbf i}(-\sin\al \,{\mathbf n}+\cos\al \,{\mathbf \Phi})>\\[3mm]
<\cos\al \,{\mathbf n}+\sin\al\, {\mathbf \Phi},{\mathbf j}(-\sin\al \,{\mathbf n}+\cos\al \,{\mathbf \Phi})>\\[3mm]
<\cos\al \,{\mathbf n}+\sin\al\, {\mathbf \Phi},{\mathbf k}(-\sin\al \,{\mathbf n}+\cos\al\, {\mathbf \Phi})>
\end{array}\rg)=
\lf(\begin{array}{c}
<\cos\al \,{\mathbf n}+\sin\al\, {\mathbf \Phi},(-\sin\al \,{\mathbf n}+\cos\al \,{\mathbf \Phi})\,{\mathbf i}>\\[3mm]
<\cos\al \,{\mathbf n}+\sin\al\, {\mathbf \Phi},(-\sin\al \,{\mathbf n}+\cos\al \,{\mathbf \Phi})\,{\mathbf j}>\\[3mm]
<\cos\al \,{\mathbf n}+\sin\al\, {\mathbf \Phi},(-\sin\al \,{\mathbf n}+\cos\al\, {\mathbf \Phi})\,{\mathbf k}>
\end{array}\rg)
\]
and for $l=1,2$
\[
{\mathcal P}(\al,x)\lf(\begin{array}{c}
<\cos\al \,{\mathbf n}+\sin\al {\mathbf \Phi},{\mathbf i}\,{\mathbf e}_l>\\[3mm]
<\cos\al \,{\mathbf n}+\sin\al {\mathbf \Phi},{\mathbf j}\,{\mathbf e}_l>\\[3mm]
<\cos\al \,{\mathbf n}+\sin\al {\mathbf \Phi},{\mathbf k}\,{\mathbf e}_l>
\end{array}\rg)
=\lf(\begin{array}{c}
<\cos\al \,{\mathbf n}+\sin\al {\mathbf \Phi},{\mathbf e}_l\,{\mathbf i}>\\[3mm]
<\cos\al \,{\mathbf n}+\sin\al {\mathbf \Phi},{\mathbf e}_l\,{\mathbf j}>\\[3mm]
<\cos\al \,{\mathbf n}+\sin\al {\mathbf \Phi},{\mathbf e}_l\,{\mathbf k}>
\end{array}\rg)
\]
This gives (\ref{polar}) and (\ref{expli-polar-form}). We consider the map into $S^3$ given by
\[
[-\pi/2,\pi/2]\times \Sigma\ \rightarrow S^3\quad (\al,x)\longrightarrow \vec{n}_\al(x):=(\cos\al\, \vec{n}(x)+\sin\al\, \vec{\Phi}(x))
\]
It projects down to a smooth map into ${\R}P^3$. We have
\[
\begin{array}{l}
\ds\mbox{deg}\,({\mathcal P})=\frac{1}{|{\R}P^3|}\int_{[-\pi/2,\pi/2]\times \Sigma} \vec{n}_{\al}^{\, \ast}\om_{S^3}=\frac{2}{3!|S^3|}\star_{{\R}^4}\int_{[-\pi/2,\pi/2]\times \Sigma}\vec{n}_{\al}\wedge d\vec{n}_{\al}\wedge d\vec{n}_{\al}\wedge d\vec{n}_{\al}\\[3mm]
\ds=\frac{1}{\pi^2}\star_{{\R}^4} \int_{[-\pi/2,\pi/2]\times \Sigma}\vec{n}_{\al}\wedge \p_{\al}\vec{n}_{\al}\wedge \p_{x_1}\vec{n}_{\al}\wedge \p_{x_2}\vec{n}_{\al}\ d\al\wedge dx_1\wedge dx_2\\[3mm]
\ds=\frac{1}{2\,\pi}\star_{{\R}^4} \int_{\Sigma}\vec{n}\wedge\vec{\Phi}\wedge\p_{x_1}\vec{n}\wedge\p_{x_2}\vec{n}\ dx_{1}\wedge dx_2+\frac{1}{2\,\pi}\star_{{\R}^4} \int_{\Sigma}\vec{n}\wedge\vec{\Phi}\wedge\p_{x_1}\vec{\Phi}\wedge\p_{x_2}\vec{\Phi}\ dx_{1}\wedge dx_2\\[5mm]
\ds=2(1-g(\Sigma)-|\Sigma|)+2|\Sigma|=2\, (1-g(\Sigma))
\end{array}
\] 
This concludes the proof of lemma~\ref{lm-polar-bubble}.\hfill $\Box$

\medskip

\noindent{\bf Proof of lemma~\ref{lm-canonical family} continued.} We recall the notations used at the beginning of the proof for points $\vec{g}$ in $B^4$ in the neighbourhood of $\Sigma$. For $\vec{g}$ contained in a neighborhood of $\vec{\Phi}(\Sigma)$ such  that $\vec{g}$ admits a unique projection onto $\vec{\Phi}(\Sigma)$ we denote by $x_{\vec{g}}$ the pre-image by $\vec{\Phi}$ on $\Sigma$
and let $t_{\vec{g}}\in (-\pi/2,\pi/2)$ such that
\[
-\vec{g}=|\vec{g}|\ \lf(\cos\,t_{\vec{g}}\ \vec{\Phi}(x_{\vec{g}})+\sin\,t_{\vec{g}}\ \vec{n}(x_{\vec{g}})\rg)
\]
Denote $d_{\vec{g}}:=\sqrt{1-|\vec{g}|+t_{\vec{g}}^2}$ and 
\[
\Om_\ep:=\lf\{\vec{g}\in B^4\quad;\quad d_{\vec{g}}<\ep\rg\}
\]
In the annulus type domain $\Om_\ep\setminus\Om_{\ep/2}$ we interpolate smoothly the identity map on $\p \Om_\ep$ and the projection map $\pi\ :\  \p\,\Om_{\ep/2}\cap \Om_\ep \rightarrow \Sigma$ such that $\pi(\vec{g}):=x_{\vec{g}}$. We denote $\Xi_\ep$ this extension of the identity in $B^4\setminus \Om_\ep$ by this interpolation in $\Om_\ep\setminus\Om_{\ep/2}$. Finally we identify the points
\[
\vec{g}\simeq\vec{g}^{\,'}\quad\mbox{ where }\quad \vec{g},\vec{g}^{\,'}\in \p B^4\cap\p \Om_{\ep/2}\quad,\quad\pi(\vec{g})=\pi(\vec{g}^{\,'}) \quad\mbox{ and }\quad t_{\vec{g}}=-t_{\vec{g}^{\,'}}=\pi/2\quad.
\]
We denote by $Z\simeq B^4$ the resulting cell obtained by taking $B^4\setminus\Om_{\ep/2}$ and assuming the above identification (see figure 3). 

Observe that $\p Z=\p B^4\sqcup S^1\times \Sigma$. From lemma~\ref{lm-polar-bubble} we have that for any point $x\in \Sigma$
\[
{\mathcal P}_{S^1\times \Sigma}(-\pi/2,x)={\mathcal P}_{S^1\times \Sigma}(\pi/2,x)
\]
\begin{figure}
\begin{center}
%\psfrag{A}{$[0,1]\times S^1$}
\includegraphics[width=12cm,height=9cm]{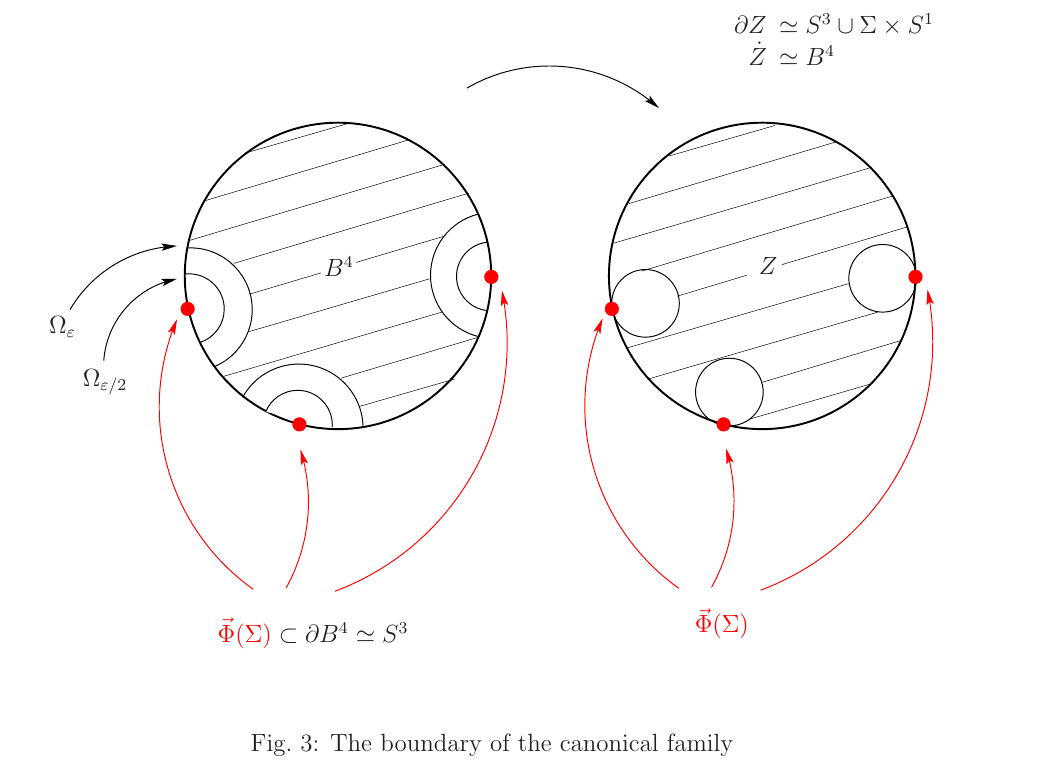}
\end{center}
\end{figure}

Hence there exists a continuous ${\mathbf F}$-limit of $\vec{\mathfrak G}_{\vec{\Phi}_{\Xi_\ep(\vec{g})}}$ on $\p Z$. We denote by $\vec{\mathfrak G}$ the resulting map in $C^0_{\mathbf F}(\overline{Z},{\mathcal M})$ which, thanks to (\ref{deg-sphere}) and (\ref{deg-sigma}), satisfies the required assumption :
\[
\vec{\mathfrak G}_\ast[\p Z]\ne 0 \quad\mbox{ in }\quad H_3(SO(3),{\Z})\quad.
\]
whenever  $g(\Sigma)\ne 0$. It is precisely $2\,g(\Sigma)$ time a generator of $H_3(SO(3),{\Z})$.  We have ensured that $\vec{\mathfrak G}(z)$ for any $z\in \p Z$ takes values into a union of spheres $S^2_{\mathcal P}$  but we still have to ensure that they are covering each sphere exactly once (i.e. in a conformal way) in order to ensure
\[
A(\vec{\mathfrak G}(z))\equiv 16\pi \quad\mbox{ on }\quad\p Z
\]
This can be obtained in slightly modifying $\vec{\mathfrak G}$ by extending $\vec{\mathfrak G}$ in a  collar neighbourhood of $\p Z$ (i.e. one add to $\p Z$ the cylinder $\p Z\times [0,1]$) taking for instance the Teichm\"uller harmonic map flow of Ding-Li-Liu/Rupflin-Topping (see \cite{RT} Theorem 1.1) of each $\vec{\mathfrak G}(z)$ into each given round 2-sphere $S^2_{\mathcal P}$ and keeping along the flow each of the  bubbles formed in finite time. This concludes the proof of lemma~\ref{lm-canonical family}.\hfill $\Box$

\subsection{A proof of the Willmore Conjecture in 3D.} In this subsection we are giving some elements to prove the following open problem
\begin{Con}
\label{th-will conj}  The value of the width of the second stage of the minmax hierarchy of the area in $S^3$ is given by
\be
\label{will-id}
8\pi^2={\mathcal W}_1\quad.
\ee
\hfill $\Box$
\end{Con}
Equality (\ref{will-id}) is implying the Willmore conjecture proved by Marques and Neves in \cite{MN}.
\begin{Co}
\label{co-willmore} Let $\vec{\Phi}$ be an immersion of an oriented closed surface $\Sigma$ in $S^3$ then the following inequality holds
\be
\label{willmore-conjecture}
2\pi^2\le W(\vec{\Phi})=\int_\Sigma (1+|\vec{H}_{\vec{\Phi}}|^2)\ dvol_{g_{\vec{\Phi}}}\quad.
\ee
Equality holds if and only if $\vec{\Phi}(\Sigma)$ is  conformally equivalent to the Clifford torus $\sqrt{2}^{-1} S^1\times S^1$.\hfill $\Box$
\end{Co}

Before going to the description of how a complete proof of open problem~\ref{th-will conj}  could be conducted we explain why the proof of the Willmore conjecture follows.

\medskip

 \noindent{\bf Proof of theorem~\ref{th-will conj} $\Rightarrow$ corollary~\ref{co-willmore}.}
 Indeed, let $\Sigma$ such that $g(\Sigma)\ne 0$. From (\ref{can-1}) we have
\[
A(\vec{\mathfrak G}_{\vec{\Phi}})=\int_\Sigma dvol_{g_{\vec{\mathfrak G}_{\vec{\Phi}}}}+8\,\pi\ \mbox{deg}\,(\vec{\mathfrak G}_{\vec{\Phi}})\le 4\ \int_\Sigma (1+|\vec{H}_{\vec{\Phi}}|^2)\ dvol_{g_{\vec{\Phi}}}=\,4\ W(\vec{\Phi})\quad.
\]
For a given surface $\Sigma$ the Willmore energy is conformally invariant. Hence we have
\be
\label{can-general}
\sup_{\vec{a}\in B^4}\ A(\vec{\mathfrak G}_{\vec{\Phi}_{\vec{a}}})\le\, 4\ W(\vec{\Phi})
\ee
Out of the $\vec{\Phi}_{\vec{a}}$ one construct as in the previous subsection an element in $\mbox{Sweep}_1({\mathfrak M}_{b(\Sigma)})$ such that
\be
\label{in-will}
{\mathcal W}_1\le{\mathcal W}_1(b(\Sigma))\le\max_{z\in Z}A(\vec{\mathfrak G}(z))\le\,4\ W(\vec{\Phi})\quad.
\ee
Combining (\ref{will-id}) and (\ref{in-will}) gives \ref{willmore-conjecture}.  

Assuming now equality holds. We consider the 4-dimensional canonical family issued from $\vec{\mathfrak G}_{\vec{\Phi}}$ . Then there must be a lagrangian minimal surface within the family  realizing the maximum otherwize we could ``push  down  infinitesimally'' the upper level $A^{-1}(8\pi^2-\delta,8\pi^2)$, which is made of smooth immersions, under  the level $8\pi^2-\delta$ for some $\delta>0$ using the {\it Lagrangian mean curvature flow} and that would lead to a contradiction. For the same reason we can moreover assume that all maxima are realized by a Lagrangian minimal surface. Starting now from this family using the pseudo-gradient flow of the viscous approximation of the area\footnote{Details are given in \cite{Riv-Lag-Visc}.} one constructs a critical point of the relaxed area
converging in flat norm to a maximum of the canonical family but realizing also at the limit, thanks to the main result in \cite{Riv-Lag-Visc} a minimal Lagrangian surface of index at most 4. As in the proof of theorem~\ref{th-will conj} below we deduce this maximum must be the Gauss map of a 
minimal surface isometric to the Clifford torus. Hence the original map is conformally equivalent to the Clifford torus. This concludes the proof of the corollary~\ref{co-willmore}.\hfill $\Box$

\medskip

\noindent{\bf Towards a proof of (\ref{will-id}).} We already know from corollary~\ref{co-lemma} that ${\mathcal W}_1\le 8\pi^2$. It remains to prove the other inequality. In order to prove this inequality the task is to show that ${\mathcal W}_1(b)$ is indeed achieved for any
$b>0$ such that $\mbox{Sweep}_1({\mathfrak M}_{b})\ne\emptyset$ by a minimal Lagrangian immersion of Hamiltonian index less or equal than the homological dimension  $4$.

In \cite{Riv-minmax} the author introduced a \underbar{PDE strategy} for producing minmax minimal surfaces  based on a relaxation procedure of the area that he called {\it viscosity method}.  After a series of works 
\cite{Riv-minmax,Riv-reg,Riv-Ind,PiRi1,PiRi2} partly in collaboration with Alessandro Pigati and also after using a work by Alexis Michelat \cite{Mich} the following result has been finally obtained
\begin{Th}
\label{th-visco}\cite{PiRi2}
Let $(N^n,g)$ be an arbitrary closed and smooth riemannian manifold, let $\Sigma$ be a smooth closed surface and ${\mathcal A}$ be an admissible homological family of ${\mathfrak M}(\Sigma)$ of dimension $d$ such that
\[
{\mathcal W}:=\inf_{\vec{\Phi}\in {\mathcal A}}\ \max_{z\in Z}\mbox{Area}(\vec{\Phi}(z))>0\quad.
\]
Then there exists a closed surface $S$ such that $g(S)\le g(\Sigma)$ and a minimal immersion $\vec{\Phi}$ of $S$ into $N$ such that
\[
{\mathcal W}=\mbox{Area}(\vec{\Phi})\quad\mbox{ and }\quad\mbox{Ind}(\vec{\Phi})\le d
\]
where $\mbox{Ind}(\vec{\Phi})$ denotes the Morse index of the minimal immersion $\vec{\Phi}$.\hfill $\Box$
\end{Th}
The goal would then be to extend the previous result to the framework of minimal Lagrangian immersions issued from Gauss maps into $S^2\times S^2$ and obtain that ${\mathcal W}_1$ is achieved by a possibly branched minimal Lagrangian Gauss map $\vec{\mathfrak G}_{\vec{\Phi}}$ of {\it Hamiltonian index} less or equal than $4$. Using lemma~\ref{reg-iso} and  proposition 6 of \cite{CaUr}\footnote{The main result in \cite{Urb} as well as proposition 6 of \cite{CaUr} are using the fact that $\vec{\Phi}$ is an immersion while the theorem~\ref{th-visco} and the viscosity method in general does not excludes the formation of branched points at the limit. Nevertheless, since the Gauss unit vectorfield is harmonic and smooth away from the branched points and most importantly in $W^{1,2}(\Sigma,S^3)$, the point removability for harmonic maps tells us that $\vec{n}$ extends smoothly throughout the branched points. Then the infinitesimal perturbation of the form $\vec{w}:= w\vec{n}$ of the surface $\vec{\Phi}$ for smooth $w$ are smooth and are ``admissible'' for the Jacobi operator. One can then follow the arguments in \cite{Urb} and \cite{CaUr} line by line (see more details in \cite{Riv-Lag-Visc}). }, since $A(\vec{\mathfrak G}_{\vec{\Phi}})>16\pi$ thanks to lemma~\ref{with-16-pi} we obtain that $\vec{\Phi}$ is either a multiple copy of a geodesic $S^2$ or is isometric to the Clifford Torus. The first option would imply $A(\vec{\mathfrak G}_{\vec{\Phi}})\ge 32\pi$ and this would contradict  (\ref{fund-ine}). Hence $\vec{\Phi}(\Sigma)$ is isometric to the Clifford Torus and we deduce (\ref{will-id}). This concludes the proof of theorem~\ref{th-will conj}.\hfill $\Box$

\subsection{Comments on the new proof of the Willmore Conjecture and possible extensions.}
\subsubsection{Some ingredients of the comparison with the proof of Marques and Neves}
While the above proof of the Willmore conjecture is clearly taking inspiration from the original one by Fernando Cod\'a Marques and Andr\'e Neves it differs by several aspects that we are going to stress now.

The canonical families used in \cite{MN} were 5 dimensional while the ones we are considering here are 4 dimensional. The 5th direction which was ``sweeping out'' the sphere $S^3$ in \cite{MN} is singular by nature :
It is given by the maps
\[
x\in \Sigma\rightarrow \cos t\,\vec{\Phi}(x)+\sin t\,\vec{n}(x)
\]
which are not defining immersions for all t but only {\it Lie-manifolds} (see \cite{Pal}). This fact has been for several years and up to now an obstruction for the author to implement the viscosity method with the 5-dimensional canonical family considered also by Ros in \cite{Ros-1}. This has been the case until he realized that by working with the Gauss maps instead of the map, this negative direction of the Jacobi field (the most negative one in fact among the 5) is ``killed'' and there  is a shift by one between the Morse Index of the underlying map and the Hamiltonian index of the Gauss map :
\be
\label{index-shift}
\mbox{Ind}(\vec{\Phi})-1=\mbox{Ham-Ind}(\vec{G}_{\vec{\Phi}})\quad.
\ee
The deformation argument in our approach in order to produce minmax solutions is mostly of \underbar{PDE nature} using a viscous regularisation of the mean-curvature equation (see \cite{Riv-minmax}) and a relaxation of  the classical Palais Smale theory. The approach of Marques and Neves is relying instead on the ``Almgren Pitts strategy'' which is using advanced {\it Geometric Measure Theory} and which is restricted up to now to the codimension 1 framework (see \cite{Riv-Bourb}). It is important to stress that, while it is restricted to surfaces, the theorem~\ref{th-visco}, and it's proof is completely independent of the co-dimension. which brings us to the next subsubsection

\subsubsection{The high-codimension Willmore Problem.}

While considering the Willmore problem in $S^n$ for arbitrary $n\ge 4$ the present work is naturally inviting us to look at the {\it Geodesic Gauss maps} of immersions into $S^n$ which to each normal direction assigns the geodesic it generates. It is defining a Lagrangian immersions 
$\vec{{\mathfrak G}}_{\vec{\Phi}}$ into the K\"ahler-Einstein Quintic $Gr(2,n+1)$ which is minimal if and only if the underlying immersion $\vec{\Phi}$ is minimal (see theorem 3.5 \cite{DrMc}). It coincides with the classical Gauss Map in codimension 1 (i.e. $n=3$). Moreover since $Gr(2,n+1)$ is K\"ahler Einstein, Lagrange Stationary is implying minimality (see \cite{SW}). Hence numerous ingredients are speaking in favour of investigating the higher codimensional Willmore problem in the framework of geodesic Gauss maps by the mean of the viscosity method. A first interesting question would consists in identifying the {\it polarisation map} ${\mathcal P}$ within the space of isometries of the Quintic and generated by the canonical family given by the composition with the M\"obius group of conformal transformations of $S^n$.

\subsubsection{The minmax free boundary minimal  surfaces and the critical catenoid.}

Another interesting direction the present work is suggesting would be to study by the mean of a minmax problem on the Gauss map the minimal area among free boundary minimal surfaces of the ball ${\mathbb B}^3$, which have more than one boundary components :
If $\vec{\Phi}$ is a minimal free boundary immersion it looks natural to introduce a minmax problem on the area of the Gauss map  $\vec{\mathfrak G}_{\vec{\Phi}}:=(\vec{\Phi},\vec{n}_{\vec{\Phi}})\in {\mathbb B}^3\times S^2$
which happens to be a {\it Legendrian minimal} map for the  contact form\footnote{Observe that $\al\wedge d\al\wedge d\al=dx_1\wedge dx_2\wedge dx_3\wedge\om_{S^2}$.}
\[
\al:=\sum_{i=1}^3 dx_i \, y_i\quad.
\] 
An element of the 3 dimensional canonical family would be given as follow :  for any free boundary embedding $\vec{\Phi}$ (not necessary minimal) with unit Gauss map $\vec{n}_{\vec{\Phi}}$
\[
\vec{a}\in {\mathbb B}^3\longrightarrow \quad (\vec{\Phi}_{\vec{a}},\vec{n}_{\vec{\Phi}_{\vec{a}}})\quad
\]
where $\vec{a}\in {\mathbb B}^3$
\[
\vec{\Phi}_{\vec{a}}=(1-|\va|^2) \frac{\iota(\vec{\Phi})-\vec{a}}{|\iota(\vec{\Phi})-\vec{a}|^2}-\vec{a}\quad\mbox{ and }\quad \iota(\vec{\Phi})=\frac{\vec{\Phi}}{|\vec{\Phi}|^2}\quad.
\]
The {\it polarization map} is now a map from $\p B^3\sqcup S^1\times\p\Sigma$ into ${\R}P^2$ which assigns the un-oriented geodesic in $S^2$ given by the limit\footnote{The product $\vec{\Phi}\times \vec{n}$ restricted to $\p \Sigma$ plays the role of $\vec{\Phi}\wedge\vec{n}$ on $\Sigma$ in the $S^3$ case. As $\vec{a}$ is going to $\p B^3\setminus\p\Sigma$ $\vec{\Phi}$ tends to a constant while $\vec{n}$ takes asymptotically values into a circle. Exactly as in the  closed case in $S^3$, this is the opposite in the bubble part as $\vec{a}$ tends to $\p \Sigma$ : $\vec{n}$ becomes constant while $\vec{\Phi}$ takes asymtoticaly values in a circle. While considering the product $\vec{\Phi}\times \vec{n}$ the two cases where $\vec{\Phi}$ and $\vec{n}$ are exchanging roles is now becoming a single one. This last fact is important in order to define the polarization map. } of $\vec{\Phi}_{\vec{a}}\times\vec{n}_{\vec{\Phi}_{\vec{a}}}$ as $\vec{a}$ converges respectively towards $\p B^3\setminus \p \Sigma$ or towards $\p\Sigma$ with the formation of bubbles similarly as in the $S^3$ case. The associated minmax should give a {\it minimal legendrian Gauss map} of {\it Legendrian Morse index} less or equal to 3 and is expected to be achieved by the Gauss map of the {\it critical catenoid}\footnote{Here again, as for the closed case, we expect the shift by one for each connected component between the Legendrian Morse Index of the minimal Gauss Map and the Morse index of the underlying minimal free boundary surface.} (see \cite{FrSc} for the introduction of the critical catenoid in the free boundary minimal surfaces context).
The viscosity method for free boundary minimal surfaces has been very recently developed by Alessandro Pigati in \cite{Pi1}.
\subsubsection{What comes next in the $S^3$ hierarchy ?}

Following the heuristic principles of the minmax hierarchy one first look at the family ${\mathcal C}_1$ of minimal surfaces realizing the previous step. This is the space of {\it oriented Clifford Tori} in $S^3$. Hence we have
\[
{\mathcal C}_1\simeq S^2\times S^2\quad.
\]
At this stage the hierarchy splits into two branches corresponding to the two generators of the $H_2(S^2\times S^2,{\Z})$. The ``first'' $S^2$ is given by the pull back of a fixed oriented closed geodesic of $S^2$ by the family of Hopf fibrations given by
\[
{\mathfrak h}_{\mathbf b}({\mathfrak q}):= {\mathfrak q} {\mathbf b} {\mathfrak q}^\ast        \quad \mbox{ where }\quad{\mathbf b}\in S^2\subset\Im {\mathbb H}
\]
the second $S^2$ corresponds to the choice of different closed geodesic in $S^2$. It would be interesting to study if these two generators would generate non empty admissible families of immersions $\mbox{Sweep}_2^1(S^3)$ and $\mbox{Sweep}_2^2(S^3)$ in the spirit of what happens for fibrations in section III. These two families would then generate Gauss Lagrangian minimal surfaces in the K\"ahler-Einstein manifold  $S^2\times S^2$ of Hamiltonian index at most 7=4+3 corresponding to minimal surfaces in $S^3$ of Morse index at most $8$.
\section{The  Definition of Homological Minmax Hierarchies.}
\reset
The iterative generation of minmax problems in the previous sections has been based on a general scheme that we are presenting now. While some of the problems we have considered were constructed out of the homotopy groups of the solutions to the previous minmax others out of the homology groups we are restricting our presentation below of the hierarchy principle to the ${\Z}_2$-cohomology of the  solutions to successive  minmax problems. A similar scheme could be developed for homotopy or homology based hierarchies.

\medskip

We shall denote ${\mathcal P}_N$ the category of $N-$dimensional compact manifolds orientable or non orientable, with or without boundary.

%of polyhedral ${\Z}_2-$chains\footnote{We recall that a polyhedral ${\Z}_2-$chain of dimension $N$ in an euclidian space ${\R}^Q$ can be identified with an integer rectifiable current given by a finite sum of the form
%\[
%\sum_{j=1}^L p_j\ {f_j}_\ast[\Delta^N]
%\]
%where $p_j\in {\Z}_2$, $f_j$ are affine maps from the oriented $N-$dimensional simplex $\Delta^N$ in ${\R}^{N-1}$ into ${\R}^Q$.}. And ${\mathcal P}_N:=\cup_{Q\in {\N}}{\mathcal P}_N({\R}^Q)$.

\subsection{The abstract scheme under the Palais Smale assumption.}

Let ${\mathfrak M}$ be an {\it Hilbert manifold} modeled on an Hilbert space ${\mathfrak H}$ or more generally a {\it Banach Manifold} modeled on a Banach space ${\mathfrak E}$  and assume that it is complete for the {\it Palais distance} $d_{\mathbf P}$ induced by the associated {\it Finsler Structure} $\|\cdot\|$ on $T\,{\mathfrak M}$. We shall restrict to the first case of
an Hilbert manifold {while considering exclusively} issues related to indices.

Let $E$ be a $C^2$ functional on ${\mathfrak M}$ which is assumed to satisfy the {\it Palais-Smale assumption} that is
\[
\forall\ \Phi_i\in {\mathfrak M}\quad\mbox{ s. t. }\quad\limsup_{i\rightarrow +\infty}E(\Phi_i)<+\infty\quad, \quad\lim_{i\rightarrow +\infty}\|DE(\Phi_i)\|_{\Phi_i}=0
\]
then
\[
\exists\, \Phi_{i'}\quad\mbox{ and }\quad \Phi_\infty\in {\mathfrak M}\quad\mbox{ s. t. }\quad d_{\mathbf P}(\Phi_{i'},\Phi_\infty)\longrightarrow 0\quad\mbox{ and }\quad DE(\Phi_\infty)=0
\]
A {\it Minmax Hierarchy} for the Functional $E$ requires the following objects

\begin{itemize}

%\item[a)] an at most countable sequence  of abstract closed surfaces $ S_0,S_1, S_2,\cdots$
\item[a)] a sequence (finite or infinite) of non zero integers $n_1,n_2,\cdots$
\item[b)] a sequence denoted 
\[
\mbox{Sweep}_{N_k}\subset\lf\{ (Y,\vec{\Phi})\ ;\ Y\in {\mathcal P}_{N_k}\quad\mbox{ and }\vec{\Phi}\in\mbox{Lip}_{\mathbf P}\lf(Y,{\mathfrak M}\rg) \rg\}\]
where $N_k:=n_1+\cdots+n_k$ and called {\it ${N_k}-$Sweepout space} of the hierarchy 
\item[c)] we have
\[
\forall\, (Y,\vec{\Phi})\in \mbox{Sweep}_{N_k}\quad\exists\, Z\in {\mathcal P}_{N_{k-1}}\quad\mbox{ s.t. }\ \p\,Y=\p( B^{n_k}\times Z)
\]
\item[d)] we have
\[
\forall\ x\in\p B^{n_k}\quad\quad \lf(Z,\vec{\Phi}(x,\cdot)\rg)\in \mbox{Sweep}_{N_{k-1}}
\]
\item[e)] a \underbar{strictly} increasing sequence of positive numbers (the {\it $k-$Widths} of the hierarchy)
\[
{W}_{k}:=\inf_{(Y,\vec{\Phi})\in \mbox{Sweep}_{N_k}}\quad\max_{y\in Y}\ E({\Phi}(y))
\]
\item[f)] for any homeomorphism $\Xi$ of $\mathfrak M$ satisfying
\[
\Xi({\Phi})={\Phi}\quad \mbox{ if }\quad E({\Phi})\le {W}_{k-1}-\eta_k
\]
for some $0<\eta_k<W_{k-1}-W_{k-2}$ we have
\[
\Xi\lf(\mbox{Sweep}_{N_k}\rg)\subset \mbox{Sweep}_{N_k}
\]
\end{itemize}

Using classical {\it Palais deformation theory} which applies to $E$  (see \cite{Riv-columb}) we obtain the existence of a sequence $\Phi_k$
such that
\[
E({\Phi}_k)={W}_k\quad.
\]
furthermore that the Banach Manifold ${\mathfrak M}$ is in fact {\it Hilbert}  and that $D^2E(\Phi)$ is Fredholm   we can ensure\footnote{We shall recall  the arguments leading to these assertions in the proof of 
theorem~\ref{th-I.1}} that $\mbox{Ind}({\Phi}_k)\le N_k$ where $\mbox{Ind}$ is the {\bf Morse Index} of $E$ (i.e. the maximal dimension of a vector space on which $D^2E$ is strictly negative).

\medskip

Of course the main issue in order to generate  such a hierarchy is to guarantee the series of strict inequalities between the successive $W_k$. We shall now explain a scheme that leads to such a series $W_k$.

 The first element in a hierarchy is an arbitrary {\it admissible familly} of the form
\[
\mbox{Sweep}_{n_1} \subset\lf\{
\begin{array}{l}
\ds(Y,\Phi)\in {\mathcal P}_{n_1}\times\mbox{Lip}(Y,{\mathfrak M})
\end{array}
\rg\}
\]
such that there exists $W_1>0$ satisfying
\[
\inf_{(Y,\Phi)\in\mbox{Sweep}_{n_1}}\max_{y\in Y}E(\Phi(y))= W_1
\]
and there exists $\eta>0$ such that for all homeomorphism $\Xi$ of ${\mathfrak M}$ equal to the identity for $E(\Phi)<W_1-\eta$ one has 
\[
\Xi(\mbox{Sweep}_{n_1})\subset \mbox{Sweep}_{n_1}\quad.
\]
Assuming now the hierarchy is constructed up to the order $k-1$, we introduce the notation for $l=1\cdots k-1$
\[
{\mathcal C}_l:=\lf\{ \Phi\in {\mathfrak M}\ ;\ \ ,\ E(\Phi)=W_{l}\quad\mbox{and }\quad DE(\Phi)=0\rg\}
\]
We are going to make the following assumption
\[
\mbox{(H1)}\quad\quad\quad {\mathcal C}_l \mbox{ \it is a smooth compact sub-manifold of }{\mathfrak M}
 \]
For any $\ep>0$ we denote
\[
{\mathcal O}_l(\ep):=\lf\{{\Phi}\in {\mathfrak M}\ ;\ d_{\mathbf P}({\Phi},{\mathcal C}_l)<\ep\rg\}
\]
Let $\ep_l>0$ be fixed such that $2\ep_l<\inf_{j<l}d_{\mathbf P}({\mathcal C}_l,{\mathcal C}_j)$ and
\[
\exists\ \pi_{l}\ \in\, \mbox{Lip}_{\mathbf P}\lf( {\mathcal O}_l(\ep_l) ,{\mathcal C}_l\rg)\quad\mbox{ s. t. }\quad \forall \ {\Phi}\in {\mathcal C}_l\quad\quad\pi_l({\Phi})={\Phi}
\]
as given by \cite{Lan}. The tubular neighborhood ${\mathcal O}_l(\ep_l) $ of ${\mathcal C}_l$  will simply be denoted ${\mathcal O}_l$. Because respectively of proposition~\ref{pr-A.1}  there exists $\delta_l>0$ such that
\be
\label{I.4}
\begin{array}{l}
\ds\forall\ (Y,\vec{\Phi})\in \mbox{Sweep}_{N_l}\quad\mbox{ , }\quad\max_{y\in Y}E({\Phi})\le {W}_l+\delta_l\\[3mm]
\ds\quad   \Longrightarrow\quad d_{\mathbf P}({\Phi}(Y),{\mathcal C}_l)<\ep_l\quad.
\end{array}
\ee
This being established we define $n_k$  as follows. Let $n_k\in{\N}^\ast$ such that\footnote{The reason why we are working
with ${\Z}_2$ cohomology comes from the fact that we are going to use Thom's resolution of Steenrod problem regarding the realization of ${\Z}_2-$homology classes
by continuous images of smooth manifolds (see \cite{Tho}).}
\[
H^{n_k-1}({\mathcal C}_{k-1},{\Z}_2)\ne 0
\]
and choose $\om_{k-1}$ being a non zero element of $H^{n_k-1}({\mathcal C}_{k-1},{\Z}_2)$. 

\medskip

Under the previous notations we define $\mbox{Sweep}_{N_k}$ to be the set of pairs $(Y,{\Phi})$ such that

\begin{itemize}
\item[i)] 
\[
Y\in {\mathcal P}_{N_{k}}\quad,\quad {\Phi}\in\mbox{Lip}_{\mathbf P}(\ov{Y},{\mathfrak M})\quad.
\]
\item[ii)] There exists $Z\in {\mathcal P}_{N_{k-1}}$ s.t.
\[
\p Y=\p\lf(B^{n_k}\times Z\rg)
\]
%where $B^{n_k}$ is the unit ball of center 0 in ${\R}^{n_k}$.
\item[iii)] We have
\[
\forall\ x\in\p B^{n_k}\quad\quad \lf(Z,\vec{\Phi}(x,\cdot)\rg)\in \mbox{Sweep}_{N_{k-1}}
\]
and
\[
\max_{y\in \p Y}E({\Phi}(y))\le {W}_{k-1}+2^{-1}\,\delta_{k-1}
\]
\item[iv)] Let 
\[
\Om_{{\Phi}}:=\lf\{  y\in \p Y\ ;\ d_{\mathbf P}({\Phi}(y),{\mathcal C}_{k-1})<\ep_{k-1} \rg\}
\]
we have\footnote{Observe that $[\Om_{{\Phi}}\cap (\{x\}\times Z)]\in H_{N_{k-1}}(\Om_{{\Phi}},{\Z}_2)$ is independent of $x\in \p B^{n_k}$ for $n_k>1$.}

\[
[\Om_{{\Phi}}\cap (\{x\}\times Z)]\in H_{N_{k-1}}(\,\Om_{{\Phi}},\p\,\Om_{{\Phi}} ,{\Z}_2)\quad\mbox{ is {\bf Poincar\'e dual} to } (\pi_{k-1}\circ{\Phi})^\ast\om_{k-1}\in H^{n_k-1}(\Om_{{\Phi}},{\Z}_2)
\]
\item[v)]  We have
\[
\max_{y\in  B^{n_k}\times \p Z} E({\Phi}(y))\le\ {W}_{k-1}-\delta_{k-1}\quad.
\]
\end{itemize}

\medskip

%In the case $\sigma=0$, a {\bf Minmax Hierarchy} is obtained in a same manner except that we replace $\mbox{Imm}(M^m)$ by $\mbox{Imm}^{g_0}(M^m)$ for some $g_0$ independent of $k$.
\begin{figure}
\begin{center}
%\psfrag{A}{$[0,1]\times S^1$}
\includegraphics[width=13cm,height=10cm]{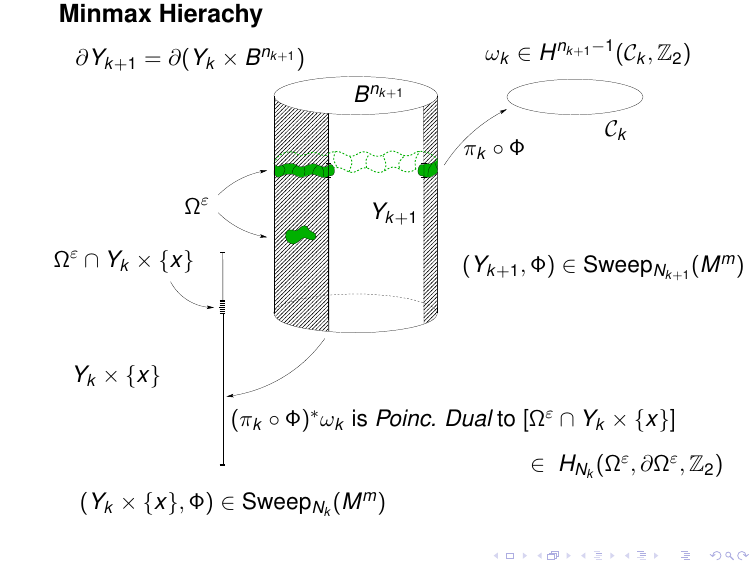}
\end{center}
\end{figure}

The main result of the present section is contained in the following theorem.
\begin{Th}
\label{th-I.1}
Under the hypothesis (H1) and assumptions i)...v) we have
\be
\label{strict-sigma}
{W}_{k-1}< {W}_k
\ee
and $\{\mbox{Sweep}_{N_l}\}_{l\le k}$ defines a {\bf Minmax Hierarchy}. Consequently, for any $l\le k$ there exists ${\Phi}_l$ such that
\[
E({\Phi}_l)={W}_l\quad,\quad DE({\Phi}_l)=0\quad
\]
Assuming furthermore that the Banach Manifold ${\mathfrak M}$ is in fact {\it Hilbert}  and that $D^2E(\Phi)$ is Fredholm for any $\Phi$, we have
\be
\label{indice-palais}
\mbox{Ind}\,({\Phi}_l)\le N_l
\ee
where $\mbox{Ind}$ is the {\bf Morse Index} of $E$ (i.e. the maximal dimension of a vector space on which $D^2E$ is strictly negative).

\hfill $\Box$
\end{Th}
\noindent{\bf Proof of theorem~\ref{th-I.1}.} We shall first prove by induction  that ${W}_{k-1}+\delta_{k-1}\le {W}_k$ where we recall that the sequence $\delta_k$ is defined in (\ref{I.4}). Assume this is not the case but assume $W_{l-1}+\delta_{l-1}\le W_l$ for $l\le k-1$. Choose $(Y,{\Phi})\in \mbox{Sweep}_{N_k}$
such that
\be
\label{I.6}
\max_{y\in Y}\, E({\Phi}(y))\le {W}_{k-1}+\delta_{k-1}
\ee
Consider 
\[
\La_{{\Phi}}=\lf\{  y\in Y\ ;\ d_{\mathbf P}({\Phi}(y),{\mathcal C}_{k-1})<\ep_{k-1} \rg\}\quad.
\]
and denote $\Xi_{{\Phi}}:=\p \La_{{\Phi}}\cap\mbox{int}(Y)$. Since $y\rightarrow d_{\mathbf P}({\Phi}(y),{\mathcal C}_{k-1})$ is lipschitz we can assume, without loss of generality, that $\Xi_{{\Phi}}$ realizes an {\it homology manifold} which makes $(\Xi_{{\Phi}},\p \Xi_{{\Phi}})$ a {\it Poincar\'e duality pair}  .

 For $\ep_{k-1}$ chosen small enough we have
\[
E(\Phi)\le W_{k-1}-\delta_{k-1}\quad\Longrightarrow\quad d_{\mathbf P}(\Phi,{\mathcal C}_{k-1})>\ep_{k-1}
\]
thus v) implies
\[
\Xi_{{\Phi}}\cap (\p B^{n_k}\times Z)=\emptyset
\]
hence
\[
\p\, \Xi_{{\Phi}}=\p\,\Om_{{\Phi}}
\]
We shall now prove the following.

\medskip

\noindent{\bf Claim 1}
\[
[\p\, \Om_{{\Phi}}\cap (\{x\}\times Z)]\in \mbox{Im}\,\p_\ast \quad,
\]
where $\p_\ast$ is the boundary operator 
\[
\p_\ast\ :\ H_{N_{k-1}}(\Xi_{{\Phi}},\p\,\Xi_{{\Phi}},{\Z}_2)\quad\longrightarrow\quad H_{N_{k-1}-1}(\p\, \Xi_{{\Phi}},{\Z}_2)
\]
We recall the following {\it relative Poincar\'e duality} commutative diagram
\[
\begin{CD}
H^{p-1}(\p\, \Xi_{{\Phi}},{\Z}_2) @>\delta^\ast>> H^p(\Xi_{{\Phi}},\p\, \Xi_{{\Phi}},{\Z}_2)@>j^\ast>> H^p(\Xi_{{\Phi}},{\Z}_2)@>i^\ast>>H^p(\p\,\Xi_{{\Phi}},{\Z}_2)\\
@VVDV    @VVDV  @VVDV  @VVDV\\
H_{N_{k}-1-p}(\p\, \Xi_{{\Phi}},{\Z}_2)        @>i_\ast>> H_{N_k-p-1}(\Xi_{{\Phi}},{\Z}_2)  @>j_\ast >> H_{N_k-p-1}(\Xi_{{\Phi}},\p \,\Xi_{{\Phi}},{\Z}_2)@>\p_\ast >> H_{N_k-p-2}(\p\, \Xi_{{\Phi}},{\Z}_2)
\end{CD}
\]
where the vertical arrows are Poincar\'e isomorphisms simply denoted by $D$ and given, modulo a sign, by cap products respectively with $[\Xi_{{\Phi}}]$ or $[\p \,\Xi_{{\Phi}}]$. We apply this diagram
to the case $p=n_k-1$. The map $\pi_{k-1}\circ{\Phi}$ is well defined on $\Xi_{{\Phi}}$ and
\[
(\pi_{k-1}\circ{\Phi})^\ast\om_{k-1}\in H^{n_k-1}(\Xi_{{\Phi}},{\Z}_2)
\]
It is clear that the image of $(\pi_{k-1}\circ{\Phi})^\ast\om_{k-1}$ by the restriction map $i^\ast$ is $(\pi_{k-1}\circ{\Phi})^\ast\om_{k-1}$
itself where ${\Phi}$ is restricted to $Z\times \p B^{n_k}$. Let $U\in H_{N_{k-1}}(\Xi_{{\Phi}},\p\,\Xi_{{\Phi}},{\Z}_2)$ be the Poincar\'e dual to
$(\pi_{k-1}\circ{\Phi})^\ast\om_{k-1}\in H^{n_k-1}(\Xi_{{\Phi}},{\Z}_2)$. Because of the last part of the above diagram we have that
\[
\p_\ast U\quad\mbox{ is the Poincar\'e dual to }(\pi_{k-1}\circ{\Phi})^\ast\om_{k-1}\in H^{n_k-1}(\p\,\Xi_{{\Phi}},{\Z}_2)
\]
We shall now make use of the following lemma
\begin{Lm}
\label{rm-I.3} The assumption
\[
[\Om_{{\Phi}}\cap (\{x\}\times Z)]\in H_{N_{k-1}}(\Om_{{\Phi}},\p\Om_{{\Phi}},{\Z}_2)\quad\mbox{ is {\bf Poincar\'e dual} to } (\pi_{k-1}\circ{\Phi})^\ast\om_{k-1}\in H^{n_k-1}(\Om_{{\Phi}},{\Z}_2)
\]
implies 
\[
[\p\,\Om_{{\Phi}}\cap (\{x\}\times Z)]\in H_{N_{k-1}-1}(\p\,\Om_{{\Phi}},{\Z}_2)\quad\mbox{ is {\bf Poincar\'e dual} to } (\pi_{k-1}\circ{\Phi})^\ast\om_{k-1}\in H^{n_k-1}(\p\Om_{{\Phi}},{\Z}_2)
\]
\hfill$\Box$
\end{Lm}
\noindent{\bf Proof of lemma~\ref{rm-I.3} .} The lemma is a direct consequence of the following {\it relative Poincar\'e duality} commutative diagram 
\[
\begin{CD}
 H^{n_k-1}(\Om_{{\Phi}},{\Z}_2)@>i^\ast>>H^{n_k-1}(\p\Om_{{\Phi}},{\Z}_2)\\
  @VVDV  @VVDV\\
 H_{N_{k-1}}(\Om_{{\Phi}},\p \Om_{{\Phi}},{\Z}_2)@>\p_\ast >> H_{N_{k-1}-1}(\p \Om_{{\Phi}},{\Z}_2)
\end{CD}
\]
\hfill $\Box$

\noindent{\bf End of the proof of claim 1}. Hence, since by the previous lemma $[\p\, \Om_{{\Phi}}\cap (\{x\}\times Z)]$ is also {\it Poincar\'e dual} to $(\pi_{k-1}\circ{\Phi})^\ast\om_{k-1}\in H^{n_k-1}(\p\,\Xi_{{\Phi}},{\Z}_2)$, by uniqueness of the {\it Poincar\'e dual} we have
\[
[\p\,\Om_{{\Phi}}\cap (\{x\}\times Z)]=\p_\ast U
\]
 and  $[\p\,\Om_{{\Phi}}\cap (\{x\}\times Z)]$ is a boundary in $\Xi_{{\Phi}}$ and the claim is proved.  Using Thom's proof of {\it Steenrod Problem} on the realization of ${\Z}_2$ homology classes by continuous images of smooth un-oriented manifolds (see theorem III.2 in \cite{Tho}), we can assume that the concrete chain $U$ in $Y$ is the image  of an element in ${\mathcal P}_{N_{k-1}}$.  By an abuse of notation we identify $U$ with this element in ${\mathcal P}_{N_{k-1}}$. By ''pushing''  $U$  inside $Y\setminus \La_{{\Phi}}$, and summing this  with the homology manifolds $Z\times \{x\}\setminus \Om_{{\Phi}}\cap (\{x\}\times Z)$ we obtain and $(V,{\Phi})\in \mbox{Sweep}_{N_{k-1}}$. Because of (\ref{I.6}) we have 
\[
\max_{y\in V} E({\Phi}(y))\le {W}_{k-1}+\delta_{k-1}
\]
Using proposition~\ref{pr-A.1} we then have the existence of $y_\ep\in V$ such that $d_{\mathbf P}({\Phi}(y),{\mathcal C}_{k-1})<\ep_{k-1}$ which is a contradiction. Hence we have
\[
{W}_{k-1}+\delta_{k-1}\le {W}_k
\]
Consider a Pseudo-gradient for $E$ on ${\mathfrak M}^\ast:={\mathfrak M}\setminus \{\Phi\ ;\ DE(\Phi)=0\}$. We choose a cut-off for the action of the  Pseudo-gradient above the energy levels
${W}_{k-1}+\delta_{k-1}/2$ in order for the flow to preserving the membership in $\mbox{Sweep}_{N_k}$. Following the classical {\it Palais deformation} arguments (see for instance \cite{Riv-columb}) we deduce the existence of $\Phi_k$ such that $E(\Phi_k)=W_k$ and $DE(\Phi_k)=0$.

Assuming now $\mathfrak M$ defines in fact an {\it Hilbert manifold} modeled on an Hilbert space ${\mathfrak H}$ and that $D^2E$ is everywhere {\it Fredholm}. Denote
\[
\p \,\mbox{Sweep}_{N_k}:=\lf\{
\begin{array}{l}
\ds (Z,\Psi)\in {\mathcal P}_{N_{k-1}}\times\mbox{Lip}(\p(B^{n_k}\times Z),{\mathfrak M})\quad;\quad\exists\, (Y,\Phi)\in \mbox{Sweep}_{N_k}\\[3mm]
\ds \p Y=\p(B^{n_k}\times Z)\quad;\quad\Phi=\Psi\quad\mbox{ on }\p Y
\end{array}
\rg\}
\]
By definition we have
\[
\mbox{Sweep}_{N_k}\subset\lf\{ C\in {\mathfrak C}_{N_k}\quad\exists (Z,\Psi)\in \p \,\mbox{Sweep}_{N_k}\quad \p C=\Psi_\ast[\p (B^{n_k}\times Z)]\rg\}
\]
where ${\mathfrak C}_{N_k}$ is the space of $N_k-$polyhedral chains in ${\mathfrak M}$. Observe that $[C]\in H_{N_k}({\mathfrak M}, \p (B^{n_k}\times Z))$ is non trivial. Indeed, assume there exists $D\subset B$ such that $\p D=\p C$ we would have $W_{k}<W_{k-1}+2^{-1} \delta_{k-1}$ which contradicts (\ref{strict-sigma}) . Hence $\mbox{Sweep}_{N_k}$ is by definition an {\it homological family} of dimension $N_k$ with boundary the cycles $\Psi_\ast[\p (B^{n_k}\times Z)]$ (see \cite{Gho}). Using corollary 10.5 of \cite{Gho} we obtain (\ref{indice-palais}).  This concludes
the proof of theorem~\ref{th-I.1}. \hfill $\Box$

\medskip

We observe that the topological condition regarding ${\Phi}$ on $\Om_{{\Phi}}$ is preserved by enlarging the set. Precisely the following lemma holds.

\begin{Lm}
\label{lm-I.4}
Let $V\subset Z\times\p B^{n_k}$ such that $\Om_{{\Phi}}\subset V$ and such that $\pi_{k-1}\circ{\Phi}$ extends continuously on $V$. Let $\om\in H^{n_{k-1}}(V)$ given by $\om:=j^\ast(\pi_{k-1}\circ{\Phi})^\ast\om_{k-1}$ where $j$ is the canonical inclusion map $j\, :\, \Om_{{\Phi}}\rightarrow V$. Assume
\[
[V\cap (\{x\}\times Z)]\in H_{N_{k-1}}(V,\p V,{\Z}_2)\quad\mbox{ is {\bf Poincar\'e dual} to } \om\in H^{n_k-1}(V,{\Z}_2)
\]
for some $x\in \p B^{n_k}$, then the condition iv) is satisfied.\hfill $\Box$
\end{Lm}
\noindent{\bf Proof of lemma~\ref{lm-I.4}.} The inclusion map $j$ induces a map on relative $N_{k-1}-$chains as follows
\[
j\,:\, C_{N_{k-1}}(V,\p V,{\Z}_2)\,\rightarrow\,C_{N_{k-1}}(\Om_{{\Phi}},\p\Om_{{\Phi}},{\Z}_2)\quad.
\] 
by restriction to $\Om_{\vec{\Phi}}$. We then have the restriction operator
\[
j_\ast\,:\, H_{N_{k-1}}(V,\p V,{\Z}_2)\,\rightarrow\,H_{N_{k-1}}(\Om_{{\Phi}},\p\Om_{{\Phi}},{\Z}_2)\quad.
\]
For any $\al\in H^{N_{k-1}}(\Om_{{\Phi}},\p\Om_{{\Phi}},{\Z})$ we have for the {\it cup product} $(\pi_{k-1}\circ{\Phi})^\ast\om_{k-1}\smile\al\in H^{N_{k-1}+n_k-1}(\Om_{{\Phi}},\p\Om_{{\Phi}},{\Z}_2)$
\[
\begin{array}{l}
\ds\lf<(\pi_{k-1}\circ{\Phi})^\ast\om_{k-1}\smile\al,\lf[\Om_\Phi\rg]\rg>=\lf<(\pi_{k-1}\circ{\Phi})^\ast\om_{k-1}\smile\al,j_\ast\lf[V\rg]\rg>\\[3mm]
\ds\quad\quad=\lf<j^\ast\lf((\pi_{k-1}\circ{\Phi})^\ast\om_{k-1}\smile\al\rg),\lf[V\rg]\rg>=\lf<\om\smile j^\ast\al,\lf[V\rg]\rg>=\lf<j^\ast\al,\om\frown\lf[V\rg]\rg>
\end{array}
\]
We are assuming that $\om$ is Dual to $[V\cap (\{x\}\times Z)]\in H_{N_{k-1}}(V,\p V,{\Z}_2)$ in other words
\[
[V\cap (\{x\}\times Z)]=\om\frown\lf[V\rg]
\]
Hence we have proved that $\forall\,\al\in H^{N_{k-1}}(\Om_{{\Phi}},\p\Om_{{\Phi}},{\Z})$
\[
\lf<(\pi_{k-1}\circ{\Phi})^\ast\om_{k-1}\smile\al,\lf[\Om_\Phi\rg]\rg>=\lf<j^\ast\al,[V\cap (\{x\}\times Z)]\rg>=\lf<\al,j_\ast[V\cap (\{x\}\times Z)]\rg>
=\lf<\al,[\Om_{{\Phi}}\cap (\{x\}\times Z)]\rg>
\]
Which implies
\[
[\Om_{{\Phi}}\cap (\{x\}\times Z)]=(\pi_{k-1}\circ{\Phi})^\ast\om_{k-1}\frown[\Om_{{\Phi}}]
\] 
Hence $[\Om_{{\Phi}}\cap (\{x\}\times Z)]$ is dual to $(\pi_{k-1}\circ{\Phi})^\ast\om_{k-1}$ and the lemma is proved.\hfill $\Box$

\begin{Rm}
\label{rm-I.3p}
For  $n_1=1$ one can afford to restrict to $(Y,\vec{\Phi})\in\mbox{Sweep}_1$ where $Y=(-1,+1)$ moreover one can replace $H^{n_2-1}({\mathcal C}_1,{\Z}_2)$ by $H^{n_2-1}({\mathcal C}_1,{\Z})$ (which is ``richer'') in the definition of $\mbox{Sweep}_{N^2}$. This is due to the fact that the chain $U$ in the proof of theorem~\ref{th-I.1} can be taken to be a segment homeomorphic to $(-1,+1)$ and that there is no orientation problem at this first level of the hierarchy in the case $n_1=1$. The whole proof
in this case for the passage from $k=1$ to $k=2$ is transposable word by word by replacing ${\Z}_2$ by ${\Z}$. 
%This fact will be used
%later in the proof we are giving of the Willmore conjecture.
\hfill$\Box$
\end{Rm}

 \renewcommand{\theequation}{A.\arabic{equation}}
\renewcommand{\theTh}{A.\arabic{Th}}
\renewcommand{\theProp}{A.\arabic{Prop}}
\renewcommand{\theLma}{A.\arabic{Lma}}
\renewcommand{\theCo}{A.\arabic{Co}}
\renewcommand{\theRm}{A.\arabic{Rm}}
\renewcommand{\theequation}{A.\arabic{equation}}
\setcounter{equation}{0} 
\reset
\appendix
\section{Appendix}
\subsection{The action of $O(4)$ on $S^2_+\times S^2_-$}
\begin{Lma}
\label{lm-iso-s2s2}
For any element ${\mathfrak I}\in O(S^2\times S^2)$ there exists ${\mathcal I}\in O({\R}^4)$ and $i,j\in\{1,2\}$ such that
\[
\forall\, \vec{a},\vec{b}\in S^3\quad\mbox{ s. t. }\quad\vec{a}\cdot\vec{b}=0\quad, \quad{\mathfrak I}\lf((\vec{a}\wedge\vec{b})^+,(\vec{a}\wedge\vec{b})^-\rg)=(-1)^i\lf(({\mathcal I}(\vec{a})\wedge{\mathcal I}(\vec{b}))^+,(-1)^j({\mathcal I}(\vec{a})\wedge{\mathcal I}(\vec{b}))^-\rg)
\]
\hfill $\Box$
\end{Lma}
\noindent{\bf Proof of lemma~\ref{lm-iso-s2s2}.}
Every element of $O(S^2\times S^2)$ is represented by one of the following matrix
\[
\lf(
\begin{array}{cc}
A& 0\\[3mm]
0& B
\end{array}
\rg)\quad\mbox{ where }\quad A,B \in O(3)\quad\mbox{ or }\quad\lf(
\begin{array}{cc}
0& C\\[3mm]
D& 0
\end{array}
\rg)\quad\mbox{ where }\quad C,D \in O(3)\quad.
\]
For any pair $(\vec{g},\vec{h})$ in $S^3\times S^3$ we associate the corresponding pair of quaternions $({\mathbf g},{\mathbf h})$ and we consider the associated element in $SO(4)$ given by
\[
{\mathcal R}_{({\mathbf g},{\mathbf h})} (\vec{ x}):= {\mathbf g}\, {\mathbf x}\, {\mathbf h}\quad,
\]
where we are mixing the ${\R}^4$-vector and the quaternion representatives. We are computing the corresponding action of ${\mathcal R}_{({\mathbf g},{\mathbf h})} \in SO(4)$ onto $S^2_+\times S^2_-$.
Let $\vec{a}$ and $\vec{b}$ be two unit vectors of ${\R}^4$ which are orthogonal to each other.  We are now computing
\[
{\mathfrak R}_{({\mathbf g},{\mathbf h})}\lf((\vec{a}\wedge\vec{b})^+,(\vec{a}\wedge\vec{b})^-\rg):=\lf(({\mathcal R}_{({\mathbf g},{\mathbf h})} (\vec{a})\wedge{\mathcal R}_{({\mathbf g},{\mathbf h})} (\vec{b}))^+,({\mathcal R}_{({\mathbf g},{\mathbf h})} (\vec{a})\wedge{\mathcal R}_{({\mathbf g},{\mathbf h})} (\vec{b}))^-\rg)
\]
First we consider ${\mathcal R}^+_{{\mathbf g}}:={\mathcal R}_{({\mathbf g},{\mathbf 1})}$. We have respectively 
\[
({\mathcal R}^+_{{\mathbf g}} (\vec{a})\wedge{\mathcal R}^+_{{\mathbf g}} (\vec{b}))^+=
\lf(
\begin{array}{c}
<\mathbf{b},{\mathbf g}^\ast{\mathbf i}{\mathbf g}\,{\mathbf a}>\\[3mm]
<\mathbf{b},{\mathbf g}^\ast{\mathbf j}{\mathbf g}\,{\mathbf a}>\\[3mm]
<\mathbf{b},{\mathbf g}^\ast{\mathbf k}{\mathbf g}\,{\mathbf a}>
\end{array}
\rg)\quad\mbox{ and }\quad
({\mathcal R}^+_{{\mathbf g}} (\vec{a})\wedge{\mathcal R}^+_{{\mathbf g}} (\vec{b}))^-=
\lf(
\begin{array}{c}
<\mathbf{b},{\mathbf a}{\mathbf i}>\\[3mm]
<\mathbf{b},{\mathbf a}{\mathbf j}>\\[3mm]
<\mathbf{b},{\mathbf a}{\mathbf k}>
\end{array}
\rg)
\]
This gives
\[
\lf(
\begin{array}{c}
<\mathbf{b},{\mathbf g}^\ast{\mathbf i}{\mathbf g}\,{\mathbf a}>\\[3mm]
<\mathbf{b},{\mathbf g}^\ast{\mathbf j}{\mathbf g}\,{\mathbf a}>\\[3mm]
<\mathbf{b},{\mathbf g}^\ast{\mathbf k}{\mathbf g}\,{\mathbf a}>
\end{array}
\rg)=
\lf(
\begin{array}{c}
<\mathbf{b},{\mathbf i}\,{\mathbf a}>\, <{\mathbf i},{\mathbf g}^\ast{\mathbf i}{\mathbf g}>+<\mathbf{b},{\mathbf j}\,{\mathbf a}>\, <{\mathbf j},{\mathbf g}^\ast{\mathbf i}{\mathbf g}>+<\mathbf{b},{\mathbf k}\,{\mathbf a}>\, <{\mathbf k},{\mathbf g}^\ast{\mathbf i}{\mathbf g}>\\[3mm]
<\mathbf{b},{\mathbf i}\,{\mathbf a}>\, <{\mathbf i},{\mathbf g}^\ast{\mathbf j}{\mathbf g}>+<\mathbf{b},{\mathbf j}\,{\mathbf a}>\, <{\mathbf j},{\mathbf g}^\ast{\mathbf j}{\mathbf g}>+<\mathbf{b},{\mathbf k}\,{\mathbf a}>\, <{\mathbf k},{\mathbf g}^\ast{\mathbf j}{\mathbf g}>\\[3mm]
<\mathbf{b},{\mathbf i}\,{\mathbf a}>\, <{\mathbf i},{\mathbf g}^\ast{\mathbf k}{\mathbf g}>+<\mathbf{b},{\mathbf j}\,{\mathbf a}>\, <{\mathbf j},{\mathbf g}^\ast{\mathbf k}{\mathbf g}>+<\mathbf{b},{\mathbf k}\,{\mathbf a}>\, <{\mathbf k},{\mathbf g}^\ast{\mathbf k}{\mathbf g}>
\end{array}
\rg)
\]
This gives
\[
{\mathfrak R}_{({\mathbf g},{\mathbf 1})}\lf((\vec{a}\wedge\vec{b})^+,(\vec{a}\wedge\vec{b})^-\rg):=\lf(  {\mathbf g}^\ast(\vec{a}\wedge\vec{b})^+{\mathbf g}\,,\,(\vec{a}\wedge\vec{b})^-\rg)
\]
Then we consider ${\mathcal R}^-_{{\mathbf h}}:={\mathcal R}_{({\mathbf 1},{\mathbf h})}$. We have respectively 
\[
({\mathcal R}^-_{{\mathbf h}} (\vec{a})\wedge{\mathcal R}^-_{{\mathbf h}} (\vec{b}))^+=
\lf(
\begin{array}{c}
<\mathbf{b},{\mathbf i}\,{\mathbf a}>\\[3mm]
<\mathbf{b},{\mathbf j}\,{\mathbf a}>\\[3mm]
<\mathbf{b},{\mathbf k}\,{\mathbf a}>
\end{array}
\rg)\quad\mbox{ and }\quad
({\mathcal R}^-_{{\mathbf h}} (\vec{a})\wedge{\mathcal R}^-_{{\mathbf h}} (\vec{b}))^-=
\lf(
\begin{array}{c}
<\mathbf{b},{\mathbf a}\,{\mathbf h}{\mathbf i}{\mathbf h}^\ast>\\[3mm]
<\mathbf{b},{\mathbf a}\,{\mathbf h}{\mathbf j}{\mathbf h}^\ast>\\[3mm]
<\mathbf{b},{\mathbf a}\,{\mathbf h}{\mathbf k}{\mathbf h}^\ast>
\end{array}
\rg)
\]
and this gives
\[
{\mathfrak R}_{({\mathbf 1},{\mathbf h})}\lf((\vec{a}\wedge\vec{b})^+,(\vec{a}\wedge\vec{b})^-\rg):=\lf(  (\vec{a}\wedge\vec{b})^+\,,\,{\mathbf h}\,(\vec{a}\wedge\vec{b})^-{\mathbf h}^\ast\rg)
\]
Thus the subgroup of the positive isometry group of $S^2\times S^2$ given by
\[
\lf(
\begin{array}{cc}
A& 0\\[3mm]
0& B
\end{array}
\rg)\quad\mbox{ where }\quad A,B \in SO(3)
\]
is fully generated by the isometries ${\mathfrak R}_{({\mathbf g},{\mathbf h})}$. Next, we observe that if we denote by ${\mathcal R}_\ast({\vec{a}}):={\mathbf a}^\ast$ and
\[
{\mathfrak R}_\ast\lf((\vec{a}\wedge\vec{b})^+,(\vec{a}\wedge\vec{b})^-\rg):=\lf(({\mathcal R}_\ast({\vec{a}})\wedge{\mathcal R}_\ast({\vec{b}}))^+,({\mathcal R}_\ast({\vec{a}})\wedge{\mathcal R}_\ast({\vec{b}}))^-\rg)
\]
We have in one hand
\[
({\mathcal R}_\ast({\vec{a}})\wedge{\mathcal R}_\ast({\vec{b}}))^+=\lf(
\begin{array}{c}
<{\mathbf b}^\ast,{\mathbf i}\,{\mathbf a}^\ast>\\[3mm]
<{\mathbf b}^\ast,{\mathbf j}\,{\mathbf a}^\ast>\\[3mm]
<{\mathbf b}^\ast,{\mathbf k}\,{\mathbf a}^\ast>\\[3mm]
\end{array}\rg)=-\lf(
\begin{array}{c}
<\mathbf{b},{\mathbf a}{\mathbf i}>\\[3mm]
<\mathbf{b},{\mathbf a}{\mathbf j}>\\[3mm]
<\mathbf{b},{\mathbf a}{\mathbf k}>
\end{array}
\rg)=-(\vec{a}\wedge\vec{b})^-
\]
and
\[
({\mathcal R}_\ast({\vec{a}})\wedge{\mathcal R}_\ast({\vec{b}}))^-=\lf(
\begin{array}{c}
<{\mathbf b}^\ast,{\mathbf a}^\ast\,{\mathbf i}>\\[3mm]
<{\mathbf b}^\ast,{\mathbf a}^\ast\,{\mathbf j}>\\[3mm]
<{\mathbf b}^\ast,{\mathbf a}^\ast\,{\mathbf k}>\\[3mm]
\end{array}\rg)=-\lf(
\begin{array}{c}
<\mathbf{b},{\mathbf i}{\mathbf a}>\\[3mm]
<\mathbf{b},{\mathbf j}{\mathbf a}>\\[3mm]
<\mathbf{b},{\mathbf k}{\mathbf a}>
\end{array}
\rg)=-(\vec{a}\wedge\vec{b})^+
\]
Hence
\[
{\mathcal R}_\ast=-\lf(
\begin{array}{cc}
0& I_2\\[3mm]
I_2& 0
\end{array}
\rg)
\]
Combining all the previous gives the lemma. \hfill $\Box$

\medskip

\begin{Lma}
\label{lma-w21} Let
\[
\pi\ :\ (z_1,z_2,z_3,z_4)\in S^3\ \longrightarrow\ \frac{(z_1,z_2,z_3)}{\sqrt{1-z_4^2}}\ \in S^2
\]
and for any $a=(a_1,a_2,a_3,0)\in B^3$
we denote
\[
\phi_a\ :\ z\in S^3\ \longrightarrow\ \phi_a(z):=(1-|a|^2)\frac{z-a}{|z-a|^2}-a\quad
\]
We have
\be
\label{A+1}
\sup_{a\in B^3}\int_{S^3}|d(\pi\circ\phi_a)|^2\ dvol_{S^3}=\int_{S^3}|d\pi|^2\ dvol_{S^3}=8\,\pi^2\quad.
\ee
\hfill $\Box$
\end{Lma}
\noindent{\bf  Proof of lemma~\ref{lma-w21}.} Formally one would like to use theorem 1.2 of \cite{ElS}, since $\pi$ is a weak harmonic map. However $\pi$ is not smooth (it is an example of singular harmonic morphism) and we cannot apply this result. We shall go instead through explicit computations. 

Recall that $\phi_a$ is a conformal transformation of $S^3$ but also a conformal transformation of $B^4$ hence $|d(\pi\circ\phi_a)|^2=|d\pi|^2(\phi_a)\ e^{2\la_{\phi_a}}$ where
$e^{\la_{\phi_a}}$ is the conformal factor given by
\[
e^{2\la_{{\phi}_a}}:=|\p_{z_i}\phi_a|^2=\frac{1-|a|^2}{|z-a|^4} \quad\quad\forall\ i=1\cdots 4
\]
We have $|d\pi|^2(z)={2}/{(1-|z|^2)}$. Hence
\[
\int_{S^3}|d(\pi\circ\phi_a)|^2\ dvol_{S^3}=2\,\int_{S^3} \frac{e^{2\la_{{\phi}_a}}}{1-|\phi^4_a|^2}\ dvol_{S^3}
\]
We  have
\[
e^{-\la_{{\phi}_a}}=\frac{|z-a|^2}{1-|a|^2}\quad\mbox{ and }\quad|\phi_a+a|=(1-|a|^2)\frac{z-a}{|z-a|^2}\quad\Rightarrow\quad e^{-\la_{{\phi}_a}}=\frac{1-|a|^2}{|\phi_a+a|^2}
\]
Combining the previous identities gives
\[
\int_{S^3}|d(\pi\circ\phi_a)|^2\ dvol_{S^3}=2\,\int_{S^3} \frac{e^{3\la_{{\phi}_a}}}{1-|\phi^4_a|^2}\ \frac{1-|a|^2}{|\phi_a+a|^2}\ dvol_{S^3}=2\,\int_{S^3} \frac{1-|a|^2}{(1-z_4^2)\, |z+a|^2}\ dvol_{S^3}
\]
Due to the obvious axial symmetry of the problem we can afford to restrict to the case $a=t\, \ep_1$. Hence in order to prove the lemma we have to establish that
\[
\begin{array}{l}
\ds t\in (-1,+1)\, \longrightarrow \int_{S^3}\frac{1-t^2}{(1-z_4^2)\,(1+t^2+2\,t\, z_1)}\ dvol_{S^3}\\[5mm]
\ds\quad\quad=\int_{-1}^{+1}ds\int_{z_4=s}\frac{1-t^2}{(1-z_4^2)\, |dz_4|_{S^3}\,(1+t^2+2\,t\, z_1)}\ d{\mathcal H}^2\res\{z_4=s\}
\end{array}
\]
achieves its maximum at $t=0$. We have $ |dz_4|_{S^3}=\sqrt{1-z_4^2}$. Introducing $\cosh \al=(1+t^2)/(1-t^2)$ and $\sinh\al=2\,t/(1-t^2)$, and restricting to $t\ge 0$ we obtain
\[
\begin{array}{l}
\ds\int_{S^3}|d(\pi\circ\phi_a)|^2\ dvol_{S^3}=2\, \int_{-1}^{+1}\frac{ds}{(1-s^2)^{3/2}}\ \frac{2\pi}{\cosh\al}\int_{-\sqrt{1-s^2}}^{\sqrt{1-s^2}}\frac{dz_1}{1+z_1\tanh \al}\\[5mm]
\ds\quad\quad=4\pi\, \int_{-1}^{+1}\frac{ds}{(1-s^2)^{3/2}}\ \frac{1}{\sinh\al} \log\lf[\frac{ 1+ \tanh\al\, \sqrt{1-s^2}}{ 1- \tanh\al\, \sqrt{1-s^2}}\rg]
\end{array}
\]
Let $A:=\tanh^{-1}\al$ and $\sigma:=\sqrt{1-s^2}$ we study now
\[
f\ :\ A\in [1,+\infty]\ \longrightarrow\ f(A):=\sqrt{A^2-1}\  \log\lf[\frac{ A+ \sigma}{ A- \sigma}\rg]\quad\mbox{ and }f(\infty)=2\,\sigma
\]
We have
\[
f'(A)=\frac{A}{\sqrt{A^2-1}}\ g(A) \quad\mbox{ where }\quad g(A):= \log\lf[\frac{A+\sigma}{A-\sigma}\rg]-2\,\sigma\, \lf(A-\frac{1}{A}\rg)\ \frac{1}{A^2-\sigma^2}
\]
Observe $g(1)=\log((1+\sigma)/(1-\sigma))>0$  and $g(\infty)=0$. A direct computation again gives
\[
g'(A)=-\frac{4\,\sigma\, A^2\ (1-\sigma^2)+2\,\sigma\, (A^2-\sigma^2)}{A^2\ (A^2-\sigma^2)}<0\quad\quad\quad\forall A\ge 1
\]
Hence $g\ge 0$ on $[1,+\infty)$. Thus $f$ is increasing on $[1,+\infty)$ and we obtain
\[
\forall t\in (-1,1)\quad\ \int_{S^3}\frac{1-t^2}{(1-z_4^2)\,(1+t^2+2\,t\, z_1)}\ dvol_{S^3}\le \int_{S^3}\frac{1}{1-z_4^2}\ dvol_{S^3}=\frac{1}{2}\int_{S^3}|d\pi|^2\ dvol_{S^3}
\]
This concludes the proof of lemma~\ref{lma-w21}.\hfill $\Box$

\medskip

We are using the following proposition which is fairly standard in {\it Palais deformation theory} but for which nevertheless we give a proof below.
\begin{Prop}
\label{pr-A.1}
Let $\{\mbox{Sweep}_{N_l}\}_{l\le k}$ be a minmax hierarchy for a $C^2$ Lagrangian $E$ on a Banach manifold ${\mathfrak M}$. Assume moreover that $E$ satisfies the {\it Palais Smale condition}. Then for any $\ep>0$ there exists
$\delta>0$ such that
\[
\forall\, (Y,{\Phi})\in \mbox{Sweep}_{N_k} \quad \max_{y\in Y} E({\Phi}(y))<W_k+\delta
\]
then there exists $y_k^\ep\in Y$ and a critical point ${\Phi}_k$ of $E$ in ${\mathfrak M}$ such that
\[
E({\Phi}_k)=W_k\quad\mbox{ and }\quad d_{\mathbf P}({\Phi}(y_k^\ep),{\Phi}_k)<\ep\quad.
\]
where $d_{\mathbf P}$ is the usual {\it Palais distance} issued from the {\it Finsler structure} induced by the {\it Banach Manifold structure} on ${\mathfrak M}$\hfill $\Box$
\end{Prop}
\noindent{\bf Proof of proposition~\ref{pr-A.1}.}  Let $\delta>0$ such that
\[
W_{k-1}+\delta<W_{k}\quad.
\]
Consider the pseudo-gradient $X_k$  defined in ${\mathfrak M}$ and locally lipschitz in ${\mathfrak M}\setminus {\mathfrak M}^\ast$ and multiply  by a cut-off function supported in $[W_{k}-\delta,W_{k}+\delta]$ and equal to one on
$[W_{k}-\delta/2,W_{k}+\delta/2]$ in such a way that we have 
\[
\forall \,{\Phi}\in {\mathfrak M} \quad \|X_k({\Phi})\|_{{\Phi}}\le \|DE({\Phi})\|_{{\Phi}}
\]
and
\[
E({\Phi})\in [W_{k}-\delta/2,W_{k}+\delta/2]\quad\Rightarrow\quad \lf<X_k({\Phi}),DE({\Phi})\rg>_{T_{{\Phi}}{\mathfrak M},T^\ast_{{\Phi}}{\mathfrak M}}>  \|DE({\Phi})\|_{{\Phi}}^2
\]
Since the pseudo-gradient is supported in the level sets larger than $W_{k-1}$ it's flow $\phi_t$ generates a family of homeomorphisms
preserving $\mbox{Sweep}_{N_k}$ due to the conditions iii) and v)  in the definition of a hierarchy. For any $\eta\in (0,\delta/2)$ and any $(Y,{\Phi})\in \mbox{Sweep}_{N_k}$ such that
\[
\max_{y\in Y}E({\Phi}(y))\le W_k+\delta/2
\]
we consider the images $(Y,\phi_t({\Phi}))\in \mbox{Sweep}_{N_k}$. Denote by $d_{\mathbf P}$ the {\it Palais distance} associated
to the Finsler structure $\|\cdot\|$ and for which $(\mathfrak M,d_{\mathbf P})$ is complete (see \cite{Riv-columb}). Following lecture 2 of \cite{Riv-columb} we have for all ${y\in Y}$ and  any $t_1<t_2<t_{max}^{y}$
\be
\label{A.I.1}
 d_{\mathbf P}\lf(\phi_{t_1}({\Phi}(y)),\phi_{t_2}({\Phi}(y))\rg)\le \ 2\ \sqrt{t_2-t_1}\ \lf[E(\phi_{t_1}(\vec{\Phi}(y)))-E(\phi_{t_2}(\vec{\Phi}(y)))\rg]^{1/2}
\ee
where $t_{max}^{y}$ is the maximal existence time of $\phi_t({\Phi}(y))$. For a given $y\in Y$, assuming $t_{max}^y<+\infty$, because of the
previous inequality $\phi_t({\Phi}(y))$ realizes a Cauchy sequence for $d_{\mathbf P}$. Since ${\mathfrak M}$ is complete for $d_{\mathbf P}$,
the only possibility is that $\lim_{t\rightarrow t_{max}^y} \phi_t({\Phi}(y))\in {\mathfrak M}^\ast$. 
Assume first that 
\[
\forall\, y\in Y\quad E({\Phi}(y))\ge W_{k}-\delta \quad\Rightarrow \quad \|DE({\Phi}(y))\|_{{\Phi}(y)}\ge\delta^{1/4}
\]
Let $T>0$ and $\al>0$ such that 
\[
\forall \, t< T\quad \forall\, y\in Y\quad E(\phi_t({\Phi}(y)))\ge W_{k}-\delta \quad\Rightarrow \quad \|DE(\phi_t({\Phi}(y)))\|_{\phi_t({\Phi}(y))}\ge\delta^{1/4}
\]
We have for any ${\Phi}\in E^{-1}([W_{k}-\delta/2,W_{k}+\delta/2])$
\[
-\,\lf.\frac{d E(\phi_t({\Phi}))}{dt}\rg|_{t=0}\ge \|DE({\Phi})\|_{{\Phi}}^2
\]
This implies
\[
\max_{y\in Y}E(\phi_T({\Phi}(y)))\le W_k+\delta/2-\sqrt{\delta}\ T
\]
Denote $T_{\max}$ the first time such that either
\[
\exists \,y\in Y\quad\mbox{ s.t. }\quad t_{max}^y=T_{\max}
\]
which implies $DE(\phi_{T_{max}}(\vec{\Phi}(y))=0$ and $E(\phi_{T_{max}}({\Phi}(y)))\ge W_{k}-\delta$ or
\[
\exists \ y\in Y\quad E(\phi_t({\Phi}(y)))\ge W_{k}-\delta\quad \mbox{ and }\quad \|DE(\phi_{T_{max}}({\Phi}(y)))\|_{\phi_{T_{max}}({\Phi}(y))}={\delta}^{1/4}
\]
We clearly have $\sqrt{\delta}\, T_{max}\le\delta/2$. This implies using (\ref{A.I.1})
\[
\max_{y\in Y}d_{\mathbf P}({\Phi}(y),\phi_{T_{max}}({\Phi}(y)))\le {\delta^{1/4}}\ \sqrt{W_{k}/2}
\]
%Using \cite{Riv-minmax} and \cite{Riv-reg} there is a sequence of positive integers $\delta_k>0$ such that
%\[
%\begin{array}{c}
%\ds\max_{s\in Y_k} \mbox{Area}(\vec{\Phi}(s))\le  W_k+2\,\delta_k\\[3mm]
%\ds \Downarrow\\[3mm]
%\exists \quad\Sigma_k\quad\mbox{minimal s.t. }\quad\mbox{Area}(\Sigma_k)=W_k\quad\mbox{ and }\quad\exists\, t\in Y_k\quad\mbox{s. t. }{\mathcal F}(\vec{\Phi}(t)-\Sigma_k)\le\delta_k
%\end{array}
%\]
Collecting the various cases and summarizing we have obtained the following :
\be
\label{A.I.2}
\begin{array}{c}
\ds\forall\ \delta>0\quad\mbox{ and }\forall\ (Y,{\Phi})\in  \mbox{Sweep}_{N_k} \quad \mbox{ if }\quad\max_{y\in Y} E({\Phi}(y))<W_k+\delta\\[3mm]
\ds\quad\mbox{ then }\exists\ y\in Y\quad\mbox{ and }{\Phi}_\delta\in {\mathfrak M}\quad\mbox{ s.t. }
 |E({\Phi}_\delta)-W_k|\le\delta\\[3mm]
 \quad\mbox{and }\lf\|DE({\Phi}_\delta)\rg\|\le \delta^{1/4}
 \mbox{ moreover }\quad d_{\mathbf P}({\Phi}_\delta,{\Phi}(y))\le {\delta^{1/4}}\ \sqrt{W_{k}/2}
\end{array}
\ee
We claim now that 
\be
\label{A.I.3}
\begin{array}{l}
\forall\ \ep>0\ \ \exists\,\delta>0\quad\mbox{ s.t. }\quad\forall \ {\Phi}_\delta\in {\mathfrak M}\quad\mbox{ s.t. }
 |E({\Phi}_\delta)-W_k|\le\delta\\[3mm]
 \quad\mbox{and }\lf\|DE({\Phi}_\delta)\rg\|\le \delta^{1/4}\quad\mbox{ then }\quad\exists\ {\Phi}_0\in{\mathfrak M} \\[3mm]
 \quad\mbox{ s. t. }\quad E({\Phi}_0)=W_k\quad,\quad DE({\Phi}_0)=0\quad\mbox{ and }\quad d_{\mathbf P}({\Phi}_\delta,{\Phi}_0)\le \ep
\end{array}
\ee
This fact is a direct consequence of the {\it Palais Smale Condition} .  Indeed if (\ref{A.I.3}) would be false there
would exist $\ep_0>0$, a sequence $\delta_i\rightarrow 0$ and a sequence ${\Phi}_{\delta_i}$ such that
\[
 |E(\vec{\Phi}_\delta)-W_k|\le\delta_i\quad\mbox{ and }\quad\lf\|DE({\Phi}_{\delta_i})\rg\|\le \delta_i^{1/4}
\]
but ${\Phi}_{\delta_i}$ would stay at a {\it Palais distance} larger than $\ep_0$ to any critical point to $E$ at the level set $W_k$ which contradicts {\it (P.S.)}.

\medskip

Combining (\ref{A.I.2}) and (\ref{A.I.3}) we obtain
\[
\begin{array}{c}
\ds\forall\ \ep>0\ \ \exists\,\delta>0\quad\mbox{ s.t. }\forall\ (Y,{\Phi})\in  \mbox{Sweep}_{N_k} \quad \mbox{ if }\quad\max_{y\in Y} E({\Phi}(y))<W_k+\delta\\[3mm]
\mbox{ then }\quad\exists\ \vec{\Phi}_0\in {\mathfrak M}
 \quad\mbox{ s. t. }\quad E({\Phi}_0)=W_k\quad,\quad DE({\Phi}_0)=0\quad\\[3mm]
 \mbox{ and }\quad\exists\ y\in Y\mbox{ s. t. }\quad d_{\mathbf P}({\Phi}(y),{\Phi}_0)\le \ep
\end{array}
\]
 This concludes the proof of proposition~\ref{pr-A.1}.\hfill $\Box$

\end{document}